\documentclass{amsart}

\usepackage{amsmath,amsthm,amssymb,mathrsfs,color,enumerate,accents,amsfonts,bbold,adjustbox,etoolbox,changepage,bbm,adjustbox,witharrows,microtype}
\usepackage{thmtools}
\usepackage[margin=1in]{geometry}

\usepackage{pgfplots}
\pgfplotsset{compat=1.15}
\usepgfplotslibrary{fillbetween} 
\usepackage{tikz-cd}
\usepackage{tikz}
\usetikzlibrary{tqft,knots,decorations.markings,arrows,external,hobby,positioning}

\usepackage{graphicx,xcolor,transparent}
\graphicspath{./Pictures/}

\usepackage{hyperref}
\hypersetup{
    colorlinks=true, 
    linktoc=all,     
    linkcolor=blue,  
}

\newtheorem{thm}{Theorem}[section]
\newtheorem*{thm*}{Theorem}
\newtheorem{cor}{Corollary}[section]
\newtheorem{lem}{Lemma}[section]
\newtheorem{prop}{Proposition}[section]

\theoremstyle{definition}
\newtheorem{dff}{Definition}[section]
\newtheorem{xmp}{Example}[section]
\newtheorem{rmk}{Remark}

\newcommand{\mbb}[1]{\mathbb{#1}}
\newcommand{\mcal}[1]{\mathcal{#1}}

\DeclareMathOperator{\Hom}{Hom}

\DeclareMathOperator{\Endo}{End}

\definecolor{darkgreen}{RGB}{26,190,26}

\newcommand{\rev}{\overset{\longrightarrow}{\mathrm{ev}}}
\newcommand{\rcoev}{\overset{\longrightarrow}{\mathrm{coev}}}
\newcommand{\lev}{\overset{\longleftarrow}{\mathrm{ev}}}
\newcommand{\lcoev}{\overset{\longleftarrow}{\mathrm{coev}}}

\newcommand{\Id}{\mathrm{Id}}
\renewcommand{\t}{\mathrm{t}}

\tikzset{
  mid arrow/.style={decoration={markings, mark=at position 0.5 with {\arrow{#1}}}, postaction={decorate}}
}

\tikzset{
  third arrow/.style={decoration={markings, mark=at position 0.75 with {\arrow{#1}}}, postaction={decorate}}
}

\newcommand{\ctr}{\t^{\mcal C}}
\newcommand{\ztr}{\t^{\mcal Z(\mcal C)}}

\title{Chromatic Spherical and Renormalized Reshetikhin-Turaev Invariants for Drinfeld Centers}

\author{Nicolas Bridges}
\address{Department of Mathematics and Statistics \\ 3900 Old Main Hill \\ Utah State University \\ Logan, UT \\ 84322}
\email{nicolas.a.bridges@gmail.com}

\author{Shawn Cui}
\address{Departments of Mathematics and Physics \\ 150 N University St \\ Purdue University \\ West Lafayette, IN \\ 47907}
\email{cui177@purdue.edu}

\author{Nathan Geer}
\address{Department of Mathematics and Statistics \\ 3900 Old Main Hill \\ Utah State University \\ Logan, UT \\ 84322}
\email{nathan.geer@gmail.com}

\date{\today}

\begin{document}

\begin{abstract}
    We give invariants of pairs $(M,\Gamma)$ consisting of a closed connected oriented three-manifold and an oriented framed ribbon graph embedded in $M$ colored by the Drinfeld center of a unimodular finite chromatic category $\mcal C$. This invariant generalizes the chromatic spherical and renormalized Reshetikhin-Turaev invariants of three-manifolds. We use Heegaard-Graph diagrams which represent the pair $(M,\Gamma)$ and use the data of a unimodular chromatic category and objects of the Drinfeld center to construct the invariant. We define a ribbon graph associated to embedded framed links to obtain an invariant of pairs $(M,L)$. We show that if $L$ is the empty link, then the invariant recovers the chromatic spherical invariant, and if $M$ is the three-sphere, then the invariant recovers the renormalized Reshetikhin-Turaev invariant from the Drinfeld center. The relationship is proven by recovering a chromatic structure for the center induced by that from $\mcal C$. The main theorem of the paper demonstrates that if the graph is a link $L$ colored by a generator in the Drinfeld center, then we recover the chromatic spherical invariant of the surgery manifold $M(L)$, contributing to a proof of the relationship between the chromatic spherical and renormalized Reshetikhin-Turaev invariants for unimodular chromatic categories, generalizing a recent result for unimodular Hopf algebras by Reina. To this end, we give a Heegaard diagram for $M(L)$ coming from the Heegaard-Link diagram representing the pair $(M,L)$.
\end{abstract}

\maketitle

\tableofcontents

\section{Introduction}
\label{Introduction_Sec}

The classification of topological quantum field theories and their invariants has been of particular interest recently. The long-standing conjecture bridging the Turaev-Viro and Reshetikhin-Turaev invariants was built piece by piece culminating recently in the complete connection established in \cite{turaev2017monoidal}. These invariants are both defined by and rely heavily on semisimplicity for their input categories. In light of this, there have been non-semisimple counterparts of the Turaev-Viro and Reshetikhin-Turaev invariants such as: the Kuperberg invariants from non-semisimple Hopf algebras defined on Heegaard diagrams of manifolds with framing \cite{kuperberg1997non}, the Hennings-Kauffman-Radford invariant defined from non-semisimple Hopf algebras on surgery links \cite{hen96} and its categorification by Lyubashenko \cite{lyubashenko1995invariants}, the chromatic spherical invariants defined using spherical finite tensor categories without any condition on a framing \cite{Geer2020kuperberg}, the renormalized Hennings/Reshetikhin-Turaev invariants from non-semisimple ribbon categories \cite{de2018renormalized,de20223dimtqft}, and many more.

While there are significant contributions to the non-semisimple counterparts to the Turaev-Viro and Reshetikhin-Turaev invariants, less is known about the relationship between them. The authors in \cite{muller2026lyubashenko} proved a relationship between string nets from pivotal tensor categories and the Lyubashenko invariant for the Drinfeld center, generalizing the semisimple Levin-Wen string net relationship. The author in \cite{reina2026chromatic} showed that for unimodular Hopf algebras, the chromatic spherical and renormalized Hennings invariants are related by the Drinfeld double. In the same article, it was conjectured that the chromatic spherical and renormalized Reshetikhin-Turaev invariants are similarly related by the Drinfeld center. In this article we provide an affirmative answer to this conjecture when the chromatic category is unimodular and finite.

In this paper, we prove a formula in Theorem \ref{thm:mtraceforcenter} for the non-degenerate modified trace in the center of a unimodular chromatic category $\mcal C$ via the left and right adjoint to the forgetful functor. This modified trace induces a chromatic morphism in the center (Theorem \ref{thm:chromaticforcenter}) induced from a chromatic morphism in $\mcal C$. The left and right adjoint can be parametrized using the central monad and comonad given by the coend and end, respectively, yielding formulae which can be computed for a given chromatic category. 

We define an invariant $Z_{\mathrm{ExCGPT}}$ of pairs $(M,\Gamma)$ where $M$ is a closed connected oriented 3-manifold and $\Gamma$ is a framed embedded $\mcal Z(\mcal C)$ ribbon graph in $M$. The invariant is computed using Heegaard diagrams. In view of relating the chromatic spherical and renormalized Reshetikhin-Turaev invariants, we define a ribbon graph $\Gamma$ from an embedded framed link $L$ in $M$ using the chromatic morphism in the center to obtain an invariant $Z_{\mathrm{ExCGPT}}$ of pairs $(M,L)$. Two features of this definition are Corollary \ref{cor:recoverCGPT} and Theorem \ref{thm:ExCGPT_is_RT}, which demonstrate that our invariant recovers both the chromatic spherical and renormalized Reshetikhin-Turaev invariants. Additionally, we use an algorithm to compute a Heegaard diagram for the manifold $M(L)$ obtained by surgery on $M$ along $L$ defined in \cite{bridges20263man} to prove $Z_{\mathrm{ExCGPT}}$ on the pair $(M,L)$ and the pair $(M(L),\emptyset)$ are equal. This implies Corollary \ref{cor:CGPT_equals_RTcenter} stating that the chromatic spherical invariant from $\mcal C$ equals the renormalized Reshetikhin-Turaev invariant from $\mcal Z(\mcal C)$.

\section{Unimodular Chromatic Categories and Graphical Calculus}
\label{ChromaticCats_Sec}

\subsection{\texorpdfstring{Finite tensor $\Bbbk$-categories}{Finite tensor k-categories}}

Let $\mcal C$ be a strict monoidal category with monoidal product $\otimes$ and monoidal unit $\mbb 1$. Then $\mcal C$ is called a \emph{pivotal} category if it has morphisms

\begin{center}
\begin{tabular}{lll}
    $\rcoev_V:\mbb 1 \to V\otimes V^*$, & & $\rev_V:V^*\otimes V \to \mbb 1 $\\
    $\lcoev_V:\mbb 1 \to V^*\otimes V$, & & $\lev_V:V\otimes V^* \to \mbb 1 $
\end{tabular}
\end{center}

satisfying the equations

\begin{align*}
    &(\Id_V\otimes \rev_V)(\rcoev_V \otimes \mathrm{Id_V}) = \Id_V\\
    &(\rev_V \otimes \Id_{V^*})(\Id_{V^*} \otimes \rcoev_V) = \Id_{V^*}\\
    &(\lev_V \otimes \Id_V)(\Id_V \otimes \lcoev_V) = \Id_V\\
    &(\Id_{V^*} \otimes \lev_V)(\lcoev_V \otimes \Id_{V^*}) = \Id_{V^*}
\end{align*}

For $\Bbbk$ a field, an (abelian) \emph{tensor $\Bbbk$-category} is a locally small pivotal abelian category $\mcal C$ such that 
\begin{itemize}
    \item its hom-sets are $\Bbbk$-vector spaces,
    \item composition of morphisms is $\Bbbk$-bilinear
    \item the monoidal product $\otimes$ of morphisms is $\Bbbk$-bilinear
    \item the algebra map $\Bbbk\to \mathrm{End}_{\mcal C}(\mbb 1)$ given by $k\mapsto k~\Id_{\mbb{1}}$ is an isomorphism
\end{itemize}

We say that a tensor $\Bbbk$-category $\mcal C$ is \emph{finite} (as in \cite{EGNO}) if its hom-sets are finite-dimensional, every object has finite length, it has enough projectives, and there are finitely many isomorphism classes of simple objects.

\subsection{Penrose graphical calculus}

For a strict finite tensor $\Bbbk$-category $\mcal C$, we will represent morphisms in $\mcal C$ using plane diagrams read from bottom to top, called \emph{$\mcal C$-colored ribbon graphs}. These diagrams are made of oriented arcs colored by objects of $\mcal C$ and of boxes colored by morphisms of $\mcal C$, called \emph{coupons}. Coupons have a distinguished bottom face indicating the domain of the morphism. The opposite face is called the top and indicates the codomain. The arcs connect the boxes only to the top or bottom faces and have no intersections or self-intersections. The identity $\Id_V$ of $V\in \mcal C$, a morphism $f:V\to U$, and the composition of two morphisms $f:V\to U$ and $g:U\to W$ are represented as follows:

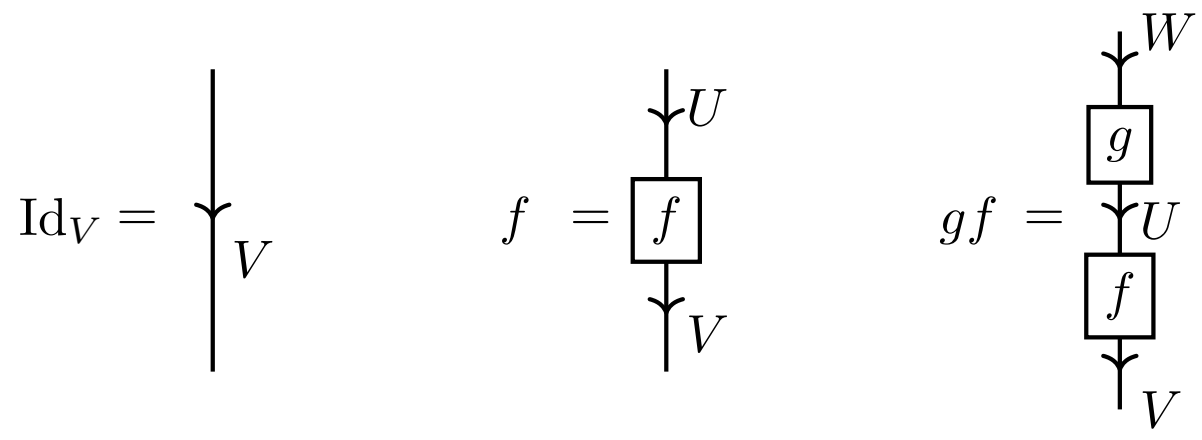
\begin{center}
    \begin{tikzpicture}
        \node at (0,0) {$\Id_V$};
        \node at (0.5,0) {$=$};
        \draw[thick,mid arrow = >] (1,1) -- (1,-1);
        \node[above right] at (1,-0.5) {$V$};

        \begin{scope}[shift = {(3,0)}]
            \node at (0,0) {$f$};
            \node at (0.5,0) {$=$};
            \draw[thick,mid arrow=>] (1,1) node[below right] {$U$} -- (1,0.25);
            \draw[thick,mid arrow=>] (1,-0.25) -- (1,-1) node[above right] {$V$};
            \node[draw,thick, fill=white, minimum height = 0.5cm] at (1,0) {$f$};
        \end{scope}

        \begin{scope}[shift = {(6,0)}]
            \node at (0,0) {$gf$};
            \node at (0.5,0) {$=$};
            \draw[thick,mid arrow = >] (1,1.25) node[right] {$W$} -- (1,0.75);
            \draw[thick,mid arrow = >] (1,0.25) -- node[right] {$U$} (1,-0.25);
            \draw[thick,mid arrow = >] (1,-0.75) -- (1,-1.25) node[right] {$V$};
            \node[draw,thick, fill=white, minimum height = 0.5cm] at (1,0.5) {$g$};
            \node[draw,thick,fill=white, minimum height = 0.5cm] at (1,-0.5) {$f$};
        \end{scope}
    \end{tikzpicture}
\end{center}

We represent the monoidal product of two morphisms $f:V\to U$ and $g:X\to Y$ by juxtaposition:

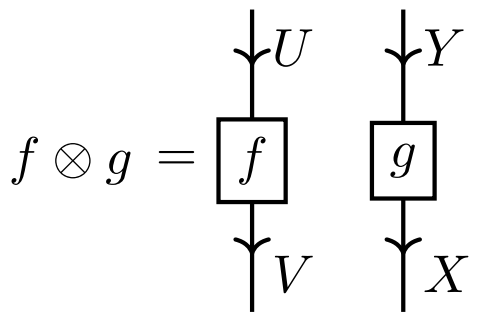
\begin{center}
    \begin{tikzpicture}
        \node at (-0.2,0) {$f\otimes g$};
            \node at (0.5,0) {$=$};
            \draw[thick,mid arrow=>] (1,1) node[below right] {$U$} -- (1,0.25);
            \draw[thick,mid arrow=>] (1,-0.25) -- (1,-1) node[above right] {$V$};
            \node[draw,thick, fill=white, minimum height = 0.5cm] at (1,0) {$f$};

            \draw[thick,mid arrow=>] (2,1) node[below right] {$Y$} -- (2,0.25);
            \draw[thick,mid arrow=>] (2,-0.25) -- (2,-1) node[above right] {$X$};
            \node[draw,thick, fill=white, minimum height = 0.5cm] at (2,0) {$g$};
    \end{tikzpicture}
\end{center}

Since $\mcal C$ is pivotal, we allow arcs directed upwards. Such an arc, colored by $V\in \mcal C$, represents $V^*$ in the source or target of the associated morphisms. For example, a morphism $f:V^*\otimes X \to U\otimes Y \otimes W^*$ may be depicted as

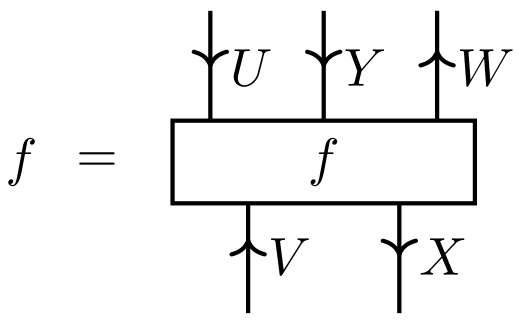
\begin{center}
    \begin{tikzpicture}
        \node at (0,0) {$f$};
        \node at (0.5,0) {$=$};
        \draw[thick,mid arrow=>] (1.25,1) -- node[right] {$U$} (1.25,0.25);
        \draw[thick,mid arrow=>] (2,1) -- node[right] {$Y$} (2,0.25);
        \draw[thick,mid arrow=<] (2.75,1) -- node[right] {$W$} (2.75,0.25);
        \draw[thick,mid arrow=<] (1.5,-0.25) -- node[right] {$V$} (1.5,-1);
        \draw[thick,mid arrow=>] (2.5,-0.25) -- node[right] {$X$} (2.5,-1);
        \node[draw,thick,fill=white,minimum height=0.5cm,minimum width = 2cm] at (2,0) {$f$};
    \end{tikzpicture}
\end{center}

We can depict the duality morphisms as

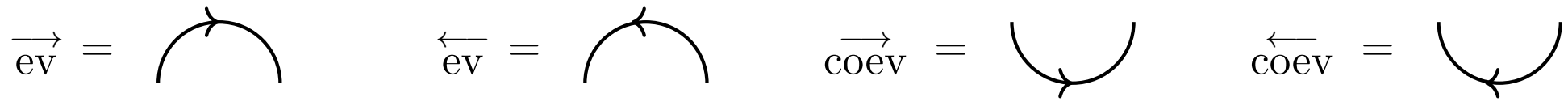
\begin{center}
    \begin{tikzpicture}
        \node at (0,0) {$\rev$};
        \node at (0.5,0) {$=$};
        \draw[thick,mid arrow=>] (1,-0.25) to[curve through = {(1.5,0.25)}] (2,-0.25);

        \begin{scope}[shift = {(3.5,0)}]
            \node at (0,0) {$\lev$};
        \node at (0.5,0) {$=$};
        \draw[thick,mid arrow=<] (1,-0.25) to[curve through = {(1.5,0.25)}] (2,-0.25);
        \end{scope}

        \begin{scope}[shift = {(7,0)}]
            \node at (-0.2,0) {$\rcoev$};
        \node at (0.5,0) {$=$};
        \draw[thick,mid arrow=>] (1,0.25) to[curve through = {(1.5,-0.25)}] (2,0.25);
        \end{scope}

        \begin{scope}[shift = {(10.5,0)}]
            \node at (-0.2,0) {$\lcoev$};
        \node at (0.5,0) {$=$};
        \draw[thick,mid arrow=<] (1,0.25) to[curve through = {(1.5,-0.25)}] (2,0.25);
        \end{scope}
    \end{tikzpicture}
\end{center}

The dual of a morphism $f:X\to Y$ can thus be depicted as

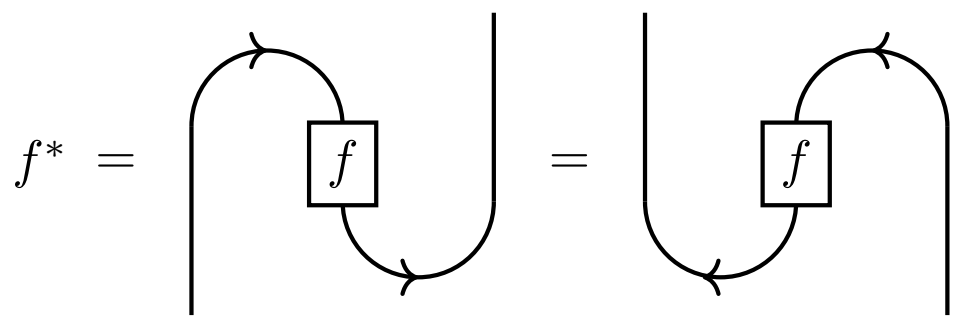
\begin{center}
    \begin{tikzpicture}
        \node at (0,0) {$f^*$};
        \node at (0.5,0) {$=$};
        \draw[thick] (1,-1) -- (1,0.25);
        \draw[thick,mid arrow=>] (1,0.25) to[curve through = {(1.5,0.75)}] (2,0.25);
        \draw[thick,mid arrow=>] (2,-0.25) to[curve through = {(2.5,-0.75)}] (3,-0.25);
        \draw[thick] (3,-0.25) -- (3,1);
        \node[draw,thick,fill=white,minimum height=0.5cm] at (2,0) {$f$};
        \node at (3.5,0) {$=$};
        \begin{scope}[shift = {(3,0)}]
            \draw[thick] (3,-1) -- (3,0.25);
        \draw[thick,mid arrow=>] (3,0.25) to[curve through = {(2.5,0.75)}] (2,0.25);
        \draw[thick,mid arrow=<] (1,-0.25) to[curve through = {(1.5,-0.75)}] (2,-0.25);
        \draw[thick] (1,-0.25) -- (1,1);
        \node[draw,thick,fill=white,minimum height=0.5cm] at (2,0) {$f$};
        \end{scope}
    \end{tikzpicture}
\end{center}

Since $\mcal C$ is pivotal, any morphism represented by a plane diagram is invariant under isotopies of the diagram fixing the top and bottom endpoints.

\subsection{Tensor ideals and m-traces}

Let $\mcal C$ be a finite tensor $\Bbbk$-category. An \emph{ideal} of $\mcal C$ is a full subcategory, $\mcal I$, of $\mcal C$ which is closed under the monoidal product and closed under retracts. I.e.,
\begin{itemize}
    \item \textbf{Closed under $\otimes$:} If $V\in \mcal I$ and $W\in \mcal C$, then $V\otimes W$ and $W\otimes V$ are objects of $\mcal I$.
    \item \textbf{Closed under retracts:} If $V\in \mcal I$ and $W\in \mcal C$ and there are morphisms $f:W\to V$, $g:V\to W$ such that $gf = \Id_W$, then $W\in \mcal I$.
\end{itemize}

An \emph{m-trace} (or \emph{modified trace}) on an ideal $\mcal I$ of $\mcal C$ is a family of linear maps 
\[\{\t_V:\Endo_{\mcal C}(V)\to \Bbbk\}_{V\in\mcal I}\]
such that the following conditions hold:
\begin{itemize}
    \item \textbf{Cyclicity:} If $U,V\in \mcal I$ then for any morphisms $f:V\to U$ and $g:U\to V$ in $\mcal C$, we have $\t_V(gf) = \t_U(fg)$. Graphically, this is
    
    
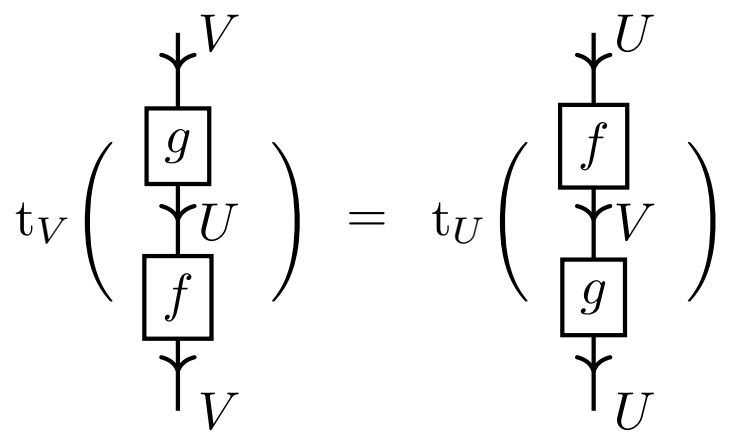
\begin{center}
        \begin{tikzpicture}
            \node at (0.25,0) {$\t_V \Bigg ($};
            \node at (1.75,0) {$\Bigg )$};
            \node at (2.25,0) {$=$};
            
            \draw[thick,mid arrow = >] (1,1.25) node[right] {$V$} -- (1,0.75);
            \draw[thick,mid arrow = >] (1,0.25) -- node[right] {$U$} (1,-0.25);
            \draw[thick,mid arrow = >] (1,-0.75) -- (1,-1.25) node[right] {$V$};
            
            \node[draw,thick, fill=white, minimum height = 0.5cm] at (1,0.5) {$g$};
            \node[draw,thick,fill=white, minimum height = 0.5cm] at (1,-0.5) {$f$};

            \begin{scope}[shift = {(2.75,0)}]
                \node at (0.25,0) {$\t_U \Bigg ($};
            \node at (1.75,0) {$\Bigg )$};
            
            \draw[thick,mid arrow = >] (1,1.25) node[right] {$U$} -- (1,0.75);
            \draw[thick,mid arrow = >] (1,0.25) -- node[right] {$V$} (1,-0.25);
            \draw[thick,mid arrow = >] (1,-0.75) -- (1,-1.25) node[right] {$U$};
            
            \node[draw,thick, fill=white, minimum height = 0.5cm] at (1,0.5) {$f$};
            \node[draw,thick,fill=white, minimum height = 0.5cm] at (1,-0.5) {$g$};
            \end{scope}
        \end{tikzpicture}
    \end{center}
    
    \item \textbf{Partial Trace:} If $U\in \mcal I$ and $W\in \mcal C$ then for any $f\in \Endo_{\mcal C}(U\otimes W)$ and $g\in \Endo_{\mcal C}(W\otimes U)$ we have
    \[\t_{U\otimes W}(f) = \t_U (\mathrm{ptr}_W(f)) =  \t_U\big ( (\Id_U\otimes \lev_W)(f\otimes \Id_{W^*})(\Id_U\otimes \rcoev_W)\big )\]
    \[\t_{W\otimes U}(g) = \t_U (\mathrm{ptr}_W(g)) = \t_U\big ( (\rev_W \otimes \Id_U)(\Id_{W^*} \otimes g)(\lcoev_W\otimes \Id_U) \big )\]
    Graphically, this is
    
    \begin{center}
        \begin{tikzpicture}
            \node at (0,0) {$\t_{U\otimes W} \Bigg ($};
            \node at (2.5,0) {$\Bigg )$};
            \node at (3,0) {$=$};

            \draw[thick,mid arrow=>] (1,1) node[right] {$U$} -- (1,0.25);
            \draw[thick,mid arrow=>] (2,1) node[right] {$W$} -- (2,0.25);
            \draw[thick,mid arrow=>] (1,-0.25) -- (1,-1) node[right] {$U$};
            \draw[thick,mid arrow=>] (2,-0.25) -- (2,-1) node[right] {$W$};
            
            \node[draw,thick,fill=white,minimum height=0.5cm,minimum width=1.5cm] at (1.5,0) {$f$};

            \begin{scope}[shift = {(3.5,0)}]
                \node at (0.25,0) {$\t_{U} \Bigg ($};
            \node at (3.5,0) {$\Bigg )$};

            \draw[thick,mid arrow=>] (1,1) node[right] {$U$} -- (1,0.25);
            \draw[thick,mid arrow=>] (2,1) node[right] {$W$} -- (2,0.25);
            \draw[thick,mid arrow=>] (1,-0.25) -- (1,-1) node[right] {$U$};
            \draw[thick,mid arrow=>] (2,-0.25) -- (2,-1) node[right] {$W$};
            \draw[thick] (2,-1) to[curve through = {(2.5,-1.5)}] (3,-1) -- (3,1) to[curve through ={(2.5,1.5)}] (2,1);
            
            \node[draw,thick,fill=white,minimum height=0.5cm,minimum width=1.5cm] at (1.5,0) {$f$};
            \end{scope}
        \end{tikzpicture}
        
        \begin{tikzpicture}
            \node at (0,0) {$\t_{W\otimes U} \Bigg ($};
            \node at (2.5,0) {$\Bigg )$};
            \node at (3,0) {$=$};

            \draw[thick,mid arrow=>] (1,1) node[right] {$W$} -- (1,0.25);
            \draw[thick,mid arrow=>] (2,1) node[right] {$U$} -- (2,0.25);
            \draw[thick,mid arrow=>] (1,-0.25) -- (1,-1) node[right] {$W$};
            \draw[thick,mid arrow=>] (2,-0.25) -- (2,-1) node[right] {$U$};
            
            \node[draw,thick,fill=white,minimum height=0.5cm,minimum width=1.5cm] at (1.5,0) {$g$};

            \begin{scope}[shift = {(4.5,0)}]
                \node at (-0.75,0) {$\t_{U} \Bigg ($};
            \node at (2.75,0) {$\Bigg )$};

            \draw[thick,mid arrow=>] (1,1) node[right] {$W$} -- (1,0.25);
            \draw[thick,mid arrow=>] (2,1) node[right] {$U$} -- (2,0.25);
            \draw[thick,mid arrow=>] (1,-0.25) -- (1,-1) node[right] {$W$};
            \draw[thick,mid arrow=>] (2,-0.25) -- (2,-1) node[right] {$U$};
            \draw[thick] (1,-1) to[curve through = {(0.5,-1.5)}] (0,-1) -- (0,1) to[curve through ={(0.5,1.5)}] (1,1);
            
            \node[draw,thick,fill=white,minimum height=0.5cm,minimum width=1.5cm] at (1.5,0) {$g$};
            \end{scope}
        \end{tikzpicture}   
    \end{center}
\end{itemize}

An m-trace on $\mcal I$ is \emph{non-degenerate} if for any $P\in \mcal I$, the following pairing is non-degenerate:
\[\Hom_{\mcal C}(\mbb 1,P) \times \Hom_{\mcal C}(P,\mbb 1) \to \Bbbk,\quad (x,y)\mapsto \t_P(xy)\]

If the above pairing is non-degenerate, then using the partial trace property and the pivotal structure of $\mcal C$ we see that for all $P\in \mcal I$ and all $V\in \mcal C$, the pairing 
\[\Hom_{\mcal C}(V,P) \times \Hom_{\mcal C}(P,V) \to \Bbbk,\quad (x,y)\mapsto \t_P(xy)\]
is non-degenerate as well.

Consider dual bases $\{x_i\}$ and $\{x^i\}$ for $\Hom_{\mcal C}(\mbb 1, P)$ and $\Hom_{\mcal C}(P,\mbb 1)$, respectively, dual under the non-degenerate m-trace $\t_P$. Then define $\Omega_P = \sum_i x_i\otimes x^i$ and $\Lambda_P = \sum_i x_i x^i \in \Endo_{\mcal C}(P)$. Then $\Omega_P$ and $\Lambda_P$ are independent of the choice of bases $\{x_i\}_i$ and $\{x^i\}_i$. 

\begin{lem}[\cite{Geer2023noncompact}]\label{lem:cuttingmorphismproperties}
    The morphisms $\Omega_P$ satisfy the following properties
    \begin{enumerate}
        \item \textbf{Duality:} If $\Omega_P = \sum_i x_i \otimes x^i$, then $\Omega_{P^*} = \sum_i (x^i)^* \otimes x_i^*$. Graphically, 
        
        
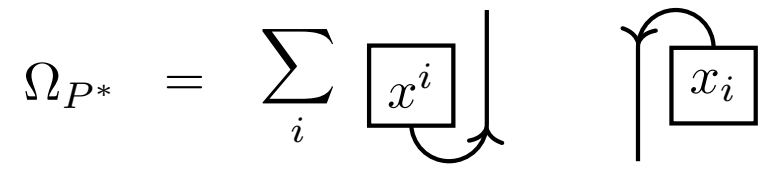
\begin{center}
            \begin{tikzpicture}
                \node at (0,0) {$\Omega_{P^*}$};
                \node at (0.75,0) {$=$};
                \node at (1.5,0) {$\displaystyle{\sum_i}$};

                \draw[thick,mid arrow=>] (2.25,-0.25) to[curve through = {(2.5,-0.5)}] (2.75,-0.25) -- (2.75,0.5);
                \draw[thick,mid arrow=>] (4.25,0.25) to[curve through = {(4,0.5)}] (3.75,0.25) -- (3.75,-0.5);
                
                \node[draw,thick,fill=white,minimum height = 0.5cm, minimum width = 0.5cm] at (2.25,0) {$x^i$};
                \node[draw,thick,fill=white,minimum height = 0.5cm, minimum width = 0.5cm] at (4.25,0) {$x_i$};
            \end{tikzpicture}
        \end{center}
        Thus $\Lambda_{P^*} = (\Lambda_P)^*$.
        
        \item \textbf{Naturality:} If $f:P\to Q$ is a morphism in $\mcal I$, then \[(f\otimes \Id_P)\circ\Omega_P = \Omega_Q\circ (\Id_Q \otimes f)\]
        \item \textbf{Rotation:} If $V\in \mcal C$ and $\Omega_{V\otimes P} = \sum_i x_i\otimes x^i$, then $\Omega_{P\otimes V} = \sum_i \widetilde{x}_i\otimes \widetilde{x}^i$ where
        \[\widetilde{x}_i = (\rev_V \otimes \Id_{P\otimes V}) (\Id_{V^*} \otimes x_i \otimes \Id_V) \lcoev_V\]
        and
        \[\widetilde{x}^i = \rev_V(\Id_{V^*} \otimes x^i \otimes \Id_V) (\lcoev_V \otimes \Id_{P\otimes V})\]
        Graphically,

        \begin{center}
            \begin{tikzpicture}
                \node at (0,0) {$\Omega_{P\otimes V}$};
                \node at (0.75,0) {$=$};
                \node at (1.5,0) {$\displaystyle{\sum_i}$};

                \draw[thick,mid arrow=>] (3.75,1)node[right] {$V$} -- (3.75,-0.25) to[curve through = {(3,-1)}] (2.25,-0.25) -- (2.25,0.25) to[curve through = {(2.5,0.5)}] (2.75,0.25);
                \draw[thick,mid arrow=>] (3.25,1)node[left] {$P$} -- (3.25,0.25);
                \draw[thick,mid arrow=>] (4.75,-0.25) to[curve through = {(4.5,-0.5)}] (4.25,-0.25) -- (4.25,0.25) to[curve through = {(5,1)}] (5.75,0.25) -- (5.75,-1)node[right] {$V$};
                \draw[thick,mid arrow=>] (5.25,-0.25) -- (5.25,-1)node[left] {$P$};

                \node[draw,thick,fill=white,minimum height = 0.5cm, minimum width = 1cm] at (3,0) {$x_i$};
                \node[draw,thick,fill=white,minimum height = 0.5cm, minimum width = 1cm] at (5,0) {$x^i$};
            \end{tikzpicture}
        \end{center}
    \end{enumerate}
\end{lem}

We call $\Lambda_P$ the \emph{cutting morphism} for $P$ in $\mcal C$.

\subsection{Projective covers and projective generators}

An object $P$ of $\mcal C$ is \emph{projective} if the functor $\Hom_{\mcal C}(P,\bullet):\mcal C\to \mathrm{Vect}_\Bbbk$ preserves epimorphisms. We say $P$ is \emph{injective} if the functor $\Hom_{\mcal C}(\bullet, P):\mcal C\to \mathrm{Vect}_\Bbbk$ preserves monomorphisms.

Denote by $\mathrm{Proj}(\mcal C)$ the full subcategory of projective objects.

\begin{prop}[\cite{geer2012traces}]
    $\mathrm{Proj}(\mcal C)$ is an ideal of $\mcal C$, and in $\mcal C$, projective and injective objects coincide.
\end{prop}

An object $G$ in an ideal $\mcal I$ is called an \emph{$\mcal I$-generator} if for each object $V$ of $\mcal C$, there is a natural number $n= n_V$ such that there is an epimorphism $G^{\oplus n}\twoheadrightarrow V$. When $\mcal I = \mathrm{Proj}(\mcal C)$, we call $G$ a \emph{projective generator}.

\begin{dff}
    Let $G$ be an $\mcal I$-generator in $\mcal C$. A \emph{chromatic morphism} for $G$ is a morphism $c:G\otimes G\to G\otimes G$ such that 
    \[ (\Id_G \otimes \lev_G \otimes \Id_G)(c\otimes\Lambda_{G^*\otimes G})(\Id_G \otimes \rcoev_G \otimes \Id_G) = \Id_{G\otimes G}\]
    or, graphically,
    \begin{center}
        \begin{tikzpicture}
            \draw[thick,mid arrow=>] (-0.25,1)node[right] {$G$} -- (-0.25,0.25);
            \draw[thick,mid arrow=>] (-0.25,-0.25) -- (-0.25,-1)node[right] {$G$};
            \draw[thick,mid arrow=>] (2.25,1)node[right] {$G$} -- (2.25,0.25);
            \draw[thick,mid arrow=>] (2.25,-0.25) -- (2.25,-1)node[right] {$G$};
            \draw[thick,mid arrow=<] (0.25,0.25) to[curve through={(1,0.75)}] (1.75,0.25);
            \draw[thick,mid arrow=>] (0.25,-0.25) to[curve through={(1,-0.75)}] (1.75,-0.25);
            \node[below] at (1,-0.75) {$G$};
            \node[above] at (1,0.75) {$G$};

            \node[draw,thick,fill=white,minimum height=0.5cm,minimum width = 1cm] at (2,0) {$\Lambda_{G^*\otimes G}$};
            \node[draw,thick,fill=white,minimum height=0.5cm,minimum width=1cm] at (0,0) {$c$};

            \node at (3.25,0) {$=$};

            \draw[thick,mid arrow=>] (4,1) -- (4,-1)node[right] {$G$};
            \draw[thick,mid arrow=>] (5,1) -- (5,-1)node[right] {$G$};
        \end{tikzpicture}
    \end{center}
\end{dff}

\begin{lem}[\cite{Geer2023noncompact}]
    Fix a projective generator $G$. For any object $X$ in $\mcal C$, and any projective object $P\in \mathrm{Proj}(\mcal C)$, there exists a morphism $c_P\in \Endo_{\mcal C}(P\otimes G)$ satisfying
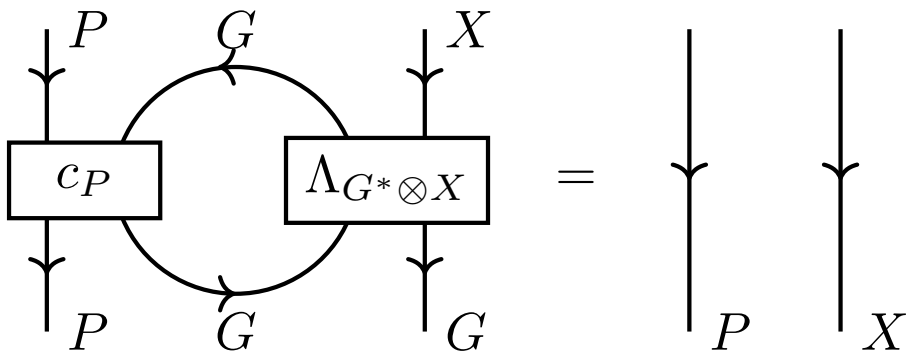
    \begin{center}
        \begin{tikzpicture}
            \draw[thick,mid arrow=>] (-0.25,1)node[right] {$P$} -- (-0.25,0.25);
            \draw[thick,mid arrow=>] (-0.25,-0.25) -- (-0.25,-1)node[right] {$P$};
            \draw[thick,mid arrow=>] (2.25,1)node[right] {$X$} -- (2.25,0.25);
            \draw[thick,mid arrow=>] (2.25,-0.25) -- (2.25,-1)node[right] {$G$};
            \draw[thick,mid arrow=<] (0.25,0.25) to[curve through={(1,0.75)}] (1.75,0.25);
            \draw[thick,mid arrow=>] (0.25,-0.25) to[curve through={(1,-0.75)}] (1.75,-0.25);
            \node[below] at (1,-0.75) {$G$};
            \node[above] at (1,0.75) {$G$};

            \node[draw,thick,fill=white,minimum height=0.5cm,minimum width = 1cm] at (2,0) {$\Lambda_{G^*\otimes X}$};
            \node[draw,thick,fill=white,minimum height=0.5cm,minimum width=1cm] at (0,0) {$c_P$};

            \node at (3.25,0) {$=$};

            \draw[thick,mid arrow=>] (4,1) -- (4,-1)node[right] {$P$};
            \draw[thick,mid arrow=>] (5,1) -- (5,-1)node[right] {$X$};
        \end{tikzpicture}
    \end{center}
    We call $c_P$ a \emph{chromatic morphism based on $P$} for the generator $G$.
\end{lem}

Clearly, a chromatic morphism based on $G$ is an ordinary chromatic morphism.

\begin{dff}
    A \emph{(finite) chromatic category} over $\Bbbk$ is a finite tensor $\Bbbk$-category $\mcal C$ endowed with a non-degenerate m-trace $\t$ on $\mathrm{Proj}(\mcal C)$ such that
    \begin{itemize}
        \item any non-zero morphism to the unit object $\mbb 1$ is an epimorphism,
        \item there exists a chromatic map for a non-zero projective generator.
    \end{itemize}
\end{dff}

We assume all chromatic categories are finite. Additionally, we will assume all chromatic categories are \emph{unimodular}, i.e. the socle of the projective cover $P_{\mbb 1}$ of the unit $\mbb 1$ is isomorphic to the unit.

\begin{xmp}
    Let $\mcal C$ be an abelian finite spherical tensor $\Bbbk$-category. Then $\mcal C$ is a chromatic category.

    Examples of these categories are the representation categories of unimodular unibalanced Hopf algebras (c.f. \cite{beliakova2021modified}).
\end{xmp}

\section{The Drinfeld Center}
\label{Center_Sec}

In this section, we review the center construction for a tensor $\Bbbk$-category $\mcal C$ and give a chromatic structure for the center inherited from $\mcal C$.

Fix a strict tensor $\Bbbk$-category $\mcal C$. Then define a new category $\mcal Z(\mcal C)$, called the \emph{Drinfeld center} of $\mcal C$, as having objects pairs $(X,\sigma)$ where $X$ is an object of $\mcal C$ and $\sigma = \{\sigma_Y:X\otimes Y \to Y\otimes X\}_{Y\in \mcal C}$ is a natural family of isomorphisms natural in $Y$ such that the following diagram commutes:

\begin{center}
    \begin{tikzcd}
     & Z\otimes X \otimes Y \arrow[rd, "\Id_Z\otimes \sigma_Y"] & \\
        X\otimes Z \otimes Y \arrow[ru, "\sigma_Z\otimes \Id_Y"] \arrow[rr,"\sigma_{Z\otimes Y}"'] & & Z\otimes Y\otimes X
    \end{tikzcd}
\end{center} 

Morphisms $f:(X,\sigma) \to (X',\sigma')$ in $\mcal Z(\mcal C)$ are morphisms $f:X \to X'$ in $\mcal C$ such that for each $Y\in \mcal C$, $(\Id_Y \otimes f)\circ \sigma_Y = \sigma'_Y\circ (f\otimes \Id_Y)$.

The category $\mcal Z(\mcal C)$ is a monoidal category with monoidal product \[(X,\sigma)\otimes (X',\sigma') := (X\otimes X', \widetilde{\sigma})\] where $\widetilde{\sigma}_Y:(X\otimes X')\otimes Y \to Y\otimes (X\otimes X')$ for $Y\in \mcal C$ is defined by the following commuting diagram:

\begin{center}
    \begin{tikzcd}
        (X\otimes X')\otimes Y \arrow[rr, "\widetilde{\sigma}"] \arrow[rd, "\Id_X\otimes \sigma'_Y"'] &  & Y\otimes (X\otimes X') \\ 
        & X\otimes Y\otimes X' \arrow[ru, "\sigma_Y\otimes \Id_{X'}"'] &
    \end{tikzcd}
\end{center}

The unit of $\mcal Z(\mcal C)$ is $(\mbb 1, \Id)$.

Note that when $\mcal C$ is strict, so is $\mcal Z(\mcal C)$. When $\mcal C$ is not strict, there are additional constraints coming from the associator and unitors in $\mcal C$.

Define the forgetful functor $\mcal U:\mcal Z(\mcal C)\to \mcal C$ by $\mcal U(X,\sigma) = X$ and $\mcal U(f) = f$. Then $\mcal U$ is a strict monoidal functor.

In addition to monoidal, the Drinfeld center $\mcal Z(\mcal C)$ is a factorizable, unimodular, braided finite tensor $\Bbbk$-category \cite{EGNO}. Since $\mcal C$ is unimodular and its m-trace is two-sided, it is spherical \cite{shibata2023modified}, and hence $\mcal Z(\mcal C)$ is a ribbon category induced from the spherical structure of $\mcal C$ \cite{shimizu2023ribbon}. In \cite{Geer2023chromatic} the authors prove that every unimodular finite tensor $\Bbbk$-category has chromatic morphisms. Later we explicitly construct a chromatic morphism for the Drinfeld center of a chromatic category.

\subsection{The non-degenerate m-trace in the center}

Since $\mcal C$ is a unimodular tensor $\Bbbk$-category, the forgetful functor $U:\mcal Z(\mcal C)\to \mcal C$ is a Frobenius functor. Hence it has a left and right adjoint $F:\mcal C\to \mcal Z(\mcal C)$. Denote by $\eta^L$ and $\epsilon^L$ the unit and counit for the adjunction $F\dashv U$ and denote by $\eta^R$ and $\epsilon^R$ the unit and counit for the adjunction $U\dashv F$.

\begin{lem}
    The functors $U$ and $F$ preserve projective objects.
\end{lem}

\begin{proof}
    This follows because $U$ and $F$ both have a left and right adjoint and thus each has an exact right adjoint.
\end{proof}


\begin{lem}\label{lem:genforcenter}
    If $G$ is a projective generator for $\mcal C$, then $F(G)$ is a projective generator for $\mcal Z(\mcal C)$.
\end{lem}

\begin{proof}
    Let $P$ be a projective object in $\mcal Z(\mcal C)$. Now, consider $U(P)$. Then since $U$ preserves projectives and $G$ is a projective generator, then there is a split epimorphism $G^{\oplus n}\twoheadrightarrow U(P)$. Applying $F$ (and using that it is an exact linear functor), we have a split epimorphism $F(G)^{\oplus n}\twoheadrightarrow FU(P)$. Since $U$ is faithful, $\epsilon^L_{P}:FU(P)\to P$ is a split epimorphism. Hence taking compositions yields a split epimorphism $F(G)^{\oplus n}\twoheadrightarrow P$, proving the claim.
\end{proof}

\begin{thm}\label{thm:mtraceforcenter}
    Let $\ctr = \{\ctr_P:\Endo(P)\to \Bbbk\}_{P\in \mathrm{Proj}(\mcal C)}$ be the unique non-degenerate m-trace on $\mathrm{Proj}(\mcal C)$. Let $F:\mcal C\to \mcal Z(\mcal C)$ be the left and right adjoint to the forgetful functor $U:\mcal Z(\mcal C)\to \mcal C$. Denote by $\eta^L:\Id_{\mcal C}\to UF$ the unit for the left adjunction $F\dashv U$. 
    
    Then there exists a right unit $\eta':\Id_{\mcal Z(\mcal C)}\to FU$  and right counit $\epsilon':UF\to \Id_{\mcal C}$ such that for each projective $P\in \mcal C$ and each $f\in \Endo_{\mcal Z(\mcal C)}(F(P))$, the function 

    \[ \ztr_{F(P)}(f) = \ctr_P(\epsilon'_P \circ U(f)\circ \eta^L_P)\]

    extends uniquely to the unique non-degenerate m-trace on $\mathrm{Proj}(\mcal Z(\mcal C))$ such that the adjoint mate $\Lambda^{\mcal Z}_{F(P)}$ of $\Lambda_P:P\to P$ takes the value $\ctr_P(\Lambda_P)$.
\end{thm}

\begin{proof}
    Denote by $\eta^R$ and $\epsilon^R$ a choice of right unit and right counit for the adjunction $U\dashv F$. For each $P\in \mathrm{Proj}(\mcal C)$, define a function $T_{F(P)}:\Endo_{\mcal Z(\mcal C)}(F(P))\to \Bbbk$ by \[T_{F(P)}(f) = \ctr_P(\epsilon^R_P\circ U(f)\circ \eta^L_P)\]

    Denote the induced pairing $\Hom_{\mcal Z(\mcal C)}(M , F(P)) \times \Hom_{\mcal Z(\mcal C)}(F(P), M)\to \Bbbk$ by \[\langle f,g \rangle^T = T_{F(P)}(fg)\]

    Notice that the map \[L_{P,M}:\Hom_{\mcal Z(\mcal C)}(F(P),M) \to \Hom_{\mcal C}(P,U(M))\] given by $L_{P,M}(f) = U(f)\eta^L_P$ is an isomorphism. Similarly, the map \[R_{P,M}:\Hom_{\mcal Z(\mcal C)}(M,F(P))\to \Hom_{\mcal C}(U(M),P)\] given by $R_{P,M}(f) = \epsilon^R_P U(f)$ is an isomorphism.

    Then for $f:M\to F(P)$ and $g:F(P)\to M$, \[\langle f,g\rangle^T = \ctr_P(R_{P,M}(f)L_{P,M}(g))\] Since $\ctr$ is a non-degenerate m-trace, its induced pairing is non-degenerate, implying the pairing $\langle f,g\rangle^T$ is non-degenerate.

    Note that $\mcal Z(\mcal C)$ has a non-degenerate m-trace on the ideal of projectives which we denote by $s = \{s_M:\Endo_{\mcal Z(\mcal C)}(M)\to \Bbbk\}_{M\in \mathrm{Proj}(\mcal Z(\mcal C))}$. By non-degeneracy of $s$, for any $f\in \Endo_{\mcal Z(\mcal C)}(F(P))$, there is a unique endomorphism $z_P:F(P)\to F(P)$ such that $T_{F(P)}(f) = s_{F(P)}(z_P\circ f)$.

    \begin{lem}
        The endomorphisms $z_P$ are invertible.
    \end{lem}

    \begin{proof}
        By non-degeneracy of the pairing $\langle f,g\rangle^T$, the pairing defined by $\langle f,g\rangle^s = s_{F(P)}(z_P\circ fg)$ is also non-degenerate. This implies $z_P$ is invertible. Indeed, suppose $z_P f = 0$, then $\langle f,g\rangle^T = 0$ for all $g$. This implies $f = 0$ by non-degeneracy, so left composition by $z_P$ is injective in the finite-dimensional algebra $\Endo_{\mcal Z(\mcal C)}(F(P))$, and hence is bijective.

        Now, consider the unique morphism $f:F(P)\to F(P)$ such that $z_P f = \Id_{F(P)}$. Then $z_P(fz_P - \Id_{F(P)}) = 0$ implying $fz_P = \Id_{F(P)}$, making $f$ the inverse of $z_P$ and hence $z_P$ is invertible.
    \end{proof}

     \begin{lem}
         For each $f:P\to Q$ in $\mathrm{Proj}(\mcal C)$ and each $g:F(P)\to F(Q)$, we have \[T_{F(P)}(g\circ F(f)) = T_{F(Q)}(F(f)\circ g)\]
     \end{lem}

     \begin{proof}
         We compute
         \begin{align*}
             T_{F(P)}(g\circ F(f)) &= \ctr_{P}(\epsilon^R_P\circ U(g)\circ UF(f)\circ \eta^L_P)\\
             &= \ctr_{P}(\epsilon^R_P\circ U(g)\circ \eta^L_Q\circ f)\\
             &= \ctr_{Q}(f\circ \epsilon^R_P\circ U(g)\circ \eta^L_Q)\\
             &= \ctr_Q(\epsilon^R_Q\circ UF(f)\circ U(g)\circ \eta^L_Q)\\
             &= T_{F(Q)}(F(f)\circ g)
         \end{align*}
     \end{proof}

     \begin{lem}
         For each $f:P\to Q$ in $\mathrm{Proj}(\mcal C)$, $z_Q\circ F(f) = F(f)\circ z_P$. 
     \end{lem}

     \begin{proof}
         Let $g:F(Q)\to F(P)$. Then
         \begin{align*}
             T_{F(Q)}(F(f)\circ g) &= s_{F(Q)}(z_Q\circ F(f)\circ g)\\
             T_{F(Q)}(F(f)\circ g) &= T_{F(P)}(g\circ F(f))\\
             T_{F(P)}(g\circ F(f)) &= s_{F(P)}(z_P\circ g\circ F(f))\\
             &= s_{F(Q)}(F(f)\circ z_P\circ g)
         \end{align*}

         Then by non-degeneracy of $s$, we must have that $z_Q\circ F(f) = F(f)\circ z_P$.
     \end{proof}

    Now, define $\epsilon'_P = \epsilon^R_P U(z_P^{-1})$ and $\eta'_M = z_{UM} \eta^R_M$

    \begin{lem}
        The maps $\epsilon'$ and $\eta'$ uniquely extend to a unit and counit for the right adjunction $U\dashv F$.
    \end{lem}

    \begin{proof}
        Since $\mcal C$ has enough projectives and $F$ is exact, the natural family $z_P:F(P)\to F(P)$ extends uniquely from $\mathrm{Proj}(\mcal C)$ to all of $\mcal C$. We now show that the snake identities hold for our proposed unit and counit.

        \begin{align*}
            \epsilon'_{U(M)} U(\eta'_M) &= \epsilon^R_{U(M)} U(z_{U(M)}^{-1}) U(z_{U(M)}\eta^R_M)\\
            &= \epsilon^R_{U(M)}U(\eta^R_M)\\
            &= \Id_{U(M)}
        \end{align*}

        \begin{align*}
            F(\epsilon'_P)\eta'_{F(P)} &= F(\epsilon^R_P U(z_P^{-1}))z_{UF(P)} \eta^R_{F(P)}\\
            &= z_P F(\epsilon^R_P)FU(z_P^{-1})\eta^R_{F(P)}\\
            &= z_P F(\epsilon^R_P) \eta^R_{F(P)}z_P^{-1}\\
            &= z_P z_P^{-1}\\
            &= \Id_{F(P)}
        \end{align*}
    \end{proof}

    Then for $f:F(P)\to F(P)$, the function $\ztr_{F(P)}(f) = \ctr_{P}(\epsilon'_P \circ U(f)\circ \eta^L_P)$ satisfies $\ztr_{F(P)}(f) = s_{F(P)}(f)$.
    
    Define the \emph{right projection isomorphism} $\theta_{P,M}:F(P\otimes U(M)) \to F(P)\otimes M$ by \[\theta_{P,M} = \epsilon^L_{F(P)\otimes M}F(\eta^L_P\otimes \Id_{U(M)})\]

    This isomorphism satisfies \[U(\theta_{P,M})\circ \eta^L_{P\otimes U(M)} = \eta^L_P\otimes \Id_{U(M)}\] and is invertible.

    \begin{lem}
        The right projection isomorphism $\theta_{P,M}$ satisfies \[\epsilon'_{P\otimes U(M)} = (\epsilon'_P \otimes \Id_{U(M)})U(\theta_{P,M})\] and it has inverse \[\theta_{P,M}^{-1} = F(\epsilon'_P\otimes \Id_{U(M)})\eta'_{F(P)\otimes M}\]
    \end{lem}

    \begin{proof}
        Let $f\in \Endo_{\mcal Z(\mcal C)}(F(P\otimes U(M)))$. Denote by $\widetilde{f} = \theta_{P,M} \circ f \circ \theta_{P,M}^{-1}$
        \begin{align*}
            \ctr_{P\otimes U(M)}((\epsilon'_P \otimes \Id_{U(M)})U(\theta_{P,M})U(f) \eta^L_{P\otimes U(M)})&= \ctr_{P\otimes U(M)}((\epsilon'_P \otimes \Id_{U(M)})U(\widetilde{f}) (\eta^L_P\otimes \Id_{U(M)}))\\
            &= \ctr_P(\epsilon'_P \mathrm{ptr}_{U(M)}(U(\widetilde{f})) \eta^L_P)\\
            &= s_{F(P)}(\mathrm{ptr}_{M}(\widetilde{f}))\\
            &= s_{F(P) \otimes M}(\widetilde{f})\\
            &= s_{F(P\otimes U(M))}(f)\\
            &= \ctr_{P\otimes U(M)}(\epsilon'_{P\otimes U(M)}\circ U(f)\circ \eta^L_{P\otimes U(M)})
        \end{align*}

        The left-hand side is the pairing of the morphism $(\epsilon'_P\otimes \Id_{U(M)})U(\theta_{P,M})$ with the morphism $L_{P\otimes U(M),F(P\otimes U(M))}(f)$. The right-hand side is the pairing of the morphism $\epsilon'_{P\otimes U(M)}$ with $L_{P\otimes U(M),F(P\otimes U(M))}(f)$. Non-degeneracy of $\ctr$ implies \[(\epsilon'_P\otimes \Id_{U(M)})U(\theta_{P,M}) = \epsilon'_{P\otimes U(M)}\]

        To prove the inverse formula, we compute the (right) adjoint mates of each side. Note first the adjoint mate of the right-hand side is $\epsilon'_P\otimes \Id_{U(M)}$. Then using the above computation,

        \begin{align*}
            \epsilon'_{P\otimes U(M)}U(\theta_{P,M}^{-1}) &= (\epsilon'_P\otimes \Id_{U(M)})U(\theta_{P,M})U(\theta_{P,M}^{-1})\\
            &= \epsilon'_P\otimes \Id_{U(M)}
        \end{align*}
    \end{proof}

    To extend to all of $\mathrm{Proj}(\mcal Z(\mcal C))$, we use that $F(G)$ is a projective generator. Then if $\iota:M \hookrightarrow F(G)^{\oplus n}$ and $\pi:F(G)^{\oplus n}\twoheadrightarrow M$ is a splitting of a projective object $M\in \mcal Z(\mcal C)$, then for an endomorphism $f:M\to M$, we define \[\ztr_M(f) = \ztr_{F(G)^{\oplus n}}(\iota \circ f\circ\pi)\] Then \[\ztr_M(f) = s_{F(G)^{\oplus n}}(\iota\circ f\circ \pi) = s_M(f)\]
     
    For the cutting morphism $\Lambda^{\mcal Z}_{F(P)}$, let $\Lambda_P = \sum_i x_i x^i\in \Endo_{\mcal C}(P)$ be the cutting morphism for $P$ defined in Section \ref{ChromaticCats_Sec}. Then consider the adjoint mate $x_i^{\mcal Z}$ of $x_i$ in $\Hom_{\mcal Z(\mcal C)}(\mbb 1, F(P))$ and the adjoint mate $x^{i,\mcal Z}$ of $x^i$ in $\Hom_{\mcal Z(\mcal C)}(F(P),\mbb 1)$. Then \[x_i^{\mcal Z} = F(x_i)\eta'_{\mbb 1},\qquad x^{i,\mcal Z} = \epsilon^L_{\mbb 1}F(x^i)\] Then define \[\Lambda^{\mcal Z}_{F(P)} := \sum_i x_i^{\mcal Z} x^{i,\mcal Z}\]

    Computing (while suppressing the summation) yields \begin{align*}
        t^{\mcal Z(\mcal C)}_{F(P)}(\Lambda^{\mcal Z}_{F(P)}) &= t^{\mcal C}_P(\epsilon'_P\circ U(\Lambda^{\mcal Z}_{F(P)})\circ \eta^L_P)\\
        &= t^{\mcal C}_P(\epsilon'_P\circ UF(x_i)\circ U(\eta'_{\mbb 1})\circ U(\epsilon^L_{\mbb 1})\circ UF(x^i)\circ \eta^L_P)\\
        &= t^{\mcal C}_P(x_i \circ \epsilon'_{U(\mbb 1)}\circ U(\eta'_{\mbb 1})\circ U(\epsilon^L_{\mbb 1})\circ \eta^L_{U(\mbb 1)}\circ x^i)\\
        &= t^{\mcal C}_P(x_i x^i)\\
        &= t^{\mcal C}_P(\Lambda_P)
    \end{align*}
\end{proof}

\subsection{A chromatic morphism for the center}

In the previous subsection, we showed that if $\Lambda_P = \sum_i x_i x^i$ is the cutting morphism for $P$ in $\mcal C$, then \[\Lambda^{\mcal Z}_{F(P)} = \sum_i x_i^{\mcal Z} x^{i,\mcal Z} = \sum_i F(x_i)\eta'_{\mbb 1} \circ \epsilon^L_{\mbb 1}F(x^i)\] is the cutting morphism for $F(P)$ in $\mcal Z(\mcal C)$.

To compute $\Lambda^{\mcal Z}_{F(G)^*\otimes F(G)}$, we first consider $\Lambda_{G^*\otimes UF(G)} = \sum_i x_i x^i$. Then consider the following chain of isomorphisms:

\[\Hom_{\mcal C}(G^*\otimes UF(G),\mbb 1) \to \Hom_{\mcal C}(UF(G),G) \to \Hom_{\mcal Z(\mcal C)}(F(G),F(G)) \to \Hom_{\mcal Z(\mcal C)}(F(G)^* \otimes F(G),\mbb 1)\] sending $f$ to 
\[A(f) = \rev_{F(G)}\left (\Id_{F(G)^*}\otimes F(\Id_{G}\otimes f)F(\rcoev_G\otimes \Id_{UF(G)})\eta'_{F(G)}\right )\]

and the chain of isomorphisms

\[\Hom_{\mcal C}(\mbb 1,G^*\otimes UF(G)) \to \Hom_{\mcal C}(G,UF(G))\to \Hom_{\mcal Z(\mcal C)}(F(G),F(G)) \to \Hom_{\mcal Z(\mcal C)}(\mbb 1, F(G)^*\otimes F(G))\]
sending $f$ to 
\[B(f) = \left (\Id_{F(G)^*}\otimes \epsilon^L_{F(G)}F(\lev_G\otimes \Id_{UF(G)})F(\Id_{G}\otimes f)\right )\lcoev_{F(G)}\]

Then define $\Lambda^{\mcal Z}_{F(G)^*\otimes F(G)}$ by \[\Lambda^{\mcal Z}_{F(G)^*\otimes F(G)} = \sum_i B(x_i)A(x^i)\]

Graphically,

\begin{center}
    \begin{tikzpicture}
        \node at (0,0) {$\Lambda^{\mcal Z}_{F(G)^*\otimes F(G)} = \displaystyle\sum_i$};
        \begin{scope}[shift = {(5,0)}]

        \begin{scope}[shift = {(0,1.5)}]
            \draw[mid arrow=>,thick] (0.25,1.75) -- (0.25,1.25);
            \draw[mid arrow=>,thick] (-0.25,1.25) .. controls +(0,0.5) and +(0,0.5) .. (-0.75,1.25) -- (-0.75,0.5);

            \draw[mid arrow=>,thick] (-0.25,2.5) -- (-0.25,1.875);

            \draw[mid arrow=>,thick] (-0.25,0.125) .. controls +(0,-2) and +(0,-2) .. (-3,0.125) -- (-3,4);

            \draw[mid arrow=>,thick] (-0.25,4) -- (-0.25,3.5);

            \node[draw, thick, fill = white, minimum height = 1cm, minimum width = 1cm] at (-0.25,3) {$\epsilon^L_{F(G)}$};

            \node[draw, thick, fill = white, minimum height = 0.5cm, minimum width = 1cm] at (0,1) {$x_i$};

            \node[draw, thick, minimum height = 1.75cm, minimum width = 1cm] at (-0.25,1) {$F\Bigg ( \qquad\qquad\qquad \Bigg )$};
        \end{scope}

        \begin{scope}[shift = {(0,-1.5)}]
            \draw[mid arrow=>,thick] (0.25,-1.25) -- (0.25,-1.75);
            \draw[mid arrow=<,thick] (-0.25,-1.25) .. controls +(0,-0.5) and +(0,-0.5) .. (-0.75,-1.25) -- (-0.75,-0.5);

            \draw[mid arrow=>,thick] (-0.25,-1.875) -- (-0.25,-2.5);

            \draw[mid arrow=<,thick] (-0.25,-0.125) .. controls +(0,2) and +(0,2) .. (-3,-0.125) -- (-3,-4);

            \draw[mid arrow=>,thick] (-0.25,-3.5) -- (-0.25,-4);

            \node[draw, thick, fill = white, minimum height = 1cm, minimum width = 1cm] at (-0.25,-3) {$\eta'_{F(G)}$};

            \node[draw, thick, fill = white, minimum height = 0.5cm, minimum width = 1cm] at (0,-1) {$x^i$};

            \node[draw, thick, minimum height = 1.75cm, minimum width = 1cm] at (-0.25,-1) {$F\Bigg ( \qquad\qquad\qquad \Bigg )$};

        \end{scope}
        \end{scope}
    \end{tikzpicture}
\end{center}

\begin{lem}
    $\Lambda^{\mcal Z}_{F(G)^*\otimes F(G)}$ is the cutting morphism for $F(G)^*\otimes F(G)$
\end{lem}

\begin{proof}
    We compute:

    \begin{align*}
        \ztr_{F(G)^*\otimes F(G)}(B(x_i)A(x^j)) &= \ztr_{F(G)}(\mathrm{ptr}_{F(G)^*}(B(x_i)A(x^j))\\
        &= \ztr_{F(G)}\Big (\epsilon^L_{F(G)}F(\lev_G\otimes \Id_{UF(G)})F(\Id_{G}\otimes x_i)\\
        &\qquad \circ F(\Id_{G}\otimes x^j)F(\rcoev_G \otimes \Id_{UF(G)})\eta'_{F(G)}\Big )\\
        &= \ztr_{F(G)}\Big (F(\Id_{G}\otimes x^j)F(\rcoev_G \otimes \Id_{UF(G)})\eta'_{F(G)}\\
        &\qquad \circ \epsilon^L_{F(G)}F(\lev_G\otimes \Id_{UF(G)})F(\Id_{G}\otimes x_i)\Big )\\
        &= \ctr_{G}\Big (\epsilon'_G UF(\Id_{G}\otimes x^j)UF(\rcoev_G \otimes \Id_{UF(G)})U(\eta'_{F(G)})\\
        &\qquad \circ U(\epsilon^L_{F(G)})UF(\lev_G\otimes \Id_{UF(G)})UF(\Id_{G}\otimes x_i)\eta^L_G\Big)\\
        &= \ctr_{G}\Big ((\Id_{G}\otimes x^j)(\rcoev_G \otimes \Id_{UF(G)})\epsilon'_{UF(G)}U(\eta'_{F(G)})\\
        &\qquad \circ U(\epsilon^L_{F(G)})\eta^L_{UF(G)}(\lev_G\otimes \Id_{UF(G)})(\Id_{G}\otimes x_i)\Big)\\
        &= \ctr_{G}\Big ((\Id_{G}\otimes x^j)(\rcoev_G \otimes \Id_{UF(G)})\\
        &\qquad \circ (\lev_G\otimes \Id_{UF(G)})(\Id_{G}\otimes x_i)\Big)\\
        &= \ctr_{UF(G)} \Big ( \mathrm{ptr}_{G^*}(x_i x^j)\Big )\\
        &= \ctr_{G^*\otimes UF(G)}(x_i x^j)\\
        &= \delta_{i,j}
    \end{align*}
\end{proof}

Recall that a chromatic morphism based on $F(G)$ in $\mcal Z(\mcal C)$ is a morphism $c^\mcal Z:F(G)\otimes F(G)\to F(G)\otimes F(G)$ such that 

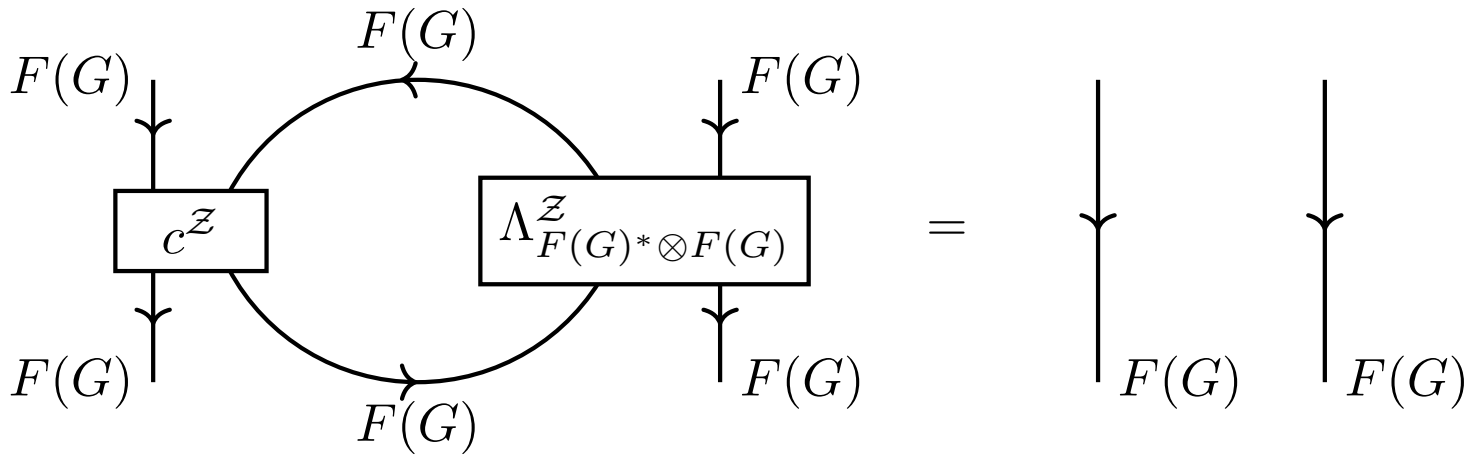
\begin{center}
        \begin{tikzpicture}
            \draw[thick,mid arrow=>] (-0.25,1)node[left] {$F(G)$} -- (-0.25,0.25);
            \draw[thick,mid arrow=>] (-0.25,-0.25) -- (-0.25,-1)node[left] {$F(G)$};
            \draw[thick,mid arrow=>] (3.5,1)node[right] {$F(G)$} -- (3.5,0.25);
            \draw[thick,mid arrow=>] (3.5,-0.25) -- (3.5,-1)node[right] {$F(G)$};
            \draw[thick,mid arrow=<] (0.25,0.25) to[curve through={(1.5,1)}] (2.75,0.25);
            \draw[thick,mid arrow=>] (0.25,-0.25) to[curve through={(1.5,-1)}] (2.75,-0.25);
            \node[below] at (1.5,-1) {$F(G)$};
            \node[above] at (1.5,1) {$F(G)$};

            \node[draw,thick,fill=white,minimum height=0.5cm,minimum width = 1cm] at (3,0) {$\Lambda^{\mcal Z}_{F(G)^*\otimes F(G)}$};
            \node[draw,thick,fill=white,minimum height=0.5cm,minimum width=1cm] at (0,0) {$c^{\mcal Z}$};

            \node at (5,0) {$=$};

            \draw[thick,mid arrow=>] (6,1) -- (6,-1)node[right] {$F(G)$};
            \draw[thick,mid arrow=>] (7.5,1) -- (7.5,-1)node[right] {$F(G)$};
        \end{tikzpicture}
    \end{center}

\begin{thm}\label{thm:chromaticforcenter}
    Let $c$ be a chromatic map for $G$. Then, a chromatic map for $F(G)$ is given by 
    \[ c^{\mcal Z} = \theta_{G,F(G)} F(\Id_G\otimes \eta^L_G)F(c)F(\Id_G\otimes \epsilon'_G)\theta_{G,F(G)}^{-1}\]
\end{thm}

\begin{proof}
    First, we note that $\theta_{P,M}$ is natural in the sense that for any $f:P\to Q$ and $g:M\to N$, we have \[(F(f)\otimes g)\theta_{P,M} = \theta_{Q,N}F(f\otimes U(g))\]
    We compute. Note that the left-hand side of the identity defining the chromatic morphism reduces graphically to:

    \begin{center}
        \begin{tikzpicture}
            \begin{scope}[shift = {(4,1.5)}]
            \draw[mid arrow=>,thick] (0.25,1.75) -- (0.25,1.25);
            \draw[mid arrow=>,thick] (-0.25,1.25) .. controls +(0,0.5) and +(0,0.5) .. (-0.75,1.25) -- (-0.75,0.5);

            \draw[mid arrow=>,thick] (-0.25,2.5) -- (-0.25,1.875);

            \draw[mid arrow=>,thick] (-0.25,0.125) .. controls +(-1,-5) and +(1,5) .. (-5,2);

            \draw[mid arrow=>,thick] (-0.25,4) -- (-0.25,3.5);

            \node[draw, thick, fill = white, minimum height = 1cm, minimum width = 1cm] at (-0.25,3) {$\epsilon^L_{F(G)}$};

            \node[draw, thick, fill = white, minimum height = 0.5cm, minimum width = 1cm] at (0,1) {$x_i$};

            \node[draw, thick, minimum height = 1.75cm, minimum width = 1cm] at (-0.25,1) {$F\Bigg ( \qquad\qquad\qquad \Bigg )$};
        \end{scope}

        \node[draw, thick, fill = white, minimum width = 4cm, minimum height = 1cm]  at (-2,0) {$F(c)$};

        \node[draw, thick, fill = white, minimum width = 4cm, minimum height = 1cm]  at (-2,1.5) {$F(\Id_G \otimes \eta^L_G)$};
        \node[draw, thick, fill = white, minimum width = 4cm, minimum height = 1cm]  at (-2,-1.5) {$F(\Id_G \otimes \epsilon'_G)$};

        \node[draw, thick, fill = white, minimum width = 4cm, minimum height = 1cm]  at (-2,3) {$\theta_{G,F(G)}$};
        \node[draw, thick, fill = white, minimum width = 4cm, minimum height = 1cm]  at (-2,-3) {$\theta_{G,F(G)}^{-1}$};

        \draw[mid arrow=>,thick] (-2,2.5) -- (-2,2);
        \draw[mid arrow=>,thick] (-2,1) -- (-2,0.5);
        \draw[mid arrow=>,thick] (-2,-0.5) -- (-2,-1);
        \draw[mid arrow=>,thick] (-2,-2) -- (-2,-2.5);
        \draw[mid arrow=>,thick] (-3,5.5) -- (-3,3.5);
        \draw[mid arrow=>,thick] (-3,-3.5) -- (-3,-5.5);

        \begin{scope}[shift = {(4,-1.5)}]
            \draw[mid arrow=>,thick] (0.25,-1.25) -- (0.25,-1.75);
            \draw[mid arrow=<,thick] (-0.25,-1.25) .. controls +(0,-0.5) and +(0,-0.5) .. (-0.75,-1.25) -- (-0.75,-0.5);

            \draw[mid arrow=>,thick] (-0.25,-1.875) -- (-0.25,-2.5);

            \draw[mid arrow=<,thick] (-0.25,-0.125) .. controls +(-1,5) and +(1,-5) .. (-5,-2);

            \draw[mid arrow=>,thick] (-0.25,-3.5) -- (-0.25,-4);

            \node[draw, thick, fill = white, minimum height = 1cm, minimum width = 1cm] at (-0.25,-3) {$\eta'_{F(G)}$};

            \node[draw, thick, fill = white, minimum height = 0.5cm, minimum width = 1cm] at (0,-1) {$x^i$};

            \node[draw, thick, minimum height = 1.75cm, minimum width = 1cm] at (-0.25,-1) {$F\Bigg ( \qquad\qquad\qquad \Bigg )$};

        \end{scope}
        \end{tikzpicture}
    \end{center}

    Let $a_i = (\lev_G \otimes \Id_{UF(G)})(\Id_G\otimes x_i)$ and $b_i = (\Id_{G}\otimes x^i)(\rcoev_G\otimes \Id_{UF(G)})$. Then using the naturality of $\theta_{G,F(G)}$, we obtain

        
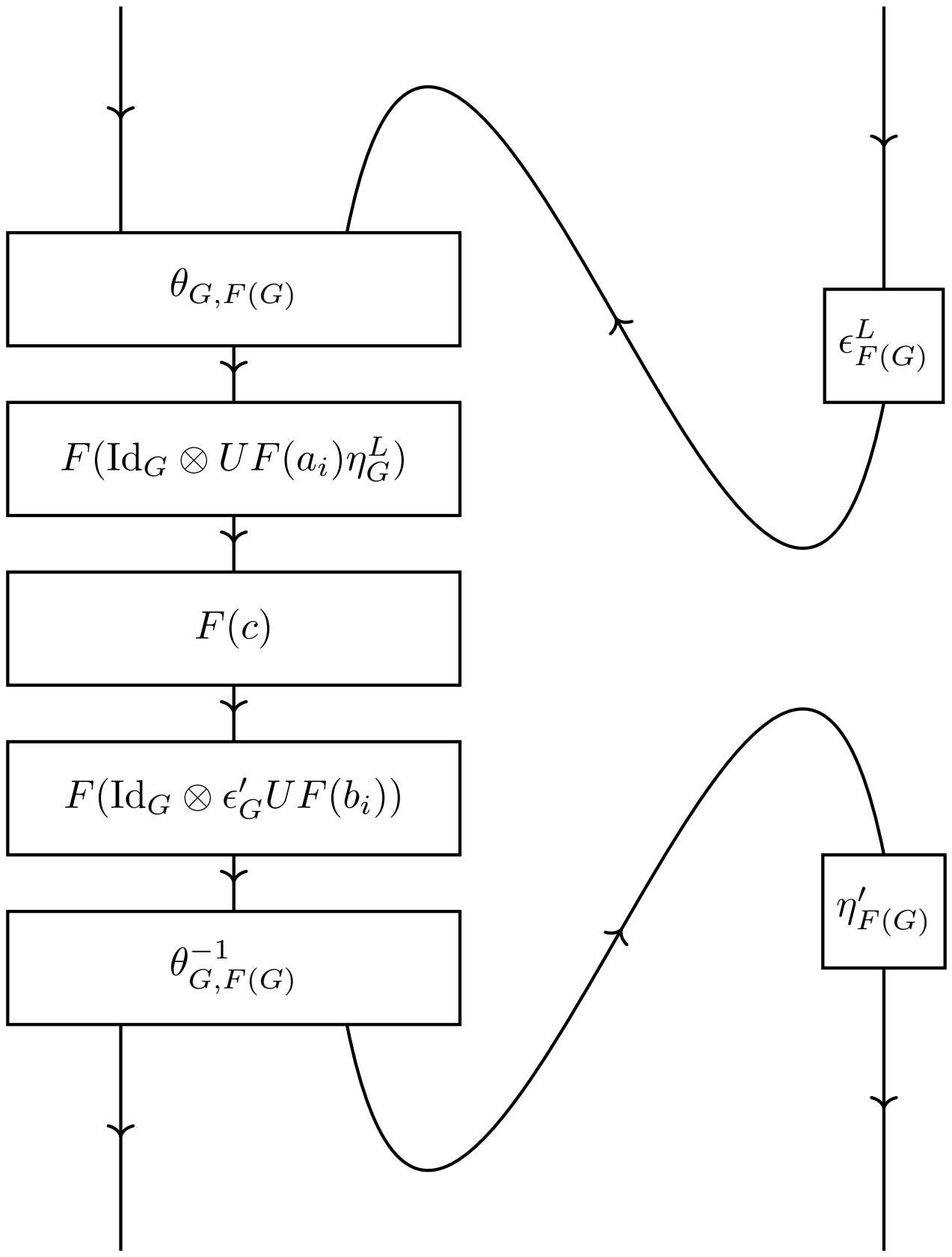
\begin{center}
        \begin{tikzpicture}
            \begin{scope}[shift = {(4,1.5)}]

            \draw[mid arrow=>,thick] (-0.25,0.5) .. controls +(-1,-5) and +(1,5) .. (-5,2);

            \draw[mid arrow=>,thick] (-0.25,4) -- (-0.25,1.5);

            \node[draw, thick, fill = white, minimum height = 1cm, minimum width = 1cm] at (-0.25,1) {$\epsilon^L_{F(G)}$};
        \end{scope}

        \node[draw, thick, fill = white, minimum width = 4cm, minimum height = 1cm]  at (-2,0) {$F(c)$};

        \node[draw, thick, fill = white, minimum width = 4cm, minimum height = 1cm]  at (-2,1.5) {$F(\Id_G \otimes UF(a_i)\eta^L_G)$};
        \node[draw, thick, fill = white, minimum width = 4cm, minimum height = 1cm]  at (-2,-1.5) {$F(\Id_G \otimes \epsilon'_G UF(b_i))$};

        \node[draw, thick, fill = white, minimum width = 4cm, minimum height = 1cm]  at (-2,3) {$\theta_{G,F(G)}$};
        \node[draw, thick, fill = white, minimum width = 4cm, minimum height = 1cm]  at (-2,-3) {$\theta_{G,F(G)}^{-1}$};

        \draw[mid arrow=>,thick] (-2,2.5) -- (-2,2);
        \draw[mid arrow=>,thick] (-2,1) -- (-2,0.5);
        \draw[mid arrow=>,thick] (-2,-0.5) -- (-2,-1);
        \draw[mid arrow=>,thick] (-2,-2) -- (-2,-2.5);
        \draw[mid arrow=>,thick] (-3,5.5) -- (-3,3.5);
        \draw[mid arrow=>,thick] (-3,-3.5) -- (-3,-5.5);

        \begin{scope}[shift = {(4,-1.5)}]

            \draw[mid arrow=<,thick] (-0.25,-0.5) .. controls +(-1,5) and +(1,-5) .. (-5,-2);

            \draw[mid arrow=>,thick] (-0.25,-1.5) -- (-0.25,-4);

            \node[draw, thick, fill = white, minimum height = 1cm, minimum width = 1cm] at (-0.25,-1) {$\eta'_{F(G)}$};

        \end{scope}
        \end{tikzpicture}
    \end{center}

    Then using naturality of the unit and counit, we have:

            
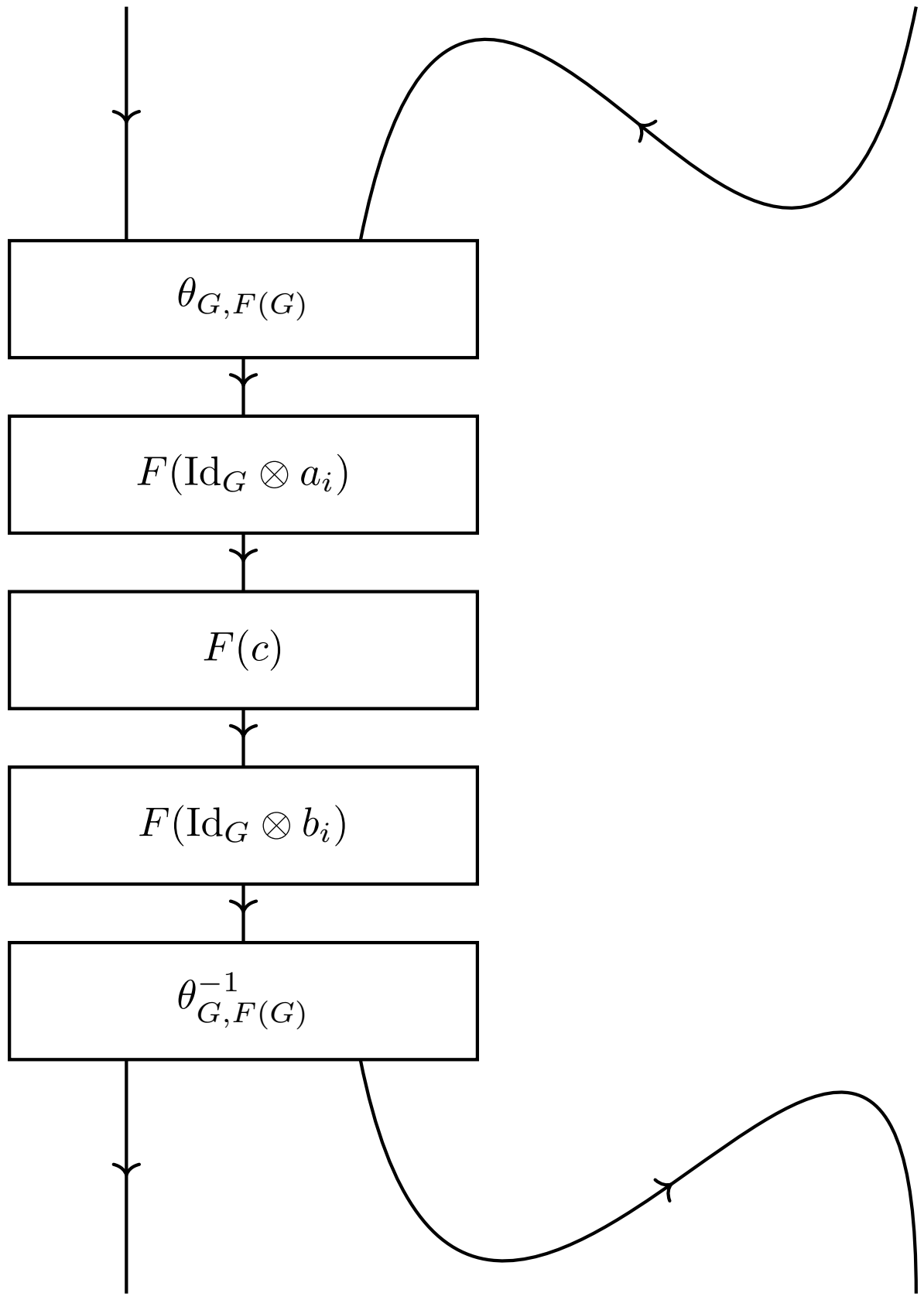
\begin{center}
        \begin{tikzpicture}
            \begin{scope}[shift = {(4,1.5)}]

            \draw[mid arrow=>,thick] (-0.25,4) .. controls +(-1,-5) and +(1,5) .. (-5,2);
        \end{scope}

        \node[draw, thick, fill = white, minimum width = 4cm, minimum height = 1cm]  at (-2,0) {$F(c)$};

        \node[draw, thick, fill = white, minimum width = 4cm, minimum height = 1cm]  at (-2,1.5) {$F(\Id_G \otimes a_i)$};
        \node[draw, thick, fill = white, minimum width = 4cm, minimum height = 1cm]  at (-2,-1.5) {$F(\Id_G \otimes b_i)$};

        \node[draw, thick, fill = white, minimum width = 4cm, minimum height = 1cm]  at (-2,3) {$\theta_{G,F(G)}$};
        \node[draw, thick, fill = white, minimum width = 4cm, minimum height = 1cm]  at (-2,-3) {$\theta_{G,F(G)}^{-1}$};

        \draw[mid arrow=>,thick] (-2,2.5) -- (-2,2);
        \draw[mid arrow=>,thick] (-2,1) -- (-2,0.5);
        \draw[mid arrow=>,thick] (-2,-0.5) -- (-2,-1);
        \draw[mid arrow=>,thick] (-2,-2) -- (-2,-2.5);
        \draw[mid arrow=>,thick] (-3,5.5) -- (-3,3.5);
        \draw[mid arrow=>,thick] (-3,-3.5) -- (-3,-5.5);

        \begin{scope}[shift = {(4,-1.5)}]

            \draw[mid arrow=<,thick] (-0.25,-4) .. controls +(0,5) and +(1,-5) .. (-5,-2);

        \end{scope}
        \end{tikzpicture}
    \end{center}

    Since $c$ is a chromatic map, $(\Id_G\otimes a_i)c(\Id_G\otimes b_i) = \Id_{G\otimes G}$, and hence this is the identity.
\end{proof}

\subsection{Ends, coends, and unimodularity}

Above we saw that when the functor $F$ has the generic structure of a left and right adjoint to the forgetful functor $U:\mcal Z(\mcal C)\to \mcal C$, then we obtained a non-degenerate trace. In this section, we describe a specific left and right adjoint using the end and coend.

\begin{dff}
    Let $\mcal C$ and $\mcal D$ be categories. Let $S$ and $T$ be functors from $\mcal C^{\mathrm{op}}\times \mcal C$ to $\mcal D$. A \emph{dinatural transformation} $\xi:S\xrightarrow{\cdot\cdot}T$ is a family \[\xi = \{\xi_X: S(X,X)\to T(X,X)\}_{X\in \mcal C}\] of morphisms in $\mcal D$ such that \[T(\Id_X, f)\circ \xi_X\circ S(f,\Id_X) = T(f,\Id_Y)\circ \xi_Y\circ S(\Id_Y,f)\] for all morphisms $f:X\to Y$ in $\mcal C$. We will consider an object $X\in \mcal D$ as a functor from $\mcal C^{\mathrm{op}}\times \mcal C$ to $\mcal D$ sending all morphisms to $\Id_X$.
\end{dff}

\begin{dff}
    An \emph{end} of $S:\mcal C^{\mathrm{op}}\times \mcal C\to \mcal D$ is a pair $(E,\pi)$ consisting of an object $E\in \mcal D$ and a dinatural transformation $\pi:E\xrightarrow{\cdot\cdot} S$ satisfying a universal property ensuring it is unique up to unique isomorphism if it exists. We say $\pi$ is a \emph{universal dinatural transformation}.

    A \emph{coend} of $T$ is a pair $(C,i)$ consisting of an object $C\in \mcal D$ and a universal dinatural transformation $i:T\xrightarrow{\cdot\cdot} C$. Then a coend is unique up to a unique isomorphism if it exists.
\end{dff}

We denote an end $(E,\pi)$ and a coend $(C,i)$, respectively, by
\[E = \int_{X\in \mcal C} S(X,X) \qquad \text{and}\qquad C = \int^{X\in \mcal C}T(X,X)\]

See \cite{mac1971categories} for a general treatment of (co)ends. 

Given a finite tensor $\Bbbk$-category $\mcal C$, the end \[R(V) = \int_{X\in \mcal C} X^*\otimes V\otimes X\] and coend \[L(V) = \int^{X\in \mcal C} X^*\otimes V\otimes X\] exist for any $V\in \mcal C$. The map $R$ is an endofunctor of $\mcal C$ and its universal dinatural transformations $\pi_{V,X}:R(V)\to X^*\otimes V\otimes X$ endow it with a comonad structure $(R,\delta,\epsilon)$. Similarly, $L$ has the structure of a monad $(L,m,\eta)$ coming from its universal dinatural transformations $i_{V,X}:X^*\otimes V\otimes X\to L(V)$. The comonad and monad structures endow $R(V)$ and $L(V)$ with canonical half braidings

\begin{align*}
    \sigma^R_{X} &= (\Id_{X\otimes R(V)}\otimes \rev_X)(\pi_{R(V),X^*}\otimes \Id_X)(\delta_V\otimes \Id_X)\\
    (\sigma^R_{X})^{-1} &= (\lev_X\otimes \Id_{R(V)\otimes X})(\Id_X \otimes \pi_{R(V),X})(\Id_X \otimes \delta_V)\\
    \sigma^L_{X} &= (\Id_X \otimes m_V)(\Id_X \otimes i_{L(V),X})(\rcoev_X \otimes \Id_{L(V)\otimes X})\\
    (\sigma^L_{X})^{-1} &= (m_V \otimes \Id_X)(i_{L(V),X^*}\otimes \Id_X)(\Id_{X\otimes L(V)}\otimes \lcoev_X)\\
\end{align*}

and hence for each $V\in \mcal C$, $R$ and $L$ define functors to the center $\mcal Z(\mcal C)$ such that $R$ is right adjoint to $U$ and $L$ is left adjoint to $U$.

The unit for $L\dashv U$ is given by $\eta^L_V = i_{V,\mbb 1}$ and the counit for $U\dashv R$ is given by $\epsilon^R_V = \pi_{V,\mbb 1}$.

For unimodular categories, $R\cong L$ (\cite{shimizu2016unimodular}). Following \cite{schweigert2023differential}, we can compute this isomorphism using the non-degenerate m-trace.

For any tensor category $\mcal C$, the maps $Y \otimes \Hom_{\mcal C}(Y,X)\to X$ for $X,Y\in\mcal C$ induce maps 
\[X^* \otimes P \otimes X \to ( Y\otimes \Hom_{\mcal C}(Y, X))^* \otimes P \otimes X \cong Y^* \otimes \Hom_{\mcal C}(Y,X)^* \otimes P \otimes X \to Y^* \otimes P \otimes N^r(Y)\] where $N^r$ denotes the Nakayama functor \[ N^rX = \int^{Y\in \mcal C}\Hom_{\mcal C}(X,Y)^*\otimes Y\] 

Then these maps descend to the coend $L(P)$ and factor through the end $R(P)$ to yield an isomorphism \[L(P) \to \int_{Y\in \mcal C} Y^*\otimes P \otimes N^rY\]

In a unimodular category the Nakayama functor $N^r$ is naturally isomorphic to the identity functor, so we obtain an isomorphism \[\Omega_P:L(P)\to R(P)\]

Note that by \cite[Proposition 5.1.7]{kerler2002non}, we have that \[L(P) \cong \int^{X\in \mathrm{Proj}(\mcal C)} X^* \otimes P\otimes X\] and \[R(P) \cong \int_{X\in \mathrm{Proj}(\mcal C)} X^*\otimes P \otimes X\] so we need only consider the projective components of our map $\Omega_P$. Then as in \cite{schweigert2023differential}, the components $\pi_{P,Y}U(\Omega_P) i_{P,X}$ are given by the composition
\begin{align*}
    X^* \otimes P \otimes X \to (Y\otimes \Hom_{\mcal C}(Y, X))^* \otimes P \otimes X &\cong Y^* \otimes \Hom_{\mcal C}(Y, X)^* \otimes P \otimes X\\
    &\cong Y^*\otimes \Hom_{\mcal C}(X,Y)\otimes P\otimes X \to Y^* \otimes P \otimes Y
\end{align*}

The isomorphism $\Hom_{\mcal C}(Y,X)^* \cong \Hom_{\mcal C}(X, Y)$ is given by the commuting square

\begin{center}
\begin{tikzcd}
{\Hom_{\mathcal C}(Y,X)^*} \arrow[rrr, "{\Hom_{\mathcal C}(Y,X)^*\otimes \Id_{X}}"] \arrow[d, "T"'] &  &  & {\Hom_{\mathcal C}(Y,X)^*\otimes \Hom_{\mathcal C}(X, X)} \arrow[d, "{\Hom_{\mathcal C}(Y,X)^* \otimes \Delta_{X,Y}}"] \\
{\Hom_{\mathcal C}(X, Y)}                                                                                                 &  &  & {\Hom_{\mathcal C}(Y,X)^* \otimes \Hom_{\mathcal C}(X,Y) \otimes \Hom_{\mathcal C}(Y,X)} \arrow[lll, "\text{evaluation}"]      
\end{tikzcd}
\end{center}

The isomorphism $T:\Hom_{\mcal C}(Y,X)^* \to \Hom_{\mcal C}(X,Y)$ is induced by the non-degenerate trace $\ctr$. Let $\{a_i\}$ and $\{b_i\}$ be bases for $\Hom_{\mcal C}(Y,X)$ and $\Hom_{\mcal C}(X,Y)$, respectively, dual under $\ctr$, and denote by $\langle f, \bullet\rangle$ in $\Hom_{\mcal C}(Y,X)^*$ the linear dual of $f\in \Hom_{\mcal C}(Y,X)$, then $T$ is defined by
\[T(\langle a_i, \bullet\rangle) = b_i\]

We can use $T$ to determine $\Delta_{X,Y}$. Say $\Delta_{X,Y}(\Id_{X}) = \alpha'_{X,Y}\otimes \alpha''_{X,Y}$ using Sweedler notation. The commuting square yields the equality
\[b_i = \langle a_i, \alpha''_{X,Y}\rangle \alpha'_{X,Y}\]

This implies $\Delta_{X,Y} = \sum_i b_i \otimes a_i$ and the components $\pi_{P,Y}\Omega_P i_{P,X}: X^*\otimes P \otimes X \to Y^*\otimes P\otimes Y$ are given by the formula \[\pi_{P,Y}U(\Omega_P) i_{P,X} = \sum_i a_i^* \otimes \Id_P \otimes b_i\] Graphically,

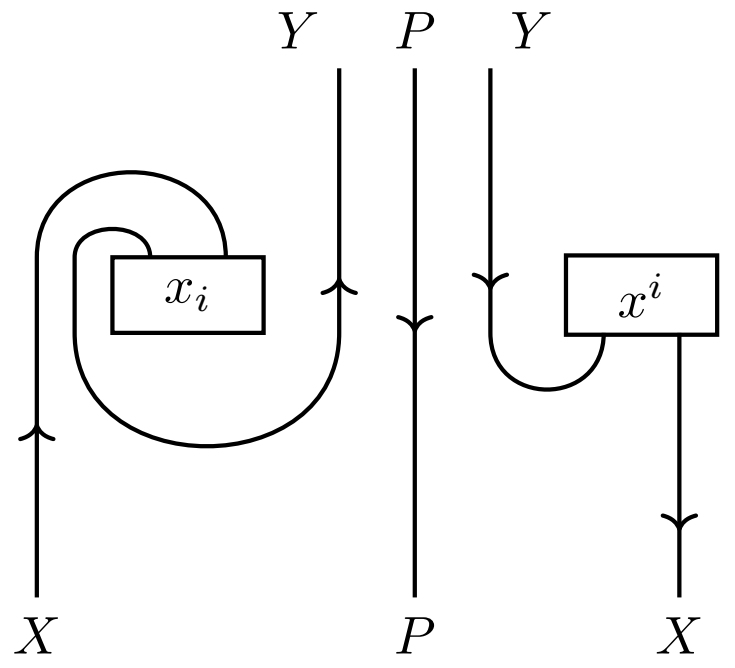
\begin{center}
        \begin{tikzpicture}
            \node at (-5,0) {$\pi_{P,Y}U(\Omega_P) i_{P,X} = $};
            \node at (-3,0) {$\displaystyle\sum_i$};
            
            \node[draw,thick,fill=white,minimum height=0.5cm,minimum width=1cm] at (-1.5,0.5) {$x_i$};
            \node[draw,thick,fill=white,minimum height=0.5cm,minimum width=1cm] at (1.5,0.5) {$x^i$};
            
            \draw[thick,mid arrow=>] (0,2)node[above] {$P$} -- (0,-1.5)node[below] {$P$};
            \draw[thick,third arrow=<] (-1.25,0.75) .. controls +(0,0.75) and +(0,0.75) .. (-2.5,0.75) -- (-2.5,-1.5)node[below] {$X$};
            \draw[thick,third arrow=>] (-1.75,0.75) .. controls +(0,0.25) and +(0,0.25) .. (-2.25,0.75) -- (-2.25,0.25) .. controls +(0,-1) and +(0,-1) .. (-0.5,0.25) -- (-0.5,2)node[above left] {$Y$};

            \draw[thick,mid arrow=>] (0.5,2)node[above right] {$Y$} -- (0.5,0.25) .. controls +(0,-0.5) and +(0,-0.5) .. (1.25,0.25);
            \draw[thick,third arrow=>] (1.75,0.25) -- (1.75,-1.5)node[below] {$X$};
        \end{tikzpicture}
    \end{center}

    where $X$ and $Y$ are projective and $a_i$ and $b_i$ are defined using $\Lambda_{Y\otimes X} = \sum_i x_i x^i$.

    \begin{lem}\label{lem:unit_counit_halfbraiding_LP} For any $X,P \in \mathrm{Proj}(\mcal C)$,
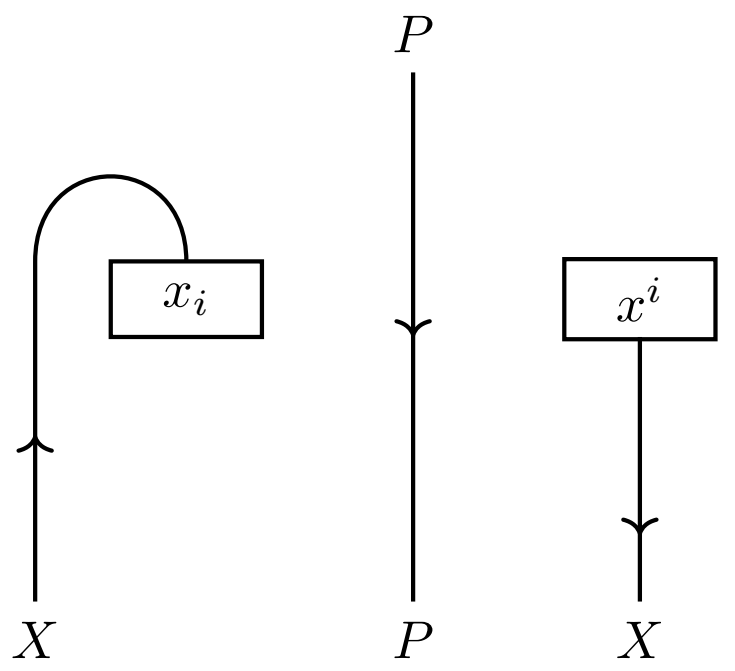
        \begin{center}
        \begin{tikzpicture}
            \node at (-5,0) {$\pi_{P,\mbb 1}U(\Omega_P) i_{P,X} = $};
            \node at (-3,0) {$\displaystyle\sum_i$};
            
            \node[draw,thick,fill=white,minimum height=0.5cm,minimum width=1cm] at (-1.5,0.5) {$x_i$};
            \node[draw,thick,fill=white,minimum height=0.5cm,minimum width=1cm] at (1.5,0.5) {$x^i$};
            
            \draw[thick,mid arrow=>] (0,2)node[above] {$P$} -- (0,-1.5)node[below] {$P$};
            \draw[thick,third arrow=<] (-1.5,0.75) .. controls +(0,0.75) and +(0,0.75) .. (-2.5,0.75) -- (-2.5,-1.5)node[below] {$X$};
            \draw[thick,third arrow=>] (1.5,0.25) -- (1.5,-1.5)node[below] {$X$};
        \end{tikzpicture}
    \end{center}

    In particular, if $\sigma_X:L(P)\otimes X\to X\otimes L(P)$ is the canonical half braiding for the coend $L(P)$, then \[(\Id_X\otimes \pi_{P,\mbb 1}U(\Omega_P))\sigma_X(i_{P,\mbb 1}\otimes \Id_X) = \sum_{i}(x_i \otimes \Id_P)(\Id_P\otimes x^i)\] 
    \end{lem}

    \begin{proof}
        Consider a projective generator $G$. Then there exists an epimorphism $e: G\twoheadrightarrow \mbb 1$ and a morphism $\delta_P:P\to G\otimes P$ such that $(e\otimes \Id_P)\delta_P = \Id_P$. Indeed, since $G\otimes G^*$ is projective and is a retract of a finite direct sum of copies of $G$, the (nonzero) evaluation epimorphism $\lev_G:G\otimes G^*\to \mbb 1$ factors through G, and so there exists a nonzero map $e:G\to \mbb 1$ which then is an epimorphism. Since $P\otimes P^*$ is projective, the evaluation morphism $\lev_P:P\otimes P^*\to \mbb 1$ factors as $\lev_P = e\circ d_{P,G}$ for some $d_{P,G}:P\otimes P^* \to G$. Then the map $\delta_{P} = (d_{P,G}\otimes \Id_P)(\Id_P\otimes \lcoev_P)$ does satisfy $\Id_P = (e\otimes \Id_P) \delta_{P}$.

        Then, consider the following:
        \begin{align*}
            \pi_{P,\mbb 1}U(\Omega_P) i_{P,X} &= (e\otimes \Id_P)\delta_P\pi_{P,\mbb 1}U(\Omega_P) i_{P,X}\\
            &= (\rev_G \otimes \Id_P)(e^* \otimes \delta_P)\pi_{P,\mbb 1}U(\Omega_P) i_{P,X}\\
            &= (\rev_G\otimes \Id_P)(\Id_{G^*}\otimes \delta_P \otimes e)\pi_{P,G}U(\Omega_P) i_{P,X}\\
            &= (\rev_G\otimes \Id_P)(\Id_{G^*}\otimes \delta_P \otimes e)(\sum_i a_i^* \otimes \Id_P \otimes b_i)\\
            &= (\rev_G \otimes \Id_P)(e^* \otimes \delta_P)(\sum_i (x_i)^* \otimes \Id_P \otimes x^i)\\
            &= \sum_i (x_i)^* \otimes \Id_P \otimes x^i
        \end{align*}
        as desired, where we used dinaturality in line 3.

        To obtain the equality with the canonical half-braiding, use that $\sigma_X(i_{P,\mbb 1}\otimes \Id_X) = (\Id_X\otimes i_{P,X})(\rcoev_X \otimes \Id_{P\otimes X})$ which follows since $i_{P,\mbb 1}$ is the unit for the monad $L$.
    \end{proof}

\subsection{The coend as a right adjoint and the modified trace}

In this section, we use the natural unimodular isomorphism $\Omega$ to describe the coend $L(P)$ as a right adjoint to the forgetful functor. Additionally, we show that this right adjoint defines the non-degenerate modified trace in the center.

Let $L:\mcal C\to \mcal Z(\mcal C)$ be the coend and $R:\mcal C\to \mcal Z(\mcal C)$ be the end.

The counit for $L\dashv U$ is the unique morphism $\epsilon^L_M:LU(M) \to M$ such that if $\sigma_X:M\otimes X \to X\otimes M$ is the half-braiding for $M$, then $U(\epsilon^L_M) i_{U(M),X} = (\rev_X \otimes \Id_{U(M)})(\Id_{X^*}\otimes \sigma_X)$ (see \cite{turaev2017monoidal} for more details). There is a similar universal property for the unit for $U\dashv R$ using the inverse half-braiding.

Let $\theta^L_{P,M}: L(P\otimes U(M))\to L(P) \otimes M$ and $\theta^R_{P,M}:R(P\otimes U(M)) \to R(P)\otimes M$ be the projection isomorphisms defined by \[\theta^L_{P,M} = \epsilon^L_{L(P)\otimes M}L(i_{P,\mbb 1} \otimes \Id_{U(M)})\] and \[(\theta^R_{P,M})^{-1} = R(\pi_{P,\mbb 1} \otimes \Id_{U(M)})\eta^R_{R(P)\otimes M}\]

Then $\theta^L$ and $\theta^R$ satisfy \[U(\theta^L_{P,M})i_{P\otimes U(M),\mbb 1} = i_{P,\mbb 1} \otimes \Id_{U(M)}\] and \[(\pi_{P,\mbb 1}\otimes \Id_{U(M)})U(\theta^R_{P,M}) = \pi_{P\otimes U(M),\mbb 1}\]

\begin{lem}
    The projection isomorphisms $\theta^L$ and $\theta^R$ satisfy
    \[U(\theta^L_{P,M}) i_{P\otimes U(M),X} = (i_{P,X}\otimes \Id_{U(M)})(\Id_{X^*\otimes P} \otimes \sigma_X)\] and 
    \[(\Id_{X^*\otimes P}\otimes \sigma_X^{-1})(\pi_{P,X}\otimes \Id_{U(M)})U(\theta^R_{P,M}) = \pi_{P\otimes U(M),X}\]
\end{lem}

\begin{proof}
    We prove the statement for $\theta^L$, then the statement for $\theta^R$ is analogous.

    \begin{align*}
        U(\theta^L_{P,M}) i_{P\otimes U(M),X} &= U(\epsilon^L_{L(P)\otimes M})UL(i_{P,\mbb 1}\otimes \Id_{U(M)}) i_{P\otimes U(M),X}\\
        &= U(\epsilon^L_{L(P)\otimes M}) i_{L(P)\otimes M,X} (\Id_{X^*}\otimes i_{P,\mbb 1} \otimes \Id_{U(M)} \otimes \Id_X)\\
        &= (\rev_X \otimes \Id_{UL(P)\otimes U(M)})(\Id_{X^*}\otimes \widetilde{\sigma}_X)\\
        &= (i_{P,X}\otimes \Id_{U(M)})(\Id_{X^*\otimes P} \otimes \sigma_X)
    \end{align*}

    where $\widetilde{\sigma}_X$ is the half-braiding for $L(P) \otimes M$. The last equality follows from the fact that $\sigma^L_X (i_{P,\mbb 1}\otimes \Id_X) = (\Id_X\otimes i_{P,X})(\rcoev_X \otimes \Id_{P\otimes X})$, where $\sigma^L_X:L(P)\otimes X \to X\otimes L(P)$ is the canonical half-braiding for $L(P)$.
\end{proof}

\begin{lem}
    The natural unimodular isomorphism $\Omega$ is a right $\mcal Z(\mcal C)$-module isomorphism in the sense that \[ \theta^R_{P,M} \Omega_{P\otimes U(M)} = (\Omega_P\otimes \Id_{U(M)})\theta^L_{P,M}\]
\end{lem}

\begin{proof}
    We check their components are equal:

    \begin{align*}
        (\pi_{P,Y}\otimes \Id_{U(M)})U(\theta^R_{P,M}) U(\Omega_{P\otimes U(M)})i_{P\otimes U(M),X} &= (\Id_{Y^*\otimes P}\otimes \sigma_Y)\pi_{P\otimes U(M),Y} U(\Omega_{P\otimes U(M)})i_{P\otimes U(M),X}\\
        &= (\pi_{P,Y}U(\Omega_P)i_{P,X} \otimes \Id_{U(M)})(\Id_{X^*\otimes P}\otimes \sigma_X)\\
        &= (\pi_{P,Y} \otimes \Id_{U(M)})(U(\Omega_P) \otimes \Id_{U(M)})U(\theta^L_{P,M})i_{P\otimes U(M),X}
    \end{align*}
\end{proof}

Now, recall the counit $\epsilon'_P:UF(P)\to P$ where $F$ is the two-sided adjoint to $U$. Then, taking $F = L$, we see that $\epsilon'_P$ is the counit for an adjunction $U\dashv L$. Since $L\cong R$, such a counit exists. Now, using the adjunction isomorphism $\Hom_{\mcal C}(UL(P),P) \to \Hom_{\mcal Z(\mcal C)}(L(P), R(P))$, we see that there is a natural isomorphism $\nu: L\to R$ such that $\epsilon'_P = \pi_{P,\mbb 1}\circ U(\nu_P)$. In particular, $\nu_P = R(\epsilon'_P)\circ \eta^R_{L(P)}$.

\begin{lem}
    The natural isomorphism $\nu$ is a right $\mcal Z(\mcal C)$-module isomorphism. I.e. \[\theta^R_{P,M}\nu_{P\otimes U(M)} = (\nu_P\otimes \Id_M)\theta^L_{P,M}\]
\end{lem}

\begin{proof}
    Recall that $\epsilon'_{P\otimes U(M)} = (\epsilon'_P\otimes \Id_{U(M)})U(\theta^L_{P,M})$. Then 
    
    \begin{align*}
        \nu_{P\otimes U(M)} &= R(\epsilon'_P\otimes \Id_{U(M)})RU(\theta^L_{P,M})\eta^R_{L(P\otimes U(M))}\\
        &= R(\epsilon'_P \otimes \Id_{U(M)})\eta^R_{L(P)\otimes M}\theta^L_{P,M}\\
        &= R(\pi_{P,\mbb 1} U(\nu_P) \otimes \Id_{U(M)})\eta^R_{L(P)\otimes M}\theta^L_{P,M}\\
        &= R(\pi_{P,\mbb 1} \otimes \Id_{U(M)}) RU(\nu_P \otimes \Id_{M})\eta^R_{L(P)\otimes M}\theta^L_{P,M}\\
        &= R(\pi_{P,\mbb 1} \otimes \Id_{U(M)})\eta^R_{R(P)\otimes M}(\nu_P \otimes \Id_{M})\theta^L_{P,M}\\
        &=(\theta^R_{P,M})^{-1}(\nu_P \otimes \Id_{M})\theta^L_{P,M}\\
    \end{align*}
\end{proof}

Thus $a = \nu^{-1} \Omega:L\to L$ is a right $\mcal Z(\mcal C)$-module isomorphism, so that \[\theta^L_{P,M} a_{P\otimes U(M)} = (a_P \otimes \Id_M)\theta^L_{P,M}\]

Now, define $b_M = U(\epsilon^L_M) U(a_{U(M)})i_{U(M),\mbb 1}$. Then since $\epsilon^L$ is a right module homomorphism, we have

\begin{align*}
    b_{N\otimes M} &= U(\epsilon^L_{N\otimes M}) U(a_{U(N)\otimes U(M)}) i_{U(N)\otimes U(M),\mbb 1}\\
    &= U(\epsilon^L_{N\otimes M} a_{U(N)\otimes U(M)}) i_{U(N)\otimes U(M),\mbb 1}\\
    &= U((\epsilon^L_N \otimes \Id_M)\theta^L_{U(N),M} a_{U(N)\otimes U(M)}) i_{U(N)\otimes U(M),\mbb 1}\\
    &= U(\epsilon^L_N a_{U(N)} \otimes \Id_M)U(\theta^L_{U(N),M})i_{U(N)\otimes U(M),\mbb 1}\\
    &= U(\epsilon^L_N a_{U(N)} \otimes \Id_M)(i_{U(N),\mbb 1}\otimes \Id_{U(M)})\\
    &= b_N \otimes \Id_{U(M)}
\end{align*}

Thus $b_M = \lambda \Id_{U(M)}$ for some non-zero scalar $\lambda$. Thus $a_P$ is uniquely determined by $a_{\mbb 1}\in \Hom_{\mcal Z(\mcal C)}(L(\mbb 1),L(\mbb 1))$ and we have that for each $M$, \[\epsilon^L_M a_{U(M)} = \lambda \epsilon^L_M\] and in particular for $M = L(\mbb 1)$, we have

\[\epsilon^L_{L(\mbb 1)}a_{UL(\mbb 1)} L(i_{\mbb 1,\mbb 1}) = \lambda \epsilon^L_{L(\mbb 1)}L(i_{\mbb 1,\mbb 1})\]

The right-hand side is $\lambda \Id_{L(\mbb 1)}$ and the left-hand side, using naturality of $a$ is $a_{\mbb 1}$. Thus for any $X\in \mcal C$, we have $a_X = \lambda \Id_{L(X)}$.

Then $\lambda \nu = \Omega$ and hence $\epsilon'_X = \frac{1}{\lambda} \pi_{X,\mbb 1}U(\Omega_X)$ for each object $X$.

Thus, rescaling by $\lambda$ yields the non-degenerate m-trace 
\[\ztr_{L(P)}(f) = \ctr_P (\pi_{P,\mbb 1}U(\Omega_P) U(f) i_{P,\mbb 1})\] Set $\theta_{P,M} = \theta^L_{P,M}$, then the previous work verifies \[\theta^{-1}_{P,M} = L(\pi_{P,\mbb 1}U(\Omega_P) \otimes \Id_{U(M)}) \Omega_{UL(P)\otimes U(M)}^{-1} \eta^R_{L(P)\otimes M}\]

Finally, we obtain a chromatic map for the center in these terms:
\[c^{\mcal Z}_{L(G)} = \theta_{G,L(G)}L(\Id_G\otimes i_{G,\mbb 1})L(c)L(\Id_G\otimes \pi_{G,\mbb 1}U(\Omega_G)) (\theta)_{G,L(G)}^{-1}\]

\begin{thm}\label{thm:chromatic_halfbraiding_identity}
    Let $c^{\mcal Z}_{L(G)}$ be the chromatic morphism above for $L(G)$ and define for each $X\in \mcal C$ and each projective $P\in\mcal C$,
    \begin{align*}
        A_{P,X} &= \sum_i(\Id_{UL(G)\otimes P}\otimes \lev_{UL(G)} \otimes \Id_X)(\Id_{UL(G)} \otimes \sigma^L_P \otimes x_i)\\
        &\qquad \circ(U(c^{\mcal Z})\otimes x^i \otimes \Id_P)(\Id_{UL(G)} \otimes \rcoev_{UL(G)} \otimes \Id_{X\otimes P})
    \end{align*}

    where $\Lambda_{UL(G)^* \otimes X} = \sum_i x_i x^i$.

    Then $A_{P,X} = \Id_{UL(G)} \otimes \sum_i (y_i \otimes \Id_X)(\Id_X\otimes y^i)$ where $\Lambda_P = \sum_i y_i y^i$.

    Graphically,

    \begin{center}
        \begin{tikzpicture}
            \node at (-1.5,0) {$\displaystyle\sum_i$};
            \draw[mid arrow=>,thick] (-0.5,2) -- (-0.5,0.5);
            \draw[mid arrow=>,thick] (-0.5,-0.5) -- (-0.5,-1.5)node[below]{$UL(G)$};
            \draw[mid arrow=>,thick] (2.75,2) -- (2.75,1.25);
            \draw[mid arrow=>,thick] (2.75,-0.75) -- (2.75,-1.5)node[below]{$X$};
            \draw[mid arrow=>,thick] (1,1) .. controls +(0,-1) and +(0,2.5) .. (4,-1.5);
            \node[below] at (4,-1.5) {$P$};
            \draw[mid arrow=<,thick] (1,1.55) .. controls +(0,0.5) and +(0,1) .. (2.25,1.25);
            \draw[mid arrow=>,thick] (0.5,-0.5) .. controls +(0,-1) and +(0,-0.5) .. (2.25,-0.75);
            \draw[mid arrow=>,thick] (0.5,2) -- (0.5,1.55);
            \draw[mid arrow=>,thick] (0.5,1) -- (0.5,0.5);

            \node[draw,thick,fill = white,minimum width = 0.75cm, minimum height = 0.6cm] at (0.75,1.25) {$\sigma^L_P$};
            \node[draw,thick,fill = white,minimum width = 2cm, minimum height = 1cm] at (0,0) {$U(c^{\mcal Z}_{L(G)})$};
            \node[draw,thick,fill = white,minimum width = 1cm, minimum height = 0.5cm] at (2.5,1) {$x_i$};
            \node[draw,thick,fill = white,minimum width = 1cm, minimum height = 0.5cm] at (2.5,-0.5) {$x^i$};

            \node at (4.5,0) {$=$};

            \begin{scope}[shift = {(5.5,0)}]
                \node at (0,0) {$\displaystyle\sum_i$};
                \draw[mid arrow=>,thick] (0.5,1.5) -- (0.5,-1.5)node[below]{$UL(G)$};
                
                \draw[mid arrow=>,thick] (1.5,1.5) -- (1.5,1);
                \draw[mid arrow=>,thick] (2.5,-1) -- (2.5,-1.5)node[below]{$P$};
                \node[draw,thick,fill=white,minimum width = 0.5cm, minimum height = 0.5cm] at (1.5,0.75) {$y_i$};
                \node[draw,thick,fill=white,minimum width = 0.5cm, minimum height = 0.5cm] at (2.5,-0.75) {$y^i$};

                \draw[mid arrow=>, thick] (2.5,1.5) .. controls +(0,-1) and +(0,1) .. (1.5,-1.5)node[below] {$X$};
                
            \end{scope}
        \end{tikzpicture}
    \end{center}

    There is an analogous statement for the inverse braiding $(\sigma^L_P)^{-1}$.
\end{thm}

\begin{proof}
    Note that the map \[a_i = (\lev_{UL(G)}\otimes \Id_X)(\Id_{UL(G)}\otimes x_i)\] is equal to \[\pi_{X,\mbb 1}U(\Omega_X) UL(a_i) U(\eta^R_{L(G)})\]
    Similarly we have that the map 
    \[b_i = (\Id_{UL(G)}\otimes x^i)(\rcoev_{UL(G)}\otimes \Id_X)\] is equal to \[U(\epsilon^L_{L(G)})UL(b_i)i_{X,\mbb 1}\]

    We note that $L(a_i)\eta^R_{L(G)}$ and $\epsilon^L_{L(G)}L(b_i)$ assemble to define the cutting morphism in the center, and hence we may slide those portions across the half braiding to obtain the chromatic identity in the center.

    Then the resulting morphism is $\Id_{UL(G)}\otimes (\Id_P \otimes \pi_{X,\mbb 1}U(\Omega_X))(\sigma^L_X)(i_{X,\mbb 1})$. By Lemma \ref{lem:unit_counit_halfbraiding_LP} this yields the right-hand side.
\end{proof}

\section{Topological Preliminaries}
\label{Topology_Sec}

In this section we develop the topological background to build an invariant of 3-manifolds with embedded ribbon graphs, and eventually embedded framed links. We do this by way of \emph{Heegaard-Graph Diagrams} and eventually \emph{Heegaard-Link Diagrams}. 

\subsection{Ribbon graphs in surfaces}

Let $\mcal C$ be a pivotal $\Bbbk$-category. To each free end of a $\mcal C$-colored ribbon graph $\Gamma$ in $\mbb R \times [0,1]$ where its free ends lie on the boundary, we associate a signed object consisting of the color of the arc incident to the free end and of a sign $\pm 1$ depending if that arc is directed up or down. Then one can view $\Gamma$ as a morphism from the sequence of signed objects associated with its bottom free ends in $\mbb R\times \{0\}$ to the sequence of signed objects associated with its top free ends in $\mbb R\times \{1\}$. This defines a monoidal category $\mathrm{Rib}_{\mcal C}$ whose objects are finite sequences of signed objects, whose morphisms are isotopy classes of $\mcal C$-colored ribbon graphs, whose composition is given by stacking, and whose monoidal product is given by concatenation. The graphical calculus gives rise to a monoidal functor \[\mcal F:\mathrm{Rib}_{\mcal C}\to \mcal C\] which we use to interpret these planar diagrams.

If the left and right traces on $\mcal C$ coincide, then $\mcal F$ induces an isotopy invariant $\mcal F: \mcal L \to \Endo_{\mcal C}(\mbb 1) = \Bbbk$, where $\mcal L$ is the class of $\mcal C$-colored ribbon graphs in $S^2 = (\mbb R\times (0,1))\cup \{\infty\}$.

The modified trace can be used to modify the isotopy invariant $\mcal F:\mcal L\to \Bbbk$ above as follows. Let $\mcal L_{\mcal I}$ denote the class of $\mcal C$-colored ribbon graphs in $S^2$ having at least one strand colored with an object in a tensor ideal $\mcal I$. In particular, each $\Gamma\in \mcal L_{\mcal I}$ is the braid closure of some $\mcal C$-colored ribbon graph $T_X$ in $\mbb R\times [0,1]$ with exactly one bottom free end and one top free end both supported by arcs oriented downward and colored by some object $X\in \mcal I$ so that $\mcal F(T_X)\in \Endo_{\mcal C}(X)$. Then, by \cite{Geer2020kuperberg}, each m-trace $\mathrm{t}$ on $\mcal I$ induces an isotopy invariant \[\mcal F':\mcal L_{\mcal I}\to \Bbbk\] sending $\Gamma \mapsto \mcal F'(\Gamma) = \mathrm{t}_X(\mcal F(T_X))$.

The Penrose graphical calculus for ribbon graphs may be extended from the plane to any surface $\Sigma$ in the following way. If an oriented arc has a free end, its free end must be placed on a boundary component of $\Sigma$. We will call such a ribbon graph $\Gamma$ (in $\Sigma$) \emph{$\mcal I$-admissible} if at least one strand of $\Gamma$ is colored by an object in $\mcal I$.

Given $\mcal I$-admissible ribbon graphs $\Gamma_1,\dots ,\Gamma_k$ in $\Sigma$ and $a_1,\dots ,a_k\in \Bbbk$, the linear combination $a_1\Gamma_1 + a_2\Gamma_2 + \cdots + a_k\Gamma_k$ is a \emph{$\mcal I$-skein relation} (in $\Sigma$) if there is a coupon $Q$ embedded in $\Sigma$ and $\mcal I$-admissible ribbon graphs $\Gamma_1',\dots,\Gamma_k'$ in $\Sigma$ such that:
\begin{itemize}
    \item $\Gamma_i'$ is isotopic to $\Gamma_i$ as a $\mcal C$-colored ribbon graph in $\Sigma$ for all $1\leq i \leq k$;
    \item the $\Gamma'_i$s coincide outside $Q: \Gamma'_i\cap (\Sigma -Q) = \Gamma'_j \cap (\Sigma - Q)$ for all $1 \leq i,j\leq k$;
    \item $\Gamma_i'$ intersects $\partial Q$ only in its bottom and top bases and transversally along the strands of $\Gamma_i'$ (so that $\Gamma'_i\cap Q$ can be seen as a $\mcal C$-colored ribbon graph in $\mbb R\times [0,1]$) for all $1\leq i\leq k$;
    \item $a_1 \mcal F(\Gamma_1'\cap Q) + \cdots + a_k \mcal F(\Gamma_k'\cap Q) = 0$ as a morphism in $\mcal C$;
    \item each $\Gamma_i'$ has an edge colored by an object in $\mcal I$ which is not entirely contained in the coupon $Q$.
\end{itemize}

Two linear combinations of $\mcal I$-admissible ribbon graphs are \emph{$\mcal I$-skein equivalent} if their difference is a $\mcal I$-skein relation.

\subsection{Heegaard-Link diagrams and surgery}

We review Heegaard diagrams:

\begin{dff}
    A (genus $g$) \emph{Heegaard diagram} is a triple $(\Sigma,\alpha,\beta)$, where $\Sigma$ is a closed connected oriented surface of genus $g$, and $\alpha = \{\alpha_1,\dots ,\alpha_g\}$ and $\beta = \{\beta_1,\dots ,\beta_g\}$ are two sets of $g$ simple closed curves embedded in $\Sigma$ such that:
    \begin{itemize}
        \item $\alpha_i\cap \alpha_j = \emptyset$ if $i\neq j$;
        \item $\beta_i \cap \beta_j = \emptyset$ if $i \neq j$;
        \item For any $\alpha_i$, $\beta_j$, their intersection is transversal;
        \item The surfaces with boundary $\Sigma - \alpha$ and $\Sigma - \beta$ are each homeomorphic to a 2-sphere with $2g$ boundary components.
    \end{itemize}
\end{dff}

A Heegaard diagram defines a closed connected oriented 3-manifold in the following way. Take $\Sigma \times [0,1]$ and attach 2-handles $D^2 \times D^1$ along each $\alpha\times \{0\}$ and 2-handles along each $\beta\times \{1\}$. Then the boundary of this manifold is two spheres, which may be filled in with 3-balls $D^3$. Any closed connected oriented 3-manifold may be obtained in this manner from a Heegaard diagram.

Given a 3-manifold $M$, a self-indexing Morse function $f:M\to \mbb R$ defines a Heegaard diagram with surface $\Sigma = f^{-1}(3/2)$.

Any two Heegaard diagrams for homeomorphic 3-manifolds may be connected via a finite sequence of \emph{Heegaard Moves}.
\begin{dff}[Heegaard Moves]~

    \begin{enumerate}[Move 1:]
        \item (Diffeomorphism of Surface) Performing a diffeomorphism of $\Sigma$ intertwining the curve sets $\alpha$ and $\beta$.
        \item (Diagram Isotopy) Performing an isotopy of the curve sets $\alpha$ and $\beta$ in $\Sigma$, not necessarily extending to an isotopy of $\Sigma$.
        \item (Handle Slide) If $\gamma$ is a curve in $\Sigma$ connecting $\alpha_i$ to $\alpha_j$ (resp. $\beta_i$ to $\beta_j$), then slide $\alpha_i$ along $\alpha_j$ via $\gamma$.
        \item (Stabilization/Destabilization) Performing a connected sum of the standard genus 1 diagram for $S^3$ (stabilization), or removing a connected sum component as such (destabilization).
    \end{enumerate}
\end{dff}

Note that the pair $(\Sigma,\alpha)$ defines a genus $g$ handlebody.

\begin{dff}
    Let $M$ be a closed connected oriented 3-manifold. A \emph{knot} is an embedding $K:S^1\to M$. A disjoint union of knots with $m$ connected components is called a \emph{link} with $m$ components. We will not distinguish a knot from its image. A \emph{framing} for a knot $K$ is a choice of trivialization of its normal bundle in $M$. This amounts to choosing a knot parallel to $K$. A framing for a link $L$ is a framing for each of its components.
\end{dff}

Suppose a link $L$ is embedded in a surface $\Sigma$ in a 3-manifold $M$. Then we say $L$ is endowed with the \emph{blackboard framing} if the choice of trivialization is a nonzero section of the tangent bundle of $\Sigma$.

\begin{dff}
    Given a closed connected oriented 3-manifold $M$ and a framed link $L$ in $M$, a \emph{Heegaard-Link diagram} for the pair $(M,L)$ is a quadruple $(\Sigma,\alpha,\beta,\tilde L)$, where $(\Sigma,\alpha,\beta)$ is a Heegaard diagram for $M$ and $\tilde L$ is a set of simple closed curves properly embedded in $\Sigma$ such that the pair $(\tilde M, \tilde L)$ comprised of the 3-manifold $\tilde M$ built from the Heegaard diagram and $\tilde L$ treated as a framed link (induced by the blackboard framing in the Heegaard surface $\Sigma$) in $\tilde M$ is diffeomorphic with the pair $(M,L)$ where the restriction to a tubular neighborhood of $\tilde L$ is a framed diffeomorphism. We additionally require that each component of the link $\tilde L$ intersect at least one $\alpha_i\in \alpha$.
    
    We say that two Heegaard-Link diagrams $(\Sigma,\alpha,\beta,L)$ and $(\Sigma',\alpha',\beta',L')$ are \emph{equivalent} if there is an orientation-preserving diffeomorphism $f:\Sigma\to \Sigma'$ such that $f(\alpha) = \alpha'$, $f(\beta) = \beta'$, and $f$ restricted to a tubular neighborhood of $L$ is a framed diffeomorphism.
\end{dff}

\begin{rmk}
    We pause to remark that a Heegaard-Link diagram is a generalization of an ordinary link diagram. An ordinary link diagram for a framed link $L$ may be considered as embedded on a sphere, with stabilizations resolving each crossing. This exactly specifies a Heegaard-Link diagram for the pair $(S^3,L)$. A general Heegaard-Link diagram is then a projection of the link on a Heegaard surface, with crossings marked (resolved with ``bridges''). 
    
    By an abuse of notation, we will use $L$ to denote the link in the original manifold and the link in a Heegaard-link diagram, as one would commonly do for links and link diagrams.
\end{rmk}

\begin{lem}\label{embedknotinsurface_lemma}
    Let $M$ be a closed connected orientable 3-manifold. Let $L$ be a framed link embedded in $M$. Then there exists a Heegaard splitting $M = H_u \cup H_l$ such that $L$ embeds in $\Sigma = H_u\cap H_l$ and the blackboard framing of $L$ agrees with the original framing.
\end{lem}

\begin{proof}
    Project the link to some Heegaard surface $\Sigma$ so that it is immersed and transverse to the attaching circles of the 1- and 2-handles and any self-intersections are at most transverse double points. At self-intersections, stabilize the surface by introducing a canceling 1-2 handle pair and push off one of the arcs across the ``bridge'' which is formed. Perform this in such a way that the projection is isotopic to the original link. The result is a Heegaard-Link diagram for the link $L$ treated as an unframed link.

    Now, let $(\Sigma,\alpha,\beta,L)$ be some (unframed) Heegaard-Link diagram representing the pair $(M,L)$. Since the set of framings on each component $L_i$ is a torsor over $\mbb Z$, we know that the blackboard framing of $L$ in $\Sigma$ differs from the endowed framing in the hypotheses by an integer $k_i$ for each component.

    We now perform the following procedure to ensure the two framings match, i.e. we wish to describe a way to take each $k_i$ to 0. For each $L_i$, stabilize the Heegaard diagram near $L_i$ and reroute it so that it traverses both the longitude of the stabilized torus once and the meridian  $|k_i|$ times. There are two ways to do this up to equivalence which correspond to taking $k_i$ to $k_i-|k_i|$ or taking $k_i$ to $k_i+|k_i|$. One of these will take $k_i$ to 0 depending on if $k_i$ is positive or negative. After doing this for each component, the result is an embedded link in $\Sigma$ whose blackboard framing matches the endowed framing.
\end{proof}

The final results of this section will be reviewing an algorithm to obtain a Heegaard diagram for the manifold obtained by surgery on $M$ along an embedded framed link $L$ and listing a complete set of moves for Heegaard-Link diagrams. Note that although the statement is for framed knots, we may extend it to links in the obvious way.

\begin{figure}
    \centering
    \begin{tikzpicture}
        \draw[darkgreen, ultra thick] (0,0) ellipse (1.5cm and 2cm);
        \draw[red, ultra thick] (-0.5,0) -- (-2.5,0);
        \draw[red, ultra thick] (-0.5,0.15) -- (-2.5,0.15);
        \draw[red, ultra thick] (-0.5,-0.15) -- (-2.5,-0.15);
        \draw[blue, ultra thick] (0.5,0) -- (2.5,0);
        \draw[blue, ultra thick] (0.5,0.15) -- (2.5,0.15);
        \draw[blue, ultra thick] (0.5,-0.15) -- (2.5,-0.15);
        \node at (1.5,-1.75) {$K$};
        \node at (-2,0.5) {$\alpha$};
        \node at (2,0.5) {$\beta$};
        \node at (4,0) {$\xrightarrow{\text{Surgery along }K}$};
    \begin{scope}[shift = {(8,0)}]
        \draw[blue, ultra thick] (0,0) ellipse (1.5cm and 2cm);
        \draw[red, ultra thick] (-0.5,0) -- (-2.5,0);
        \draw[red, ultra thick] (-0.5,0.15) -- (-2.5,0.15);
        \draw[red, ultra thick] (-0.5,-0.15) -- (-2.5,-0.15);
        \draw[blue, ultra thick] (0.5,0) -- (1,0);
        \draw[blue, ultra thick] (0.5,0.15) -- (1,0.15);
        \draw[blue, ultra thick] (0.5,-0.15) -- (1,-0.15);
        \draw[blue, ultra thick] (2,0) -- (2.5,0);
        \draw[blue, ultra thick] (2,0.15) -- (2.5,0.15);
        \draw[blue, ultra thick] (2,-0.15) -- (2.5,-0.15);
        \draw[draw=black, fill=lightgray] (1,0) circle (0.25);
        \draw[draw=black, fill=lightgray] (2,0) circle (0.25);
        \draw[red, ultra thick] (2,0) circle (0.25);
        \node at (1.5,-1.75) {$K_{\beta}$};
        \node at (-2,0.5) {$\alpha$};
        \node at (0.75,0.5) {$\beta$};
        \node at (2,0.5) {$K_\alpha$};
        \end{scope}
    \end{tikzpicture}
    \caption{A Heegaard Diagram for the surgery manifold $M(K)$}
    \label{fig:SurgeryHeegaardDiagram}
\end{figure}
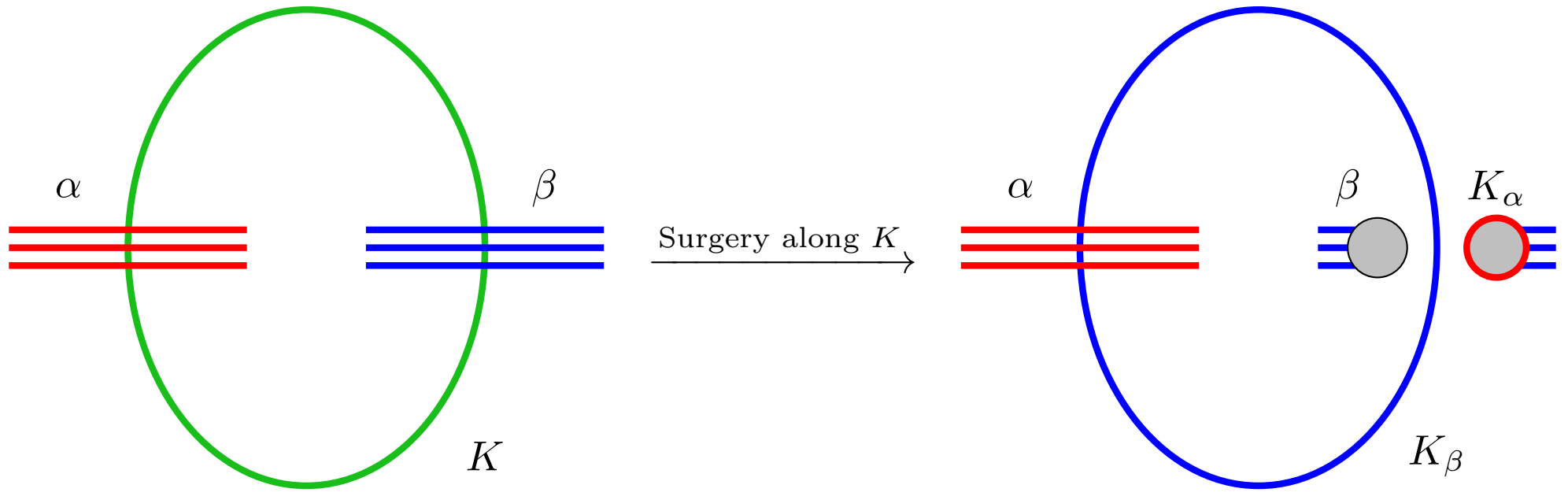

\begin{thm}[\cite{bridges20263man}]\label{thm:surgerydiagram}
    Let $K$ be a framed knot embedded in a 3-manifold $M$. Let $E=(\Sigma,\alpha,\beta,K)$ be a Heegaard-Link diagram for the pair $(M,K)$ such that the blackboard framing of $K$ in $\Sigma$ agrees with its endowed framing. Then there exists a genus $g(\Sigma)+1$ Heegaard diagram for the 3-manifold $M(K)$ obtained by surgery of $M$ along $K$ which is obtained from $E$ in the following way (see Figure \ref{fig:SurgeryHeegaardDiagram}):
    \begin{enumerate}[1)]
        \item Perform an isotopy of the curve sets $\alpha$ and $\beta$ so that there is a separation of $K$ into two curves $K_1$ and $K_2$ with $|\alpha\cap K_1| = 0$ and $|\beta\cap K_2| = 0$. 
        \item Choose a tubular neighborhood $N$ of $K_1$ in $\Sigma$ and label the components of $N\setminus K_1$ $N'$ and $N''$. Remove an open disk away from $\beta$ from both $N'$ and $N''$.
        \item Attach $I\times S^1$ to $\Sigma$ by identifying the two boundaries with the two boundaries introduced in the previous step.
        \item Add a new circle $K_{\alpha}$ to $\alpha$ and $K_{\beta}$ to $\beta$ such that $K_\beta$ replaces the knot $K$, and $K_\alpha$ is the belt of the attached $I\times S^1$.
        \item Reroute the $\beta$ curves which intersected $K$ across the attached $I\times S^1$ so that they do not traverse the meridian.
    \end{enumerate}
\end{thm}

\begin{thm}[Heegaard-Link Moves, \cite{bridges20263man}] \label{thm:extendedHeegaardMoves}
    Let $M$ be a closed connected oriented 3-manifold and let $L$ be an $m$-component framed link in $M$. Then any two Heegaard-Link diagrams for the pair $(M,L)$ are diffeomorphic after a finite sequence of the following moves:
    \begin{enumerate}[Move 1:]
        \item (Diagram Isotopy) An isotopy of the curve sets $\alpha$ and $\beta$, or an isotopy of the link $L$.
        \item (Handle Slide) A handle slide of the Heegaard curve sets away from the link or a slide of the link across a Heegaard circle which does not intersect any link components.
        \item (Stabilization/Destabilization) A stabilization or destabilization of the Heegaard diagram at a disk away from the link $L$.
    \end{enumerate}
\end{thm}

\subsection{Heegaard-Graph diagrams}

For a ribbon $\Bbbk$-category $\mcal D$, we say a $\mcal D$-colored ribbon graph $\Gamma$ is \emph{framed} if each arc is framed. We use the ribbon structure of $\mcal D$ to distinguish framings. We may embed a framed $\mcal D$-colored ribbon graph $\Gamma$ in a 3-manifold $M$. This embedding is well-defined because $\mcal \mcal D$ is ribbon. We say that two pairs $(M,\Gamma)$ and $(M',\Gamma')$ consisting of 3-manifolds and embedded framed $\mcal D$-colored ribbon graphs are \emph{diffeomorphic} if there is an orientation-preserving diffeomorphism $f:M\to M'$ such that its restriction to a tubular neighborhood of $\Gamma$ is a framed diffeomorphism to a tubular neighborhood of $\Gamma'$. We say they are \emph{equivalent} if there is a $\Gamma''$ skein equivalent to $\Gamma'$ in $M$ such that $(M,\Gamma)$ is diffeomorphic to $(M',\Gamma'')$.

As in Lemma \ref{embedknotinsurface_lemma}, we may embed $\Gamma$ in a Heegaard diagram $(\Sigma,\alpha,\beta)$ for $M$ such that its intersections with $\alpha$ and $\beta$ are transversal and occur only with the arcs of $\Gamma$. We can guarantee that $\Gamma$ has no self-intersections via stabilizations, but for our purposes we do not require this, as we may replace any self-intersections with half-braidings respecting the over--crossing information. 

\begin{dff}
    Let $\mcal D$ be a ribbon $\Bbbk$-category and $\mcal I$ a tensor ideal of $\mcal D$. We call a tuple $(\Sigma,\alpha,\beta,\Gamma)$ consisting of a Heegaard diagram $(\Sigma,\alpha,\beta)$ and $\Gamma$ an embedded framed $\mcal D$-colored $\mcal I$-admissible ribbon graph a \emph{Heegaard-Graph diagram}. We additionally require that a Heegaard-Graph diagram satisfies the condition that each $\alpha_i\in \alpha$ intersects with at least one $\beta_j$ or it intersects at least one arc in $\Gamma$ colored by an object in $\mcal I$.

    Two Heegaard-Graph diagrams $(\Sigma,\alpha,\beta,\Gamma)$ and $(\Sigma',\alpha',\beta',\Gamma')$ are (oriented) \emph{diffeomorphic} if there is an orientation-preserving diffeomorphism $f:\Sigma \to \Sigma'$ such that $f(\alpha) = \alpha'$, $f(\beta) = \beta'$, $f$ restricted to a tubular neighborhood of $\Gamma$ is a framed diffeomorphism. 

    We say $(\Sigma,\alpha,\beta,\Gamma)$ and $(\Sigma',\alpha',\beta',\Gamma')$ are \emph{equivalent} if there is a $\Gamma''$ $\mcal I$-skein equivalent to $\Gamma'$ in $\Sigma'$ such that $(\Sigma,\alpha,\beta,\Gamma)$ is diffeomorphic to $(\Sigma',\alpha',\beta',\Gamma'')$.
\end{dff}

As for Heegaard-Link diagrams, we have a complete collection of moves between equivalent Heegaard-Graph diagrams.

\begin{thm}\label{thm:Heegaard_Graph_Moves}
    Any two Heegaard-Graph diagrams representing equivalent pairs $(M,\Gamma)$ and $(M',\Gamma')$ are equivalent after a finite sequence of the following moves (restricted so that each indeed takes a Heegaard-Graph diagram to a Heegaard-Graph diagram):
    \begin{enumerate}[Move 1:]
        \item (Diagram Isotopy) An isotopy of the curve sets, or an isotopy of the graph.
        \item (Handle Slide) A handle slide of the Heegaard curve sets away from the graph or a slide of a framed arc of the graph across a Heegaard circle which does not intersect the graph anywhere.
        \item (Stabilization/Destabilization) A stabilization or destabilization of the Heegaard diagram at a disk away from the graph.
    \end{enumerate}
\end{thm}

\begin{proof}
    The proof of this theorem is analogous to the proof of Theorem \ref{thm:extendedHeegaardMoves}.
\end{proof}

\section{Invariants of Heegaard-Graph Diagrams and 3-Manifolds}
\label{GraphInvariant_Sec}

In this section we define an invariant of Heegaard-Graph diagrams using the chromatic structure of $\mcal C$. We then give a family of associated Heegaard-Graph diagrams to a given Heegaard-Link diagram, yielding an invariant of Heegaard-Link diagrams.

For embedded framed $\mcal D$-colored ribbon graphs in $S^3$, the trace defines an isotopy invariant $\mcal G:\mcal L \to \Bbbk$ if $\mcal D$ is spherical, where $\mcal L$ here is the set of framed $\mcal D$-colored ribbon graphs in $S^3$. As above, the modified trace $\mathrm{t}$ on $\mcal I$ in $\mcal D$ defines an isotopy invariant $\mcal G':\mcal L_{\mcal I}\to \Bbbk$ defined on framed $\mcal D$-colored $\mcal I$-admissible ribbon graphs in $S^3$, where $\mcal G'(\Gamma) = \mathrm{t}_X(\mcal G(T_X))$ where $T_X$ is a framed $\mcal D$-colored ribbon graph in $\mbb R^2\times [0,1]$ with one free bottom end and one free top end supported by framed arcs with no twists oriented downward and colored by $X\in \mcal I$.

\subsection{From Heegaard-Graph diagrams to ribbon graphs on handlebodies}

Let $\mcal D = \mcal Z(\mcal C)$ and let $\mcal I = \mathrm{Proj}(\mcal Z(\mcal C))$ be the tensor ideal of projective objects. Let $U:\mcal Z(\mcal C) \to \mcal C$ be the forgetful functor. We give a procedure to turn any Heegaard-Graph diagram $(\Sigma,\alpha,\beta,\Gamma)$ into an embedded graph in the boundary of the handlebody $(\Sigma,\alpha)$, called \emph{$\mcal C$-coloring $\beta$}. Let $G$ be a projective generator for $\mcal C$ and let $c_P:P\otimes G\to P\otimes G$ be a chromatic morphism based on $P$ for each projective $P$ in $\mcal C$. Suppose there is a curve $\gamma$ connecting an arc in $\Gamma$ colored by a projective object $M\in \mcal Z(\mcal C)$ to $\beta_i$. Then we may $\mcal C$-color each $\beta_i$ by first turning $\Gamma$ into an embedded $\mcal C$-colored ribbon graph in $\Sigma$ by applying $U$ to $\Gamma$, then performing

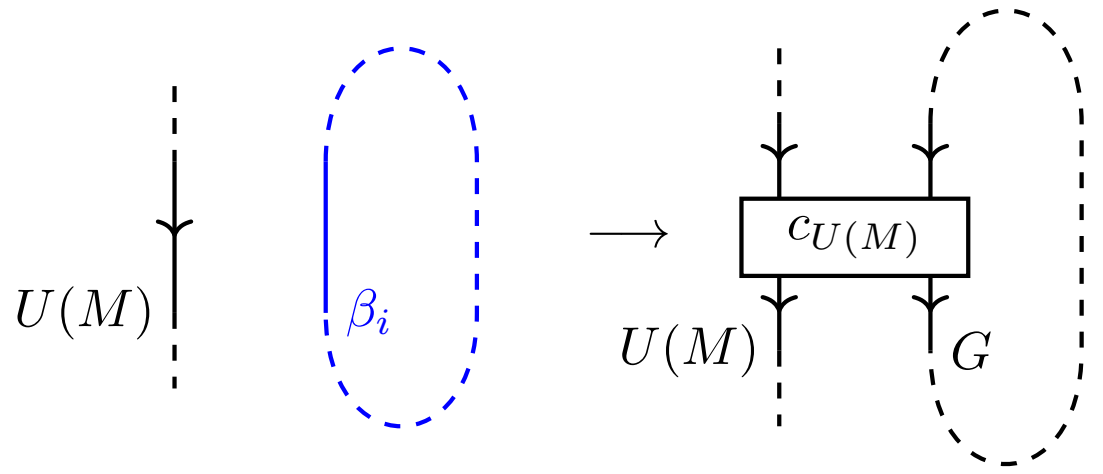
\begin{center}
    \begin{tikzpicture}
        \draw[mid arrow=>,thick] (0,0.5) -- (0,-0.5)node[left]{$U(M)$};
        \draw[dashed,thick] (0,1) -- (0,0.5);
        \draw[dashed,thick] (0,-0.5) -- (0,-1);

        \draw[thick, blue] (1,0.5) -- (1,-0.5)node[right]{$\beta_i$};
        \draw[thick,dashed,blue] (1,0.5) .. controls +(0,1) and +(0,1) .. (2,0.5) -- (2,-0.5) .. controls +(0,-1) and +(0,-1) .. (1,-0.5);

        \node at (3,0) {$\longrightarrow$};

        \begin{scope}[shift = {(4,0)}]
            \node[draw,thick,fill = white,minimum width = 1.5cm, minimum height = 0.5cm] at (0.5,0) {$c_{U(M)}$};
            \draw[mid arrow=>,thick] (0,0.75) -- (0,0.25);
            \draw[mid arrow=>,thick] (0,-0.25) -- (0,-0.75)node[left]{$U(M)$};
            \draw[thick,dashed] (0,1.25) -- (0,0.75);
            \draw[thick,dashed] (0,-0.75) -- (0,-1.25);

            \draw[mid arrow=>,thick] (1,0.75) -- (1,0.25);
            \draw[mid arrow=>,thick] (1,-0.25) -- (1,-0.75)node[right]{$G$};
            \draw[thick,dashed] (1,0.75) .. controls +(0,1) and +(0,1) .. (2,0.75) -- (2,-0.75) .. controls +(0,-1) and +(0,-1) .. (1,-0.75);
        \end{scope}
    \end{tikzpicture}
\end{center}

If $\beta_i$ intersects the graph $\Gamma$ in an arc colored by $M$, then replace the intersection with the half braiding $\sigma_G:M\otimes G\to G\otimes M$ or its inverse using the following prescription if the plane is given the standard orientation.

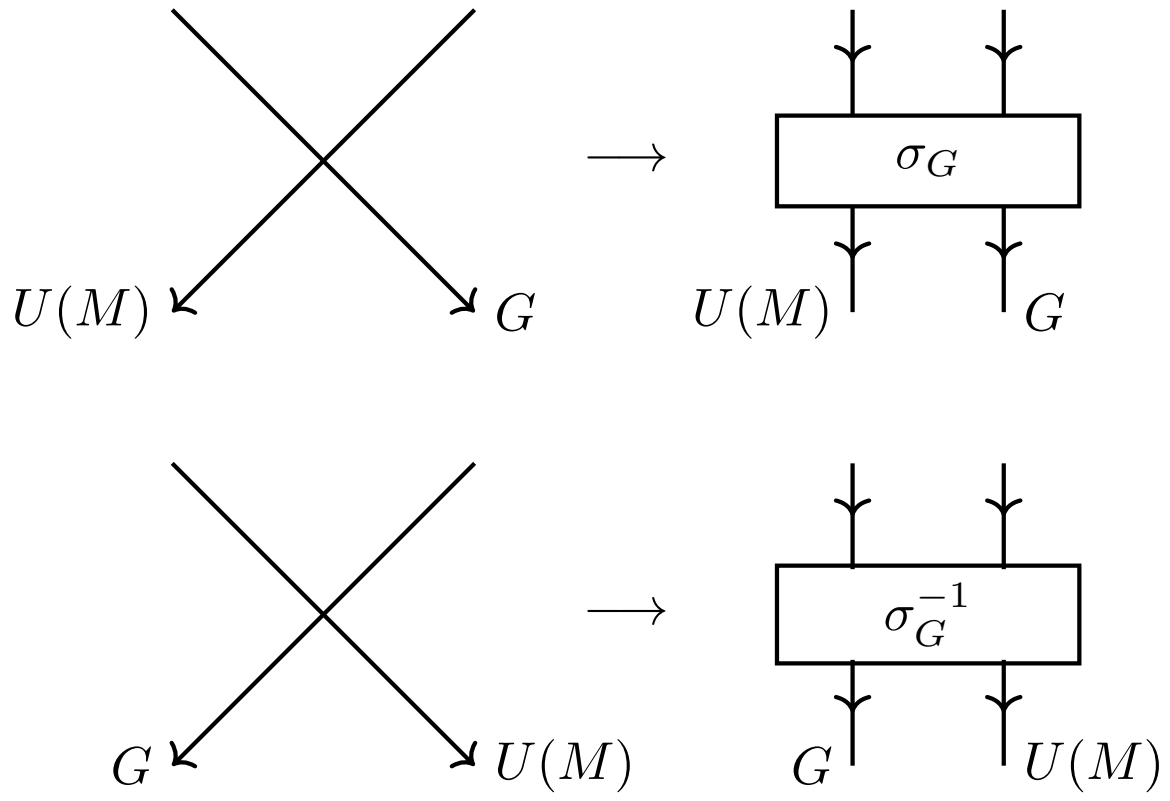
\begin{center}
    \begin{tikzpicture}
        \draw[thick,->] (-1,1) -- (1,-1)node[right]{$U(M)$};
        \draw[thick,->] (1,1) -- (-1,-1)node[left]{$G$};
        \node at (2,0) {$\longrightarrow$};
        \node[draw,thick,fill = white,minimum width = 2cm,minimum height = 0.6cm] at (4,0) {$\sigma_G^{-1}$};
        \draw[thick,mid arrow=>] (3.5,1) -- (3.5,0.3);
        \draw[thick,mid arrow=>] (4.5,1) -- (4.5,0.3);
        \draw[thick,mid arrow=>] (3.5,-0.3) -- (3.5,-1)node[left]{$G$};
        \draw[thick,mid arrow=>] (4.5,-0.3) -- (4.5,-1)node[right]{$U(M)$};

        \begin{scope}[shift = {(0,3)}]
            \draw[thick,->] (-1,1) -- (1,-1)node[right]{$G$};
        \draw[thick,->] (1,1) -- (-1,-1)node[left]{$U(M)$};
        \node at (2,0) {$\longrightarrow$};
        \node[draw,thick,fill = white,minimum width = 2cm,minimum height = 0.6cm] at (4,0) {$\sigma_G$};
        \draw[thick,mid arrow=>] (3.5,1) -- (3.5,0.3);
        \draw[thick,mid arrow=>] (4.5,1) -- (4.5,0.3);
        \draw[thick,mid arrow=>] (3.5,-0.3) -- (3.5,-1)node[left]{$U(M)$};
        \draw[thick,mid arrow=>] (4.5,-0.3) -- (4.5,-1)node[right]{$G$};
        \end{scope}
    \end{tikzpicture}
\end{center}

This procedure is well-defined and does not depend on the choice of generator $G$ or the choice of chromatic morphisms $c_P$ (see \cite{Geer2020kuperberg} for details). Denote the resulting $\mcal C$-colored ribbon graph by $\widetilde{\Gamma}$.

Finally, we may turn this graph in $\Sigma$ into a $\mcal C$-colored $\mathrm{Proj}(\mcal C)$-admissible ribbon graph in $S^2$, which, together with the modified invariant $\mcal F'$, yields a well-defined map from the set of Heegaard-Graph diagrams to $\Bbbk$.

To such a ribbon graph in the boundary of the handlebody $(\Sigma,\alpha)$, we perform the following. First arbitrarily orient each $\alpha_i$ and endow each with a base point. From the base point of each $\alpha_i$, record the signed intersections with arcs in the ribbon graph $\widetilde{\Gamma}$ and form the object given by the ordered tensor product of the objects coloring each arc, signed according to the parity of the intersection. If the crossing $p$ of $\alpha_i$ with an arc $a$ colored by the object $X$ is positive with respect to the orientation of $\Sigma$, meaning the ordered pair $(\mathrm{d}_p a, \mathrm{d}_p \alpha)$ forms a positively oriented basis for the tangent space of $\Sigma$ at $p$, then this crossing contributes $X$ to the tensor product, and if it is negative, then this crossing contributes $X^*$. After this, we obtain an associated object $X_{\alpha_i} = X_1 \otimes X_2 \otimes \cdots \otimes X_c$, where $c$ is the number of crossings with $\alpha_i$. Then we ``cut'' along $\alpha_i$ by replacing it with $\Lambda_{X_{\alpha_i}}$, where the coupon is read with the orientation of $\alpha_i$. Note that since $\alpha_i$ always has at least one intersection, this is defined. Graphically:

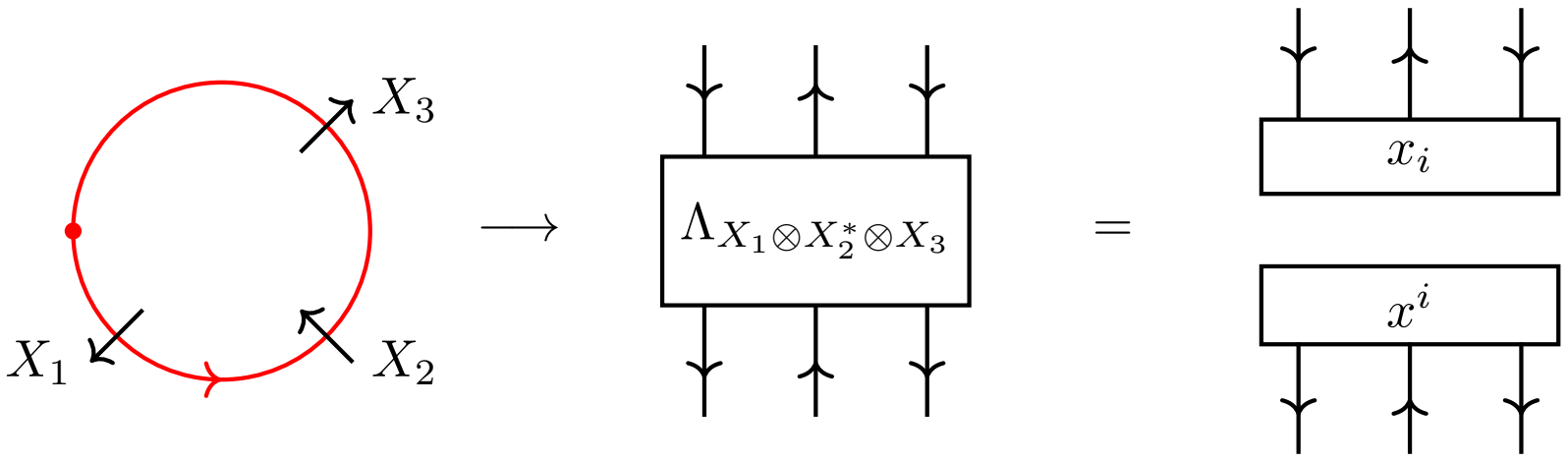
\begin{center}
    \begin{tikzpicture}
        \draw[red,thick,decoration={markings, mark=at position 0.75 with {\arrow{>}}}, postaction={decorate}] (0,0) circle (1cm);
        \draw[red,fill = red] (-1,0) circle (0.05cm);

        \foreach \x in {-135,45} {
            \draw[thick,->] (\x:0.75) -- (\x:1.25);
        }
        \draw[thick,<-] (-45:0.75) -- (-45:1.25);
        \node[left] at (-135:1.25) {$X_1$};
        \node[right] at (-45:1.25) {$X_2$};
        \node[right] at (45:1.25) {$X_3$};

        \node at (2,0) {$\longrightarrow$};

        \begin{scope}[shift = {(3,0)}]
            \node[draw,thick,fill=white,minimum width = 2cm,minimum height = 1cm] at (1,0) {$\Lambda_{X_1\otimes X_2^*\otimes X_3}$};
            \draw[mid arrow=>,thick] (0.25,1.25) -- (0.25,0.5);
            \draw[mid arrow=<,thick] (1,1.25) -- (1,0.5);
            \draw[mid arrow=>,thick] (1.75,1.25) -- (1.75,0.5);
            \draw[mid arrow=<,thick] (0.25,-1.25) -- (0.25,-0.5);
            \draw[mid arrow=>,thick] (1,-1.25) -- (1,-0.5);
            \draw[mid arrow=<,thick] (1.75,-1.25) -- (1.75,-0.5);

            \node at (3,0) {$=$};
        \end{scope}

        \begin{scope}[shift = {(7,0)}]
            \node[draw,thick,fill=white,minimum width = 2cm,minimum height = 0.5cm] at (1,0.5) {$x_i$};
            \node[draw,thick,fill=white,minimum width = 2cm,minimum height = 0.5cm] at (1,-0.5) {$x^i$};
            \draw[mid arrow=>,thick] (0.25,1.5) -- (0.25,0.75);
            \draw[mid arrow=<,thick] (1,1.5) -- (1,0.75);
            \draw[mid arrow=>,thick] (1.75,1.5) -- (1.75,0.75);
            \draw[mid arrow=<,thick] (0.25,-1.5) -- (0.25,-0.75);
            \draw[mid arrow=>,thick] (1,-1.5) -- (1,-0.75);
            \draw[mid arrow=<,thick] (1.75,-1.5) -- (1.75,-0.75);
        \end{scope}
    \end{tikzpicture}
\end{center}

 We call this procedure \emph{cutting along $\alpha_i$}. Then since $\Sigma - \alpha$ is a sphere with $2g$ boundary components, this new ribbon graph embeds in a sphere. Again, this procedure is well-defined as in it does not depend on the orientation or choice of base point of each $\alpha_i$ thanks to Lemma \ref{lem:cuttingmorphismproperties}.

\begin{thm}
    Let $\mcal C$ be a chromatic category. Then the mapping $\mcal F'$ extends uniquely to Heegaard-Graph diagrams in a way that is preserved under $\mcal C$-coloring $\beta$ and cutting along $\alpha$.
\end{thm}

\begin{proof}
    This follows from Theorem 2.3 of \cite{Geer2020kuperberg}.
\end{proof}

\begin{xmp}
    Consider the following Heegaard-Graph diagram where the graph $\Gamma$ defined from $f:U(M)\to U(M)$ is colored green.

    
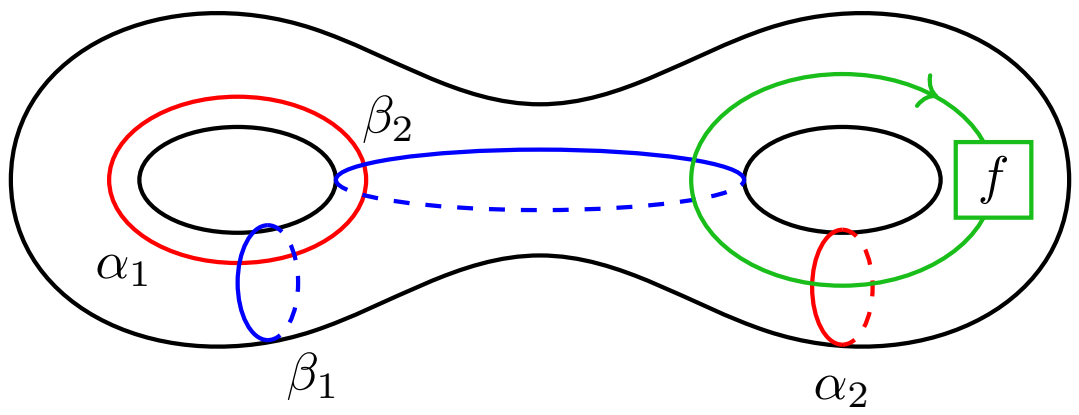
\begin{center}
        \begin{tikzpicture}
        \draw[closed, thick] (0,0) to[curve through = {(2,-1) .. (3.5,-0.5) .. (5,-1) .. (7,0) .. (5,1) .. (3.5,0.5) .. (2,1)}] (0,0);
        \draw[thick] (1.5,0) ellipse (0.65 and 0.35);
        \draw[thick] (5.5,0) ellipse (0.65 and 0.35);
        \draw[red, thick] (1.5,0) ellipse (0.85 and 0.55);
        \node at (0.75,-0.6) {$\alpha_1$};
        \draw[blue, thick] (1.7,-0.3) arc(90:270:0.2 and 0.38);
        \draw[blue, thick, dashed] (1.7,-0.3) arc(90:-90:0.2 and 0.38);
        \node at (2,-1.3) {$\beta_1$};
        \draw[red, thick] (5.5,-0.33) arc(90:270:0.2 and 0.38);
        \draw[red, thick, dashed] (5.5,-0.33) arc(90:-90:0.2 and 0.38);
        \node at (5.5,-1.4) {$\alpha_2$};
        \draw[blue, thick] (2.15,0) arc(180:0:1.35 and 0.2);
        \draw[blue, thick, dashed] (2.15,0) arc(-180:0:1.35 and 0.2);
        \node at (2.5,0.4) {$\beta_2$};

        \draw[darkgreen, thick, decoration={markings, mark=at position 0.15 with {\arrow{<}}}, postaction={decorate}] (5.5,0) ellipse (1 and 0.7);
        \node[draw, thick, darkgreen, fill = white, minimum width = 0.5cm, minimum height = 0.5cm] at (6.5,0) {~};
        \node at (6.5,0) {$f$};
    \end{tikzpicture}
    \end{center}

    Then $\mcal C$-coloring $\beta$ yields the following ribbon graph on $(\Sigma,\alpha)$ this time colored blue and embedded in the plane with gray circles indicating genus.

    
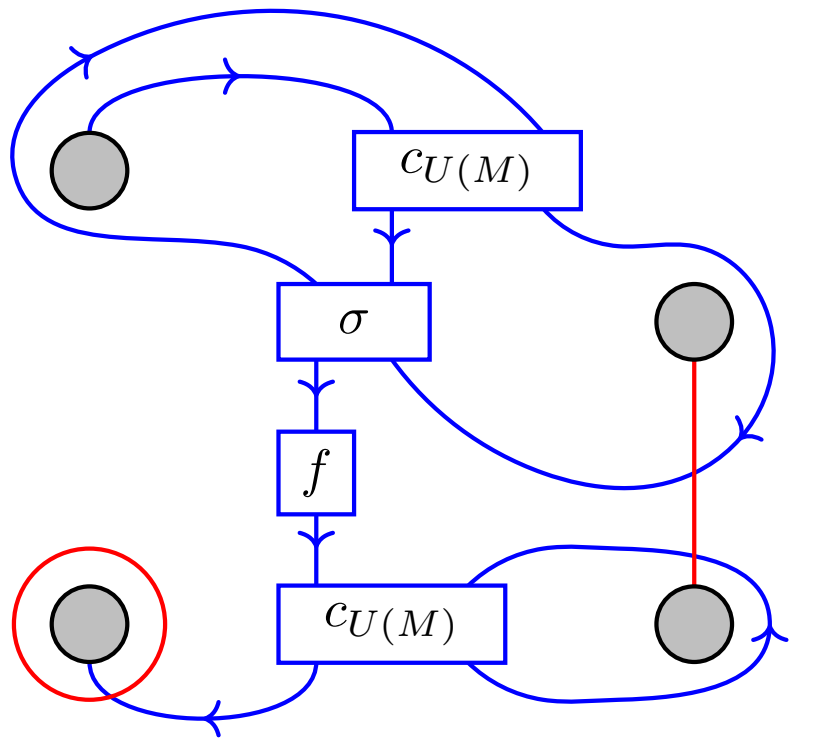
\begin{center}
        \begin{tikzpicture}
            \draw[blue,thick,mid arrow=>] (-0.5,0.75) -- (-0.5,0.25);
            \draw[blue,thick,mid arrow=>] (-0.5,-0.25) -- (-0.5,-0.75);
            \draw[blue,thick,mid arrow=>] (-0.5,-1.25) .. controls +(0,-0.5) and +(0,-0.5) .. (-2,-1.25);
            \draw[blue,thick,mid arrow=>] (-2,2.25) .. controls +(0,0.5) and +(0,0.5) .. (0,2.25);
            \draw[blue,thick,mid arrow=>] (0,1.75) -- (0,1.25);
            \draw[blue,thick,mid arrow=>] (-0.5,1.25) to[curve through={(-1,1.5) (-2.5,2) (-2,2.75) (0.75,2.5)}] (1,2.25);
            \draw[blue,thick,mid arrow=>] (1,1.75) to[curve through={(1.5,1.5) (2,1.5) (2.5,1) (2,0) (1,0)}] (0,0.75);
            \draw[blue,thick,mid arrow=>] (0.5,-1.25) to[curve through = {(1,-1.5) (1.5,-1.5) (2.5,-1) (1.5,-0.5) (1,-0.5)}] (0.5,-0.75);

            \draw[red,thick] (-2,-1) circle (0.5cm);
            \draw[red,thick] (2,0.75) -- (2,-0.75);
            
            \draw[black,thick,fill = lightgray] (-2,-1) circle (0.25cm);   
            \draw[black,thick,fill = lightgray] (-2,2) circle (0.25cm); 
            \draw[black,thick,fill = lightgray] (2,-1) circle (0.25cm);   
            \draw[black,thick,fill = lightgray] (2,1) circle (0.25cm);
            
            \node[draw=blue,thick,fill = white,minimum width = 1.5cm, minimum height = 0.5cm] at (0,-1) {$c_{U(M)}$};
            \node[draw=blue,thick,fill=white, minimum width = 0.5cm, minimum height = 0.5cm] at (-0.5,0) {$f$};
            \node[draw=blue,thick,fill = white, minimum width =1cm, minimum height = 0.5cm] at (-0.25,1) {$\sigma$};
            \node[draw=blue,thick,fill = white,minimum width = 1.5cm, minimum height = 0.5cm] at (0.5,2) {$c_{U(M)}$};
        \end{tikzpicture}
    \end{center}

    Then cutting along $\alpha$ yields the diagram in the plane

    
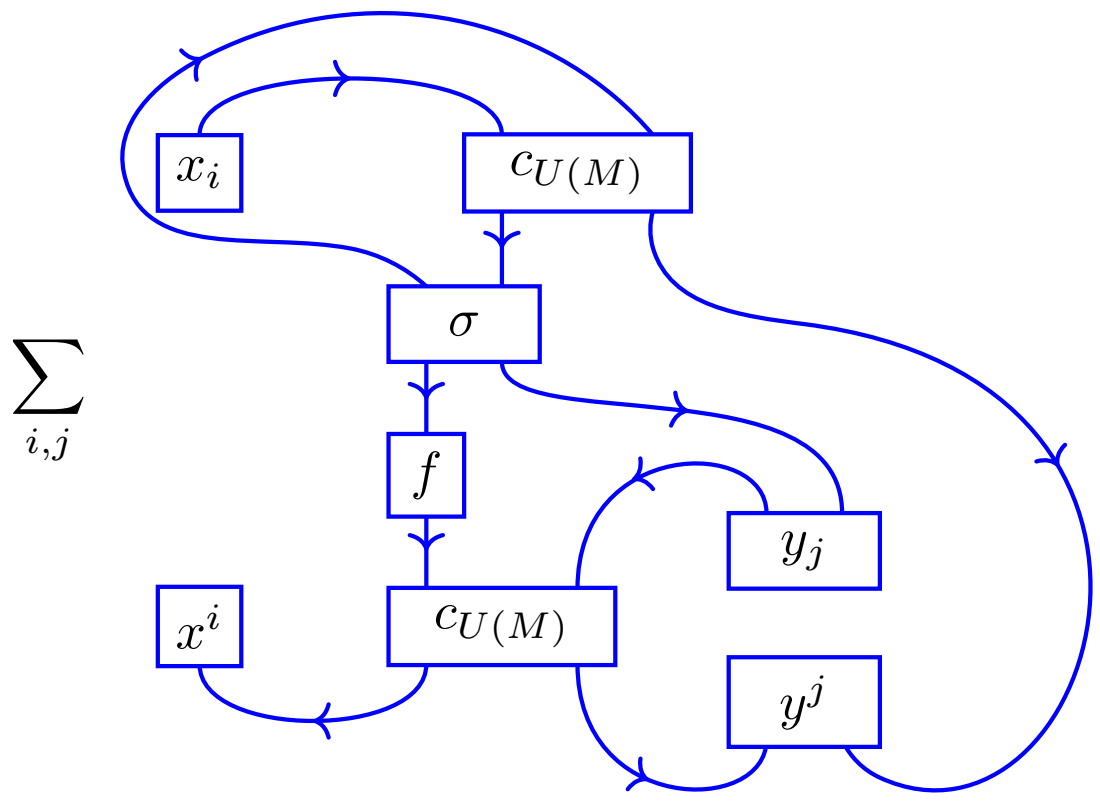
\begin{center}
        \begin{tikzpicture}
            \draw[blue,thick,mid arrow=>] (-0.5,0.75) -- (-0.5,0.25);
            \draw[blue,thick,mid arrow=>] (-0.5,-0.25) -- (-0.5,-0.75);
            \draw[blue,thick,mid arrow=>] (-0.5,-1.25) .. controls +(0,-0.5) and +(0,-0.5) .. (-2,-1.25);
            \draw[blue,thick,mid arrow=>] (-2,2.25) .. controls +(0,0.5) and +(0,0.5) .. (0,2.25);
            \draw[blue,thick,mid arrow=>] (0,1.75) -- (0,1.25);
            \draw[blue,thick,mid arrow=>] (-0.5,1.25) to[curve through={(-1,1.5) (-2.5,2) (-2,2.75) (0.75,2.5)}] (1,2.25);
            \draw[blue,thick,mid arrow=>] (1,1.75) to[curve through={(1,1.5) (2,1) (2.5,-2)}] (2.25,-1.75);
            \draw[blue,thick,mid arrow=>] (0,0.75) .. controls +(0,-0.5) and +(0,1) .. (2.25,-0.25);
            \draw[blue,thick,mid arrow=<] (0.5,-0.75) .. controls +(0,1) and +(0,0.5) .. (1.75,-0.25);
            \draw[blue,thick,mid arrow=>] (0.5,-1.25) .. controls +(0,-1) and +(0,-0.5) .. (1.75,-1.75);
            
            \node at (-3,0.5) {$\displaystyle\sum_{i,j}$};
            
            \node[draw=blue,thick,fill = white,minimum width = 0.5cm, minimum height =0.5cm] at (-2,-1) {$x^i$};   
            \node[draw=blue,thick,fill = white,minimum width = 0.5cm, minimum height =0.5cm] at (-2,2) {$x_i$};
            \node[draw=blue,thick,fill = white,minimum width = 1cm, minimum height =0.5cm] at (2,-1.5) {$y^j$};   
            \node[draw=blue,thick,fill = white,minimum width = 1cm, minimum height =0.5cm] at (2,-0.5) {$y_j$};
            
            \node[draw=blue,thick,fill = white,minimum width = 1.5cm, minimum height = 0.5cm] at (0,-1) {$c_{U(M)}$};
            \node[draw=blue,thick,fill=white, minimum width = 0.5cm, minimum height = 0.5cm] at (-0.5,0) {$f$};
            \node[draw=blue,thick,fill = white, minimum width =1cm, minimum height = 0.5cm] at (-0.25,1) {$\sigma$};
            \node[draw=blue,thick,fill = white,minimum width = 1.5cm, minimum height = 0.5cm] at (0.5,2) {$c_{U(M)}$};
        \end{tikzpicture}
    \end{center}
\end{xmp}

\begin{thm}\label{thm:HeegaardGraphInvariant}
    The extension of $\mcal F'$ to Heegaard-Graph diagrams is invariant under the Heegaard-Graph moves, and hence defines an invariant of pairs $(M,\Gamma)$ where $M$ is a closed connected oriented 3-manifold and $\Gamma$ is a framed $\mcal D$-colored $\mcal I$-admissible ribbon graph in $M$. Denote this invariant by \[Z_{\mathrm{ExCGPT}}(M,\Gamma; \mcal C)\]
\end{thm}

\begin{proof}
    We show that $\mcal F'$ is invariant under each move.

    For a diagram isotopy, the following will be useful:

    \begin{lem}
        Let $\Gamma$ be a closed 1-manifold in a surface $\Sigma$. Let $i_0$ and $i_1$ be homotopic immersions of $\Gamma$ into $\Sigma$ such that each component of $\Gamma$ is embedded in $\Sigma$ by $i_0$ and $i_1$, and the components of $i_0(\Gamma)$ and $i_1(\Gamma)$ intersect transversely, then $i_0$ and $i_1$ are diffeomorphic after a sequence of two-point and three-point moves (see Figures \ref{fig:twopoint} and \ref{fig:threepoint}).
    \end{lem}

    \begin{proof}
        By the bigon criterion (see \cite[Lemma~3.1]{Hass85}), after a sequence of two-point moves, $i_0$ and $i_1$ may be assumed to be in minimal position. Then by \cite[Lemma~3.4]{PATERSON2002205} any two minimal immersions are ambiently isotopic after a sequence of three-point moves.
    \end{proof}

    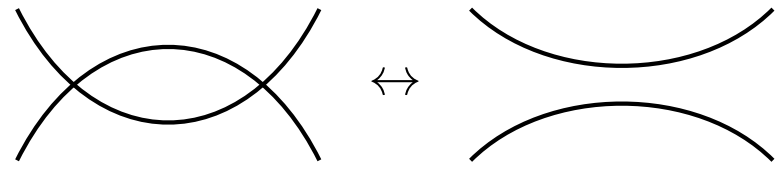
\begin{figure}
        \centering
        \begin{tikzpicture}
            \draw[thick] (-1,0.5) .. controls +(0.5,-1) and +(-0.5,-1) .. (1,0.5);
            \draw[thick] (-1,-0.5) .. controls +(0.5,1) and +(-0.5,1) .. (1,-0.5);
            \node[thick] at (1.5,0) {$\leftrightarrow$};
            \draw[thick] (2,0.5) .. controls +(0.5,-0.5) and +(-0.5,-0.5) .. (4,0.5);
            \draw[thick] (2,-0.5) .. controls +(0.5,0.5) and +(-0.5,0.5) .. (4,-0.5);
        \end{tikzpicture}
        \caption{Two-Point Move}
        \label{fig:twopoint}
    \end{figure}

    \begin{figure}
        \centering
        \begin{tikzpicture}
            \draw[thick] (-1.25,-0.5) -- (1.25,-0.5);
            \draw[thick] (-0.5,1) -- (1,-1);
            \draw[thick] (0.5,1) -- (-1,-1);
            \node at (2,0) {$\leftrightarrow$};
            \begin{scope}[shift = {(4,0)}]
            \draw[thick] (-1.25,-0.5) .. controls +(1,1) and +(-1,1) .. (1.25,-0.5);
            \draw[thick] (-0.5,1) .. controls +(0,-1) and +(-1,0.25) .. (1,-1);
            \draw[thick] (0.5,1) .. controls +(0,-1) and +(1,0.25) .. (-1,-1);
            \end{scope}
        \end{tikzpicture}
        \caption{Three-Point Move}
        \label{fig:threepoint}
    \end{figure}
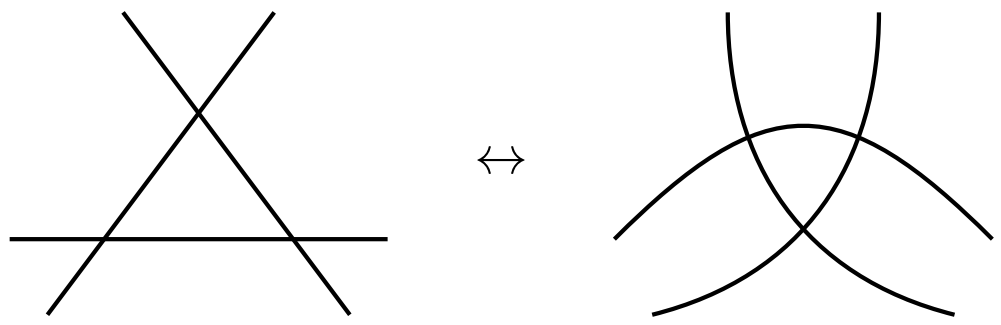

    Thus we need only to demonstrate invariance under two-point and three-point moves. We treat the graph $\Gamma$ as a stratified 1-manifold where its framed arcs are edges and its coupons are vertices.

    A two-point move between $\Gamma$ and a $\beta$ curve introduces $\sigma_G \sigma_G^{-1}$ and hence is invariant.

    A two-point move between $\Gamma$ and an $\alpha$ curve yields \[(\Id_{X_1}\otimes \cdots \Id_{X_k} \otimes \rev_{U(M)} \otimes \Id_{X_{k+1}} \otimes \cdots \otimes \Id_{X_c}) \Lambda_{X_1 \otimes \cdots \otimes X_k \otimes \Id_{U(M)^*\otimes U(M)}\otimes X_{k+1}\otimes  \cdots \otimes X_c}\] or a similar expression. Then since $c\neq 0$, and one of $X_i$ is always projective, we may use naturality of $\Lambda$ to obtain invariance.

    For a two-point move between an $\alpha$ and $\beta$ curve, we have the same as for a two-point move between $\Gamma$ and $\alpha$ except replace $U(M)$ with $G$.

    A three-point move among all three curves follows from naturality of $\Lambda$, moving the half-braiding across. This proves invariance under diagram isotopy (Move 1).

    For a handle slide of $\Gamma$ or $\beta$ over a $\beta$ curve, this follows from \cite[Lemma 2.8]{Geer2026skein3plus1}.

    For a handle slide of $\alpha$ across an $\alpha$ curve, we have the following:

    
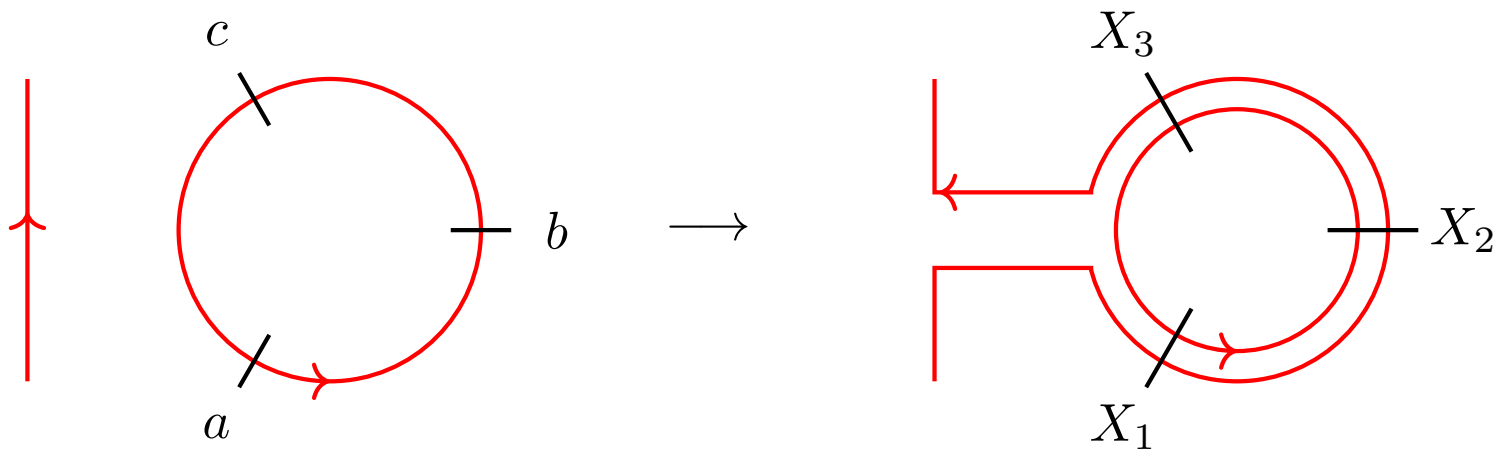
\begin{center}
    \begin{tikzpicture}
        \draw[red,thick,mid arrow=<] (0,1) -- (0,-1);
        \begin{scope}[shift = {(2,0)}]
            \draw[red,thick,decoration={markings, mark=at position 0.75 with {\arrow{>}}}, postaction={decorate}] (0,0) circle (1cm);
            \foreach \angle in {0,-120,120} {
            \draw[thick] (\angle : 0.8) -- (\angle : 1.2);
        }
        \node at (0 : 1.5) {$b$};
        \node at (120 : 1.5) {$c$};
        \node at (-1200 : 1.5) {$a$};
        \end{scope}
        \node at (4.5,0) {$\longrightarrow$};
        \begin{scope}[shift = {(6,0)}]
        \begin{scope}[shift = {(2,0)}]
        \draw[red,thick,decoration={markings, mark=at position 0.75 with {\arrow{>}}}, postaction={decorate}] (0,0) circle (0.8cm);
        \draw[red,thick] (0,0) circle (1cm);
        \foreach \angle in {0,-120,120} {
            \draw[thick] (\angle : 0.6) -- (\angle : 1.2);
        }
        \node at (0 : 1.5) {$X_2$};
        \node at (120 : 1.5) {$X_3$};
        \node at (-1200 : 1.5) {$X_1$};
        \draw[fill=white,draw=white] (-1.2,-0.25) rectangle (-0.9,0.25);
        \draw[red,thick] (-2,-1) -- (-2,-0.25) -- (-0.95,-0.25);
        \draw[red,thick,mid arrow=<] (-2,1) -- (-2,0.25) -- (-0.95,0.25);
        \end{scope}
        \end{scope}
    \end{tikzpicture}
    \end{center}

    Then the after picture corresponds to 

    
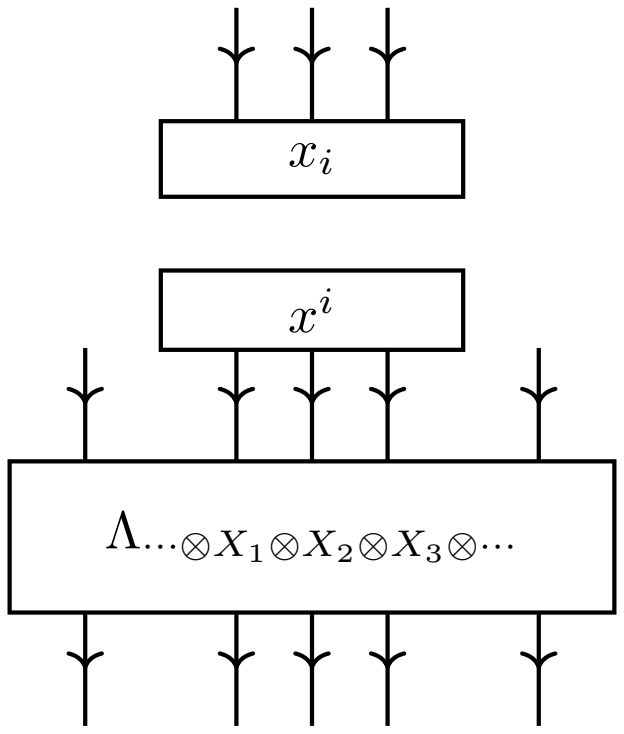
\begin{center}
        \begin{tikzpicture}
            \node[draw,thick,fill=white,minimum width = 2cm,minimum height = 0.5cm] at (0,0.5) {$x_i$};
            \node[draw,thick,fill=white,minimum width = 2cm,minimum height = 0.5cm] at (0,-0.5) {$x^i$};
            \node[draw,thick,fill=white,minimum width = 4cm,minimum height = 1cm] at (0,-2) {$\Lambda_{\cdots \otimes X_1\otimes X_2\otimes X_3\otimes \cdots}$};
            \draw[thick,mid arrow=>] (-0.5,1.5) -- (-0.5,0.75);
            \draw[thick,mid arrow=>] (0,1.5) -- (0,0.75);
            \draw[thick,mid arrow=>] (0.5,1.5) -- (0.5,0.75);
            \draw[thick,mid arrow=>] (-0.5,-0.75) -- (-0.5,-1.5);
            \draw[thick,mid arrow=>] (0,-0.75) -- (0,-1.5);
            \draw[thick,mid arrow=>] (0.5,-0.75) -- (0.5,-1.5);
            \draw[thick,mid arrow=>] (-1.5,-0.75) -- (-1.5,-1.5);
            \draw[thick,mid arrow=>] (1.5,-0.75) -- (1.5,-1.5);
            \begin{scope}[shift = {(0,-1.75)}]
                \draw[thick,mid arrow=>] (-0.5,-0.75) -- (-0.5,-1.5);
            \draw[thick,mid arrow=>] (0,-0.75) -- (0,-1.5);
            \draw[thick,mid arrow=>] (0.5,-0.75) -- (0.5,-1.5);
            \draw[thick,mid arrow=>] (-1.5,-0.75) -- (-1.5,-1.5);
            \draw[thick,mid arrow=>] (1.5,-0.75) -- (1.5,-1.5);
            \end{scope}
        \end{tikzpicture}
    \end{center}

    Sliding $x^i$ through the lower $\Lambda$ yields invariance. Hence, we have invariance under handle slides (Move 2).

    Finally, stabilization is given by the chromatic identity.
\end{proof}

\subsection{Invariants of Heegaard-Link diagrams}

From a Heegaard-Link diagram $(\Sigma,\alpha,\beta,L)$, we obtain a Heegaard-Graph diagram via the following procedure. Recall that $L(G)$ is a projective generator for $\mcal Z(\mcal C)$ and let $c^{\mcal Z}$ be a chromatic morphism for $L(G)$ in $\mcal Z(\mcal C)$. Recall the non-degenerate modified trace $\ztr$ for $\mcal Z(\mcal C)$ from Section \ref{Center_Sec}. Then since $\ztr$ is non-degenerate, there is an endomorphism $h:L(G)\to L(G)$ such that $\ztr_{L(G)}(h) =1$. Then embed the ribbon graph $O_{L(G)} = \mathrm{tr}(\pi_{G,\mbb 1}U(\Omega_G)U(h)i_{G,\mbb 1})$ in a disk in $\Sigma$ away from $\alpha$, $\beta$, and $L$, and ``color the link'' with the chromatic morphism by

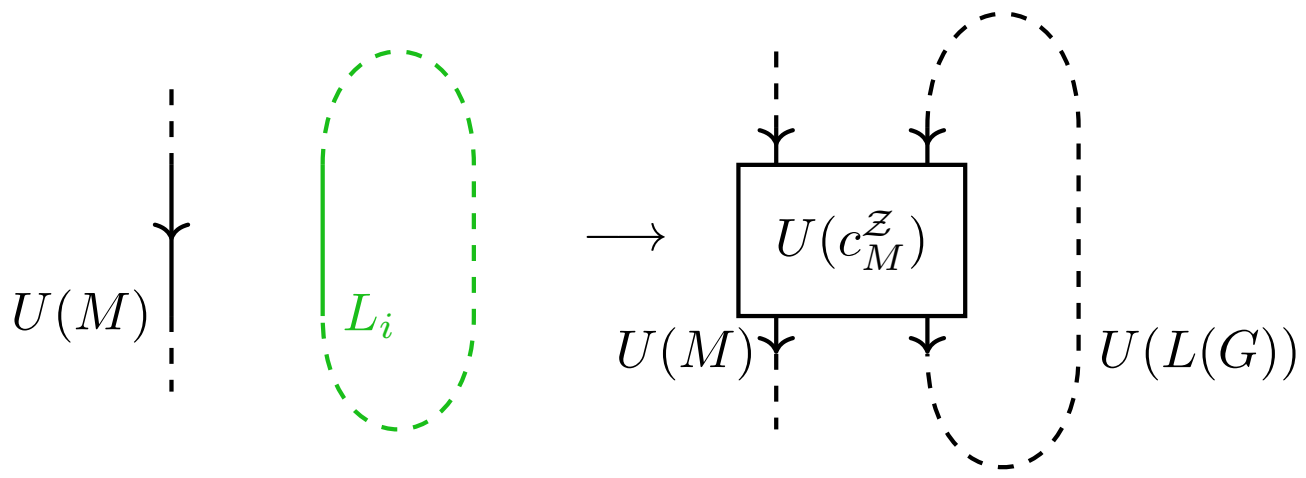
\begin{center}
    \begin{tikzpicture}
        \draw[mid arrow=>,thick] (0,0.5) -- (0,-0.5)node[left]{$U(M)$};
        \draw[dashed,thick] (0,1) -- (0,0.5);
        \draw[dashed,thick] (0,-0.5) -- (0,-1);

        \draw[thick, darkgreen] (1,0.5) -- (1,-0.5)node[right]{$L_i$};
        \draw[thick,dashed,darkgreen] (1,0.5) .. controls +(0,1) and +(0,1) .. (2,0.5) -- (2,-0.5) .. controls +(0,-1) and +(0,-1) .. (1,-0.5);

        \node at (3,0) {$\longrightarrow$};

        \begin{scope}[shift = {(4,0)}]
            \node[draw,thick,fill = white,minimum width = 1.5cm, minimum height = 1cm] at (0.5,0) {$U(c^{\mcal Z}_{M})$};
            \draw[mid arrow=>,thick] (0,0.75) -- (0,0.5);
            \draw[->,thick] (0,-0.5) -- (0,-0.75)node[left]{$U(M)$};
            \draw[thick,dashed] (0,1.25) -- (0,0.75);
            \draw[thick,dashed] (0,-0.75) -- (0,-1.25);

            \draw[mid arrow=>,thick] (1,0.75) -- (1,0.5);
            \draw[->,thick] (1,-0.5) -- (1,-0.75);
            \draw[thick,dashed] (1,0.75) .. controls +(0,1) and +(0,1) .. (2,0.75) -- (2,-0.75)node[right]{$U(L(G))$} .. controls +(0,-1) and +(0,-1) .. (1,-0.75);
        \end{scope}
    \end{tikzpicture}
\end{center}

The result is a Heegaard-Graph diagram $(\Sigma,\alpha,\beta,\widetilde{L})$, where $\widetilde{L}$ is the ribbon graph obtained by using $O_{L(G)}$ to color the link $L$. 

\begin{thm}
    The mapping $\mcal F'$ induces an invariant of pairs $(M,L)$ consisting of a closed connected oriented 3-manifold $M$ with an embedded framed link $L$. We denote this invariant by \[Z_{\mathrm{ExCGPT}}(M,L;\mcal C)\]
\end{thm}

\begin{proof}
    If we show that embedding the ribbon graph $O_{L(G)}$ is well-defined, then Theorem \ref{thm:HeegaardGraphInvariant} together with Theorem \ref{thm:extendedHeegaardMoves} gives the result. This follows because $\Sigma -\beta$ is connected, so we may use naturality of cutting along $\alpha$ and the centrality of the colored link $L$ to ensure that $O_{L(G)}$ may be placed anywhere.
\end{proof}

\begin{cor}\label{cor:recoverCGPT}
    When $L$ is empty, then $Z_{\mathrm{ExCGPT}}(M,\emptyset;\mcal C)$ recovers the 3-manifold invariant in \cite{Geer2020kuperberg}, which we denote by $Z_{\mathrm{CGPT}}(M;\mcal C)$.
\end{cor}

\begin{proof}
    We only need to show that $\ctr_G(\pi_{G,\mbb 1}U(\Omega_G)U(h)i_{G,\mbb 1}) = 1$. This follows from the definition of $\ztr$.
\end{proof}

\section{Relationship with the Renormalized Reshetikhin-Turaev Invariant}
\label{RT_Sec}

In the previous section, we showed that the renormalized graphical calculus functor $\mcal F'$ extends uniquely to an invariant of embedded framed ribbon graphs and framed links in a 3-manifold $M$. Here we will show that when $M = S^3$, the invariant of Heegaard-Link diagrams recovers the renormalized Reshetikhin-Turaev invariant for $M = S^3(L)$.

Let $\mcal D$ be a twist nondegenerate ribbon chromatic category (see \cite{Geer2026skein3plus1}). Then $\mcal D$ defines a 3-manifold invariant $B_{\mcal D}(M) = \frac{\mcal F'(\widetilde{L})}{\Delta^r_+ \Delta^s_-}$ called the \emph{renormalized Reshetikhin-Turaev Invariant} defined by using a chromatic morphism in $\mcal D$ to color a surgery link for $M$ analogously to our construction in Section \ref{GraphInvariant_Sec}, then using the modified Reshetikhin-Turaev ribbon functor $\mcal F'$ to obtain a scalar. Here, $(r,s)$ is the signature of the linking matrix for $L$, and the scalars $\Delta_+$ and $\Delta_-$ are defined as follows.

Let $\epsilon:P_{\mbb 1}\to \mbb 1$ be the projective cover for the unit $\mbb 1\in \mcal D$ and let $\eta:\mbb 1\to P_{\mbb 1}$ be the socle. Then $\epsilon$ and $\eta$ are bases for $\Hom_{\mcal D}(P_{\mbb 1}, \mbb 1)$ and $\Hom_{\mcal D}(\mbb 1, P_{\mbb 1})$, respectively. For $G$ a generator for $\mcal D$ and $c$ a chromatic morphism, we have (interpreted through the ribbon functor $\mcal F$)

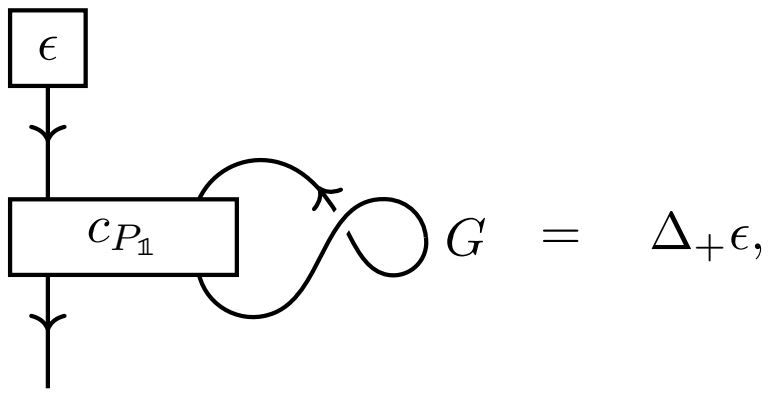
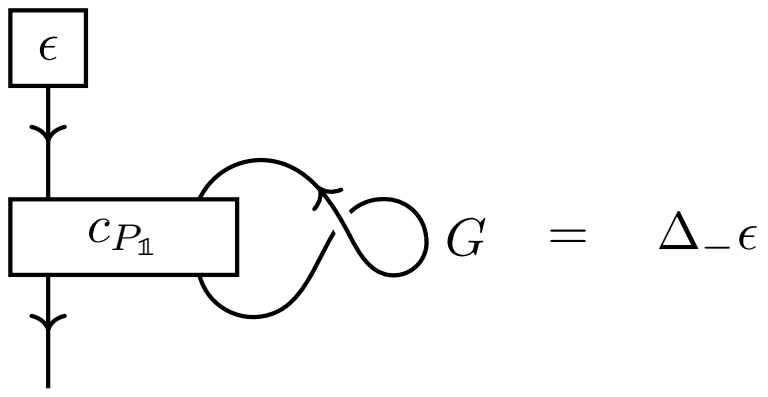
\begin{center}
    \begin{tikzpicture}
        \draw[thick,mid arrow=>] (-0.5,1) -- (-0.5,0.25);
        \draw[thick,mid arrow=>] (-0.5,-0.25) -- (-0.5,-1);
        \draw[thick,mid arrow=<] (0.5,0.25) to[curve through = {(1,0.5) (1.5,0) (1.75,-0.25)}] (2,0);
        \draw[white,line width = 1ex] (2,0) to[curve through = {(1.75,0.25) (1,-0.5)}] (0.5,-0.25);
        \draw[thick] (2,0) to[curve through = {(1.75,0.25) (1,-0.5)}] (0.5,-0.25);
        \node[right] at (2,0) {$G$};

        \node[draw,thick,minimum width = 1.5cm, minimum height = 0.5cm] at (0,0) {$c_{P_{\mbb 1}}$};
        \node[draw,thick,minimum width = 0.5cm, minimum height = 0.5cm] at (-0.5,1.25) {$\epsilon$};

        \node at (3.5,0) {$= \quad\Delta_+ \epsilon$,};

        \begin{scope}[shift = {(7,0)}]
            \draw[thick,mid arrow=>] (-0.5,1) -- (-0.5,0.25);
        \draw[thick,mid arrow=>] (-0.5,-0.25) -- (-0.5,-1);
        \draw[thick] (2,0) to[curve through = {(1.75,0.25) (1,-0.5)}] (0.5,-0.25);
        \draw[white,line width = 1ex] (0.5,0.25) to[curve through = {(1,0.5) (1.5,0) (1.75,-0.25)}] (2,0);
        \draw[thick,mid arrow=<] (0.5,0.25) to[curve through = {(1,0.5) (1.5,0) (1.75,-0.25)}] (2,0);
        \node[right] at (2,0) {$G$};

        \node[draw,thick,minimum width = 1.5cm, minimum height = 0.5cm] at (0,0) {$c_{P_{\mbb 1}}$};
        \node[draw,thick,minimum width = 0.5cm, minimum height = 0.5cm] at (-0.5,1.25) {$\epsilon$};
        \node at (3.5,0) {$= \quad\Delta_- \epsilon$};
        \end{scope}
    \end{tikzpicture}
\end{center}

\begin{lem}\label{lem:Deltas_are_1}
    When $\mcal D = \mcal Z(\mcal C)$, and the chromatic morphism $c^{\mcal Z}$ is chosen as above, then $\Delta_+ = \Delta_- = 1$.
\end{lem}

\begin{proof}
    We compute $\Delta_+$. The proof for $\Delta_-$ holds analogously. Set $H = UL(G)$. Let $\sum_i x_i x^i = \Lambda_{G^*\otimes H^*}$ and $\sum_i y_i y^i = \Lambda_G$.

    \begin{align}
        \Delta_+U(\epsilon) &= (U(\epsilon) \otimes \rev_{H})(\Id_{U(P_{\mbb 1})}\otimes \sigma^L_{H^*}) (U(c^{\mcal Z}_{P_{\mbb 1}})\otimes \Id_{H^*})(\Id_{U(P_{\mbb 1})}\otimes \rcoev_H)\\
        &= \sum_i (U(\epsilon) \otimes \rev_H)\\
        &\qquad \circ(\Id_{U(P_{\mbb 1})}\otimes \lev_G \otimes \Id_{H^*\otimes H})(c_{U(P_{\mbb 1})}\otimes x_i x^i \otimes \Id_H)(\Id_{U(P_{\mbb 1})}\otimes \rcoev_G \otimes \sigma^L_{H^*}) \nonumber\\
        &\qquad \circ(U(c^{\mcal Z}_{P_{\mbb 1}})\otimes \Id_{H^*})(\Id_{U(P_{\mbb 1})}\otimes \rcoev_H) \nonumber\\
        &= \sum_i (U(\epsilon) \otimes \rev_H)\\
        &\qquad \circ(\Id_{U(P_{\mbb 1})}\otimes \lev_G \otimes \Id_{H^*\otimes H})(c_{U(P_{\mbb 1})}\otimes x_i \otimes \Id_H)(\Id_{U(P_{\mbb 1})} \otimes \sigma^L_G \otimes x^i) \nonumber\\
        &\qquad \circ (\Id_{U(P_{\mbb 1})\otimes H} \otimes \rcoev_G \otimes \Id_{H^*})(U(c^{\mcal Z}_{P_{\mbb 1}})\otimes \Id_{H^*})(\Id_{U(P_{\mbb 1})}\otimes \rcoev_H) \nonumber\\
        &= \sum_i (U(\epsilon)\otimes \lev_G)(c_{U(P_{\mbb 1})}\otimes \rev_H \otimes \Id_{G^*})\\
        &\qquad \circ (\Id_{U(P_{\mbb 1})}\otimes \sigma^L_G \otimes x_i)(U(c^{\mcal Z}_{P_{\mbb 1}}) \otimes x^i \otimes \Id_G)(\Id_{U(P_{\mbb 1})} \otimes \rcoev_H \otimes \lcoev_G)\nonumber\\
        &= \sum_i (U(\epsilon) \otimes \lev_G)(c_{U(P_{\mbb 1})}\otimes \Id_{G^*})\\
        &\qquad \circ (\Id_{U(P_{\mbb 1})}\otimes y_i \otimes \Id_{G^*})(\Id_{U(P_{\mbb 1})}\otimes \Id_{G^*} \otimes y^i)(\Id_{U(P_{\mbb 1})} \otimes \lcoev_G)\nonumber\\
        &= U(\epsilon)
    \end{align}

    where line 2 introduces a chromatic identity, line 3 uses the naturality of the half braiding, line 4 uses the spin identity for the cutting morphism, line 5 uses Theorem \ref{thm:chromatic_halfbraiding_identity}, and line 6 uses the chromatic identity.
\end{proof}

\begin{thm}\label{thm:ExCGPT_is_RT}
    Let $M$ be presented via surgery on $S^3$ along a framed link $L$. Then \[Z_{\mathrm{ExCGPT}}(S^3,L;\mcal C) = B_{\mcal Z(\mcal C)}(M)\]
\end{thm}

\begin{proof}
    Since Lemma \ref{lem:Deltas_are_1} shows $\Delta_+ = \Delta_- = 1$, then we just need to show $Z_{\mathrm{ExCGPT}}(S^3,L;\mcal C) = \mcal F'(\widetilde{L})$ as described above. To this end, we give a Heegaard-Link diagram for the link $L\subseteq S^3$.

    Embed $L$ on a Heegaard surface for $S^3$ in the following way. Immerse the link in $S^2$ and to each crossing, stabilize as in Figure \ref{fig:resolvingcrossings}. If the immersion has no crossings, introduce a pair of crossings via a Reidemeister 2 move, then stabilize. The result is a Heegaard-Link diagram for the pair $(S^3,L)$ where the blackboard framing of $L$ matches the endowed framing.

    \begin{figure}
        \centering
        \begin{tikzpicture}
            \draw[darkgreen,ultra thick] (-1,1) -- (1,-1);
            \draw[darkgreen,ultra thick] (1,1) -- (0.5,0.5);
            \draw[darkgreen,ultra thick] (-0.5,-0.5) -- (-1,-1);
            \draw[blue,ultra thick] (0.5,0.5) -- (-0.5,-0.5);
            \draw[draw=black, fill=lightgray] (0.5,0.5) circle (0.25);
            \draw[draw=black, fill=lightgray] (-0.5,-0.5) circle (0.25);
            \draw[draw=red,ultra thick] (-0.5,-0.5) circle (0.25);
        \end{tikzpicture}
        \caption{Stabilizing bridge at a crossing}
        \label{fig:resolvingcrossings}
    \end{figure}
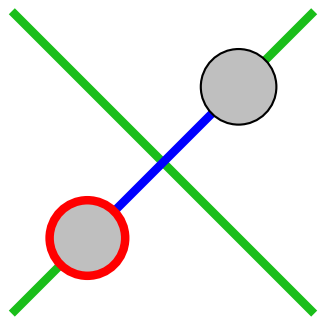

    The ribbon graph about a crossing is given by:

    
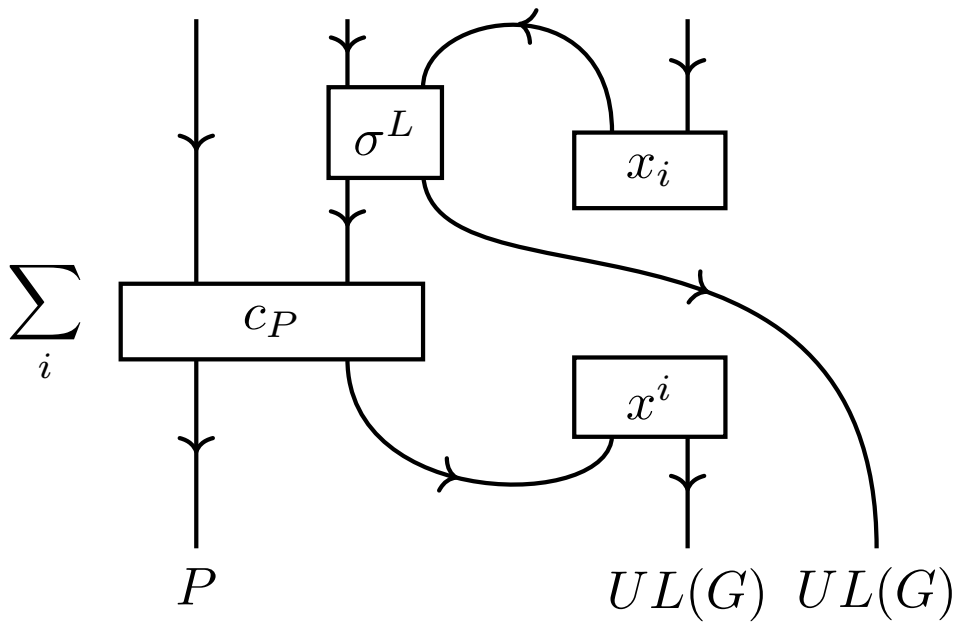
\begin{center}
        \begin{tikzpicture}
            \node at (-1.5,0) {$\displaystyle\sum_i$};
            \draw[mid arrow=>,thick] (-0.5,2) -- (-0.5,0.25);
            \draw[mid arrow=>,thick] (-0.5,-0.25) -- (-0.5,-1.5)node[below]{$P$};
            \draw[mid arrow=>,thick] (2.75,2) -- (2.75,1.25);
            \draw[mid arrow=>,thick] (2.75,-0.75) -- (2.75,-1.5)node[below]{$UL(G)$};
            \draw[mid arrow=>,thick] (1,1) .. controls +(0,-1) and +(0,2.5) .. (4,-1.5);
            \node[below] at (4,-1.5) {$UL(G)$};
            \draw[mid arrow=<,thick] (1,1.55) .. controls +(0,0.5) and +(0,1) .. (2.25,1.25);
            \draw[mid arrow=>,thick] (0.5,-0.25) .. controls +(0,-1) and +(0,-0.5) .. (2.25,-0.75);
            \draw[mid arrow=>,thick] (0.5,2) -- (0.5,1.55);
            \draw[mid arrow=>,thick] (0.5,1) -- (0.5,0.25);

            \node[draw,thick,fill = white,minimum width = 0.75cm, minimum height = 0.6cm] at (0.75,1.25) {$\sigma^L$};
            \node[draw,thick,fill = white,minimum width = 2cm, minimum height = 0.5cm] at (0,0) {$c_P$};
            \node[draw,thick,fill = white,minimum width = 1cm, minimum height = 0.5cm] at (2.5,1) {$x_i$};
            \node[draw,thick,fill = white,minimum width = 1cm, minimum height = 0.5cm] at (2.5,-0.5) {$x^i$};
        \end{tikzpicture}
    \end{center}

    In this example, using the naturality of the half braiding and the chromatic identity yields the braiding of $UL(G)$ with itself or its dual, depending on the link. If the other strand crossed over we would obtain the inverse braiding of $UL(G)$ with itself or its dual.

    To see that the invariants are equal, we use that $Z_{\mathrm{ExCGPT}}$ is defined by taking the resulting ribbon graph in $S^2$, viewing it as the braid closure of some projective object, and taking the trace $\ctr$. We use the embedded ribbon graph $O_{L(G)}$ to do this.

    Using the cyclic property of $\ctr$ we can ensure that the invariant is of the form \[Z_{\mathrm{ExCGPT}}(S^3,L;\mcal C) = \ctr_G(\pi_{G,\mbb 1}U(\Omega_G)U(\ell)U(h)i_{G,\mbb 1})\] where $\ell$ is the morphism in $\mcal Z(\mcal C)$ defined by coloring the link using $c^{\mcal Z}$. Then by the definition of $\ztr$, we have \[Z_{\mathrm{ExCGPT}}(S^3,L;\mcal C) = \ztr_{L(G)}(\ell\circ h)\]

    To obtain $\mcal F'(\widetilde{L})$, we use that the projective cover $P_{\mbb 1}$ of $\mbb 1\in \mcal Z(\mcal C)$ may be inserted using that $L(G)$ is a projective generator.
\end{proof}

\section{Relationship with Surgery}
\label{Surgery_Sec}

Let $(M,L)$ be a closed connected oriented 3-manifold with an embedded framed link $L$. Then recall Theorem \ref{thm:surgerydiagram} presenting from a Heegaard-Link diagram for $(M,L)$ a Heegaard diagram for the surgery manifold $M(L)$, where on each component, we perform:

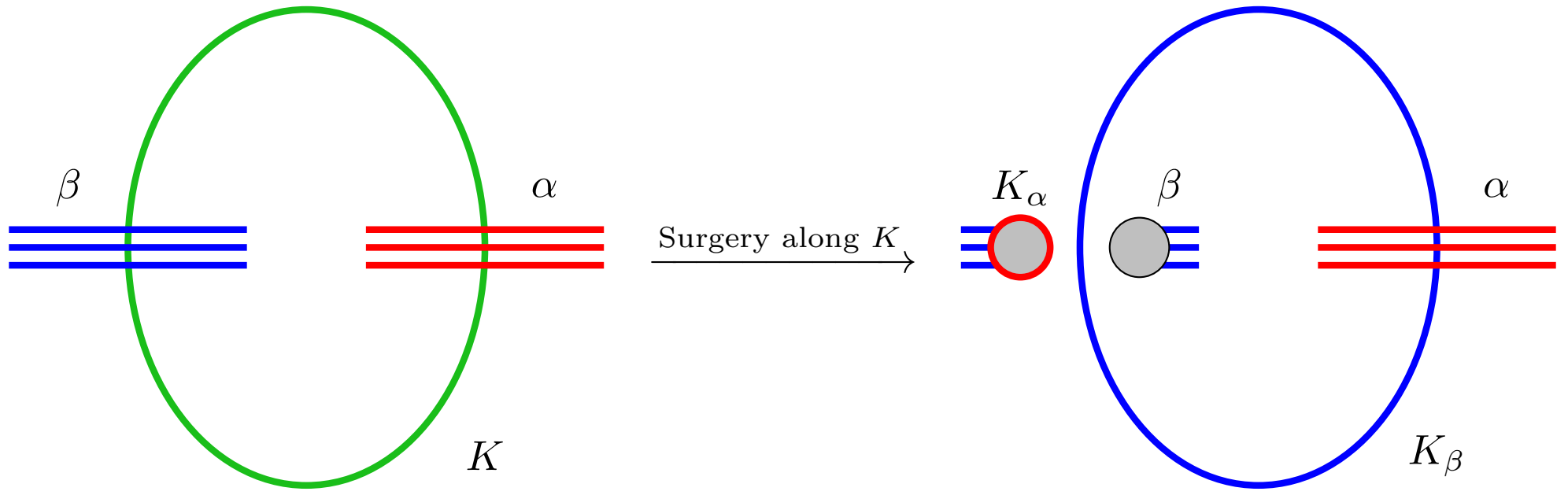
\begin{center}
    \begin{tikzpicture}
        \draw[darkgreen, ultra thick] (0,0) ellipse (1.5cm and 2cm);
        \draw[blue, ultra thick] (-0.5,0) -- (-2.5,0);
        \draw[blue, ultra thick] (-0.5,0.15) -- (-2.5,0.15);
        \draw[blue, ultra thick] (-0.5,-0.15) -- (-2.5,-0.15);
        \draw[red, ultra thick] (0.5,0) -- (2.5,0);
        \draw[red, ultra thick] (0.5,0.15) -- (2.5,0.15);
        \draw[red, ultra thick] (0.5,-0.15) -- (2.5,-0.15);
        \node at (1.5,-1.75) {$K$};
        \node at (-2,0.5) {$\beta$};
        \node at (2,0.5) {$\alpha$};
        \node at (4,0) {$\xrightarrow{\text{Surgery along }K}$};
    \begin{scope}[shift = {(8,0)}]
        \draw[blue, ultra thick] (0,0) ellipse (1.5cm and 2cm);
        \draw[red, ultra thick] (0.5,0) -- (2.5,0);
        \draw[red, ultra thick] (0.5,0.15) -- (2.5,0.15);
        \draw[red, ultra thick] (0.5,-0.15) -- (2.5,-0.15);
        \draw[blue, ultra thick] (-0.5,0) -- (-1,0);
        \draw[blue, ultra thick] (-0.5,0.15) -- (-1,0.15);
        \draw[blue, ultra thick] (-0.5,-0.15) -- (-1,-0.15);
        \draw[blue, ultra thick] (-2,0) -- (-2.5,0);
        \draw[blue, ultra thick] (-2,0.15) -- (-2.5,0.15);
        \draw[blue, ultra thick] (-2,-0.15) -- (-2.5,-0.15);
        \draw[draw=black, fill=lightgray] (-1,0) circle (0.25);
        \draw[draw=black, fill=lightgray] (-2,0) circle (0.25);
        \draw[red, ultra thick] (-2,0) circle (0.25);
        \node at (1.5,-1.75) {$K_{\beta}$};
        \node at (2,0.5) {$\alpha$};
        \node at (-0.75,0.5) {$\beta$};
        \node at (-2,0.5) {$K_\alpha$};
        \end{scope}
    \end{tikzpicture}
\end{center}

which is the rotated version of Figure \ref{fig:SurgeryHeegaardDiagram}.

\begin{thm}\label{thm:surgery_invariant}
    For a pair $(M,L)$, \[Z_{\mathrm{ExCGPT}}(M,L;\mcal C) = Z_{\mathrm{ExCGPT}}(M(L),\emptyset;\mcal C)\]
\end{thm}

\begin{proof}
    The left-hand side of the formula above yields the ribbon graph locally

    
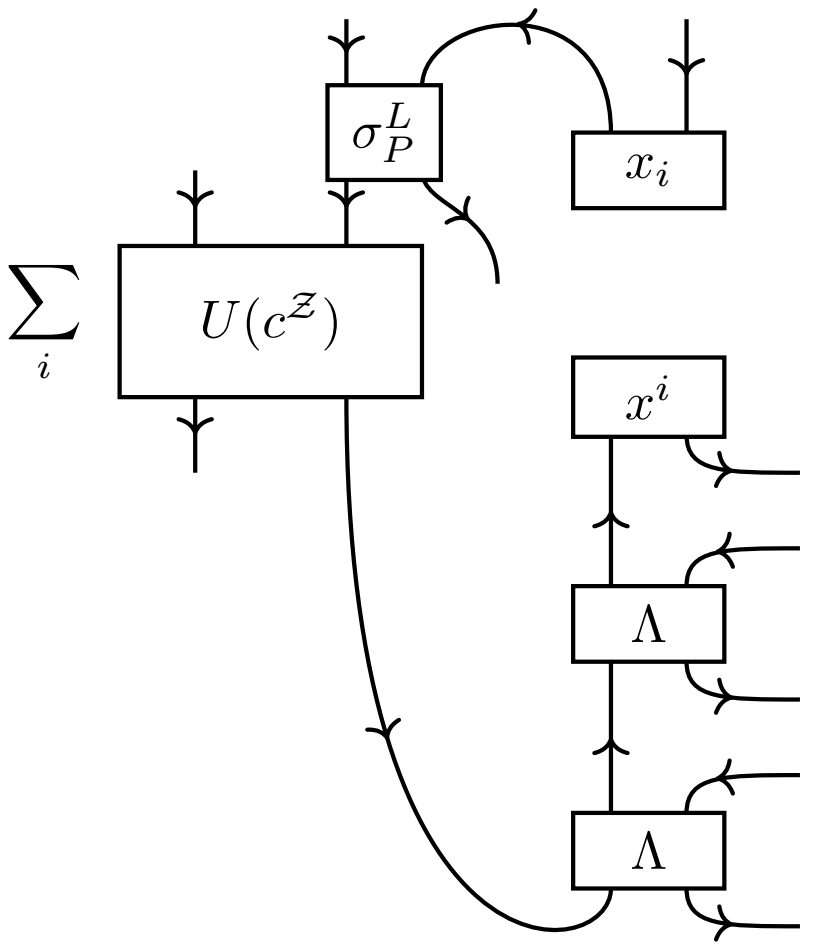
\begin{center}
        \begin{tikzpicture}
            \node at (-1.5,0) {$\displaystyle\sum_i$};
            \draw[mid arrow=>,thick] (-0.5,1) -- (-0.5,0.5);
            \draw[mid arrow=>,thick] (-0.5,-0.5) -- (-0.5,-1);
            \draw[mid arrow=>,thick] (2.75,2) -- (2.75,1.25);
            
            \draw[mid arrow=>,thick] (1,1) .. controls +(0,-0.25) and +(0,0.5) .. (1.5,0.25);
            \draw[mid arrow=<,thick] (1,1.55) .. controls +(0,0.5) and +(0,1) .. (2.25,1.25);
            \draw[mid arrow=>,thick] (0.5,-0.5) .. controls +(0,-4) and +(0,-0.5) .. (2.25,-3.75);
            \draw[mid arrow=>,thick] (0.5,2) -- (0.5,1.55);
            \draw[mid arrow=>,thick] (0.5,1) -- (0.5,0.5);

            \draw[mid arrow=>,thick] (2.75,-0.75) .. controls +(0,-0.25) and +(-0.5,0) .. (3.5,-1);
            \draw[mid arrow=<,thick] (2.75,-1.75) .. controls +(0,0.25) and +(-0.5,0) .. (3.5,-1.5);

            \draw[mid arrow=>,thick] (2.75,-2.25) .. controls +(0,-0.25) and +(-0.5,0) .. (3.5,-2.5);
            \draw[mid arrow=<,thick] (2.75,-3.25) .. controls +(0,0.25) and +(-0.5,0) .. (3.5,-3);

            \draw[mid arrow=>,thick] (2.75,-3.75) .. controls +(0,-0.25) and +(-0.5,0) .. (3.5,-4);

            \draw[mid arrow=>,thick] (2.25,-3.25) -- (2.25,-2.25);
            \draw[mid arrow=>,thick] (2.25,-1.75) -- (2.25,-0.75);

            \node[draw,thick,fill = white,minimum width = 0.75cm, minimum height = 0.6cm] at (0.75,1.25) {$\sigma^L_P$};
            \node[draw,thick,fill = white,minimum width = 2cm, minimum height = 1cm] at (0,0) {$U(c^{\mcal Z})$};
            \node[draw,thick,fill = white,minimum width = 1cm, minimum height = 0.5cm] at (2.5,1) {$x_i$};
            \node[draw,thick,fill = white,minimum width = 1cm, minimum height = 0.5cm] at (2.5,-0.5) {$x^i$};

            \node[draw,thick,fill = white,minimum width = 1cm, minimum height = 0.5cm] at (2.5,-2) {$\Lambda$};

            \node[draw,thick,fill = white,minimum width = 1cm, minimum height = 0.5cm] at (2.5,-3.5) {$\Lambda$};
        \end{tikzpicture}   
    \end{center}

    where $P$ is a projective object defined by the $\beta$ curves.

    Then, using naturality of $\Lambda$, we can move $x^i$ past each $\Lambda$ cutting morphism. Additionally, we can perform an isotopy of the graph to bring $x^i$ and $x_i$ to the left so that locally in a disk about the intersections with $L$ and the $\beta$ curves, we have 

    
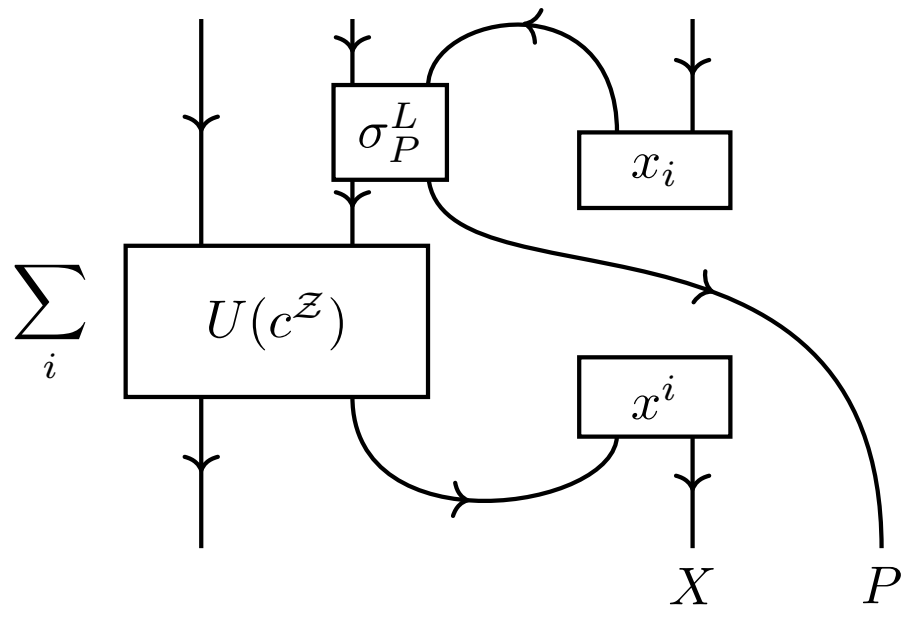
\begin{center}
        \begin{tikzpicture}
            \node at (-1.5,0) {$\displaystyle\sum_i$};
            \draw[mid arrow=>,thick] (-0.5,2) -- (-0.5,0.5);
            \draw[mid arrow=>,thick] (-0.5,-0.5) -- (-0.5,-1.5);
            \draw[mid arrow=>,thick] (2.75,2) -- (2.75,1.25);
            \draw[mid arrow=>,thick] (2.75,-0.75) -- (2.75,-1.5)node[below]{$X$};
            \draw[mid arrow=>,thick] (1,1) .. controls +(0,-1) and +(0,2.5) .. (4,-1.5);
            \node[below] at (4,-1.5) {$P$};
            \draw[mid arrow=<,thick] (1,1.55) .. controls +(0,0.5) and +(0,1) .. (2.25,1.25);
            \draw[mid arrow=>,thick] (0.5,-0.5) .. controls +(0,-1) and +(0,-0.5) .. (2.25,-0.75);
            \draw[mid arrow=>,thick] (0.5,2) -- (0.5,1.55);
            \draw[mid arrow=>,thick] (0.5,1) -- (0.5,0.5);

            \node[draw,thick,fill = white,minimum width = 0.75cm, minimum height = 0.6cm] at (0.75,1.25) {$\sigma^L_P$};
            \node[draw,thick,fill = white,minimum width = 2cm, minimum height = 1cm] at (0,0) {$U(c^{\mcal Z})$};
            \node[draw,thick,fill = white,minimum width = 1cm, minimum height = 0.5cm] at (2.5,1) {$x_i$};
            \node[draw,thick,fill = white,minimum width = 1cm, minimum height = 0.5cm] at (2.5,-0.5) {$x^i$};
        \end{tikzpicture}
    \end{center}

    where $X$ connects to the bottom-most $\Lambda$ in the previous diagram.

    Then using Theorem \ref{thm:chromatic_halfbraiding_identity}, we have that locally this is skein equivalent to 

    
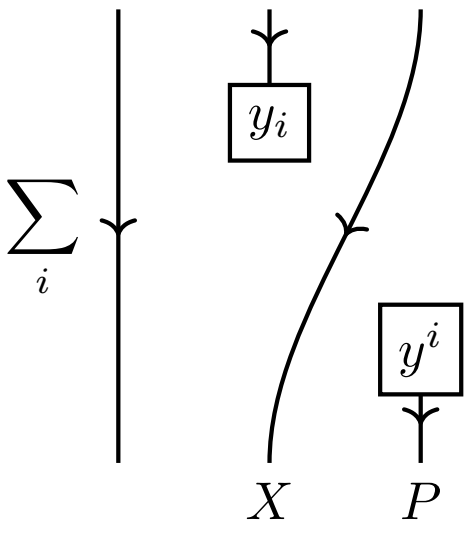
\begin{center}
        \begin{tikzpicture}
            \node at (0,0) {$\displaystyle\sum_i$};
                \draw[mid arrow=>,thick] (0.5,1.5) -- (0.5,-1.5);
                
                \draw[mid arrow=>,thick] (1.5,1.5) -- (1.5,1);
                \draw[mid arrow=>,thick] (2.5,-1) -- (2.5,-1.5)node[below]{$P$};
                \node[draw,thick,fill=white,minimum width = 0.5cm, minimum height = 0.5cm] at (1.5,0.75) {$y_i$};
                \node[draw,thick,fill=white,minimum width = 0.5cm, minimum height = 0.5cm] at (2.5,-0.75) {$y^i$};

                \draw[mid arrow=>, thick] (2.5,1.5) .. controls +(0,-1) and +(0,1) .. (1.5,-1.5)node[below] {$X$};
        \end{tikzpicture}
    \end{center}

    Putting this back into the first diagram yields:

    
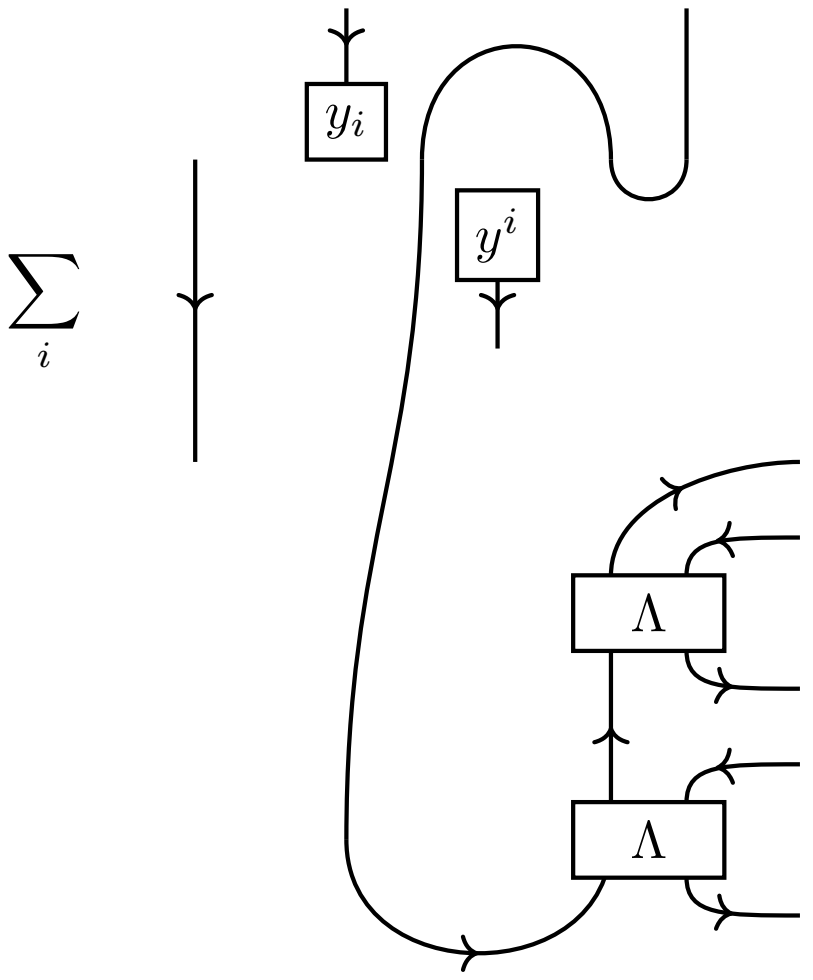
\begin{center}
        \begin{tikzpicture}
            \node at (-1.5,0) {$\displaystyle\sum_i$};
            \draw[mid arrow=>,thick] (-0.5,1) -- (-0.5,-1);
            \draw[thick] (2.75,2) -- (2.75,1);
            \draw[thick] (2.75,1) .. controls +(0,-0.35) and +(0,-0.35) .. (2.25,1);
            
            \draw[mid arrow=>,thick] (1.5,0.25) -- (1.5,-0.25);
            \draw[thick] (2.25,1) .. controls +(0,1) and +(0,1) .. (1,1);
            \draw[thick] (1,1) .. controls +(0,-2) and +(0,2) .. (0.5,-3.5);
            \draw[mid arrow=>,thick] (0.5,-3.5) .. controls +(0,-1) and +(0,-1) .. (2.25,-3.5);
            \draw[mid arrow=>,thick] (0.5,2) -- (0.5,1.5);

            \draw[mid arrow=>,thick] (2.25,-1.75) .. controls +(0,0.5) and +(-0.5,0) .. (3.5,-1);
            \draw[mid arrow=<,thick] (2.75,-1.75) .. controls +(0,0.25) and +(-0.5,0) .. (3.5,-1.5);

            \draw[mid arrow=>,thick] (2.75,-2.25) .. controls +(0,-0.25) and +(-0.5,0) .. (3.5,-2.5);
            \draw[mid arrow=<,thick] (2.75,-3.25) .. controls +(0,0.25) and +(-0.5,0) .. (3.5,-3);

            \draw[mid arrow=>,thick] (2.75,-3.75) .. controls +(0,-0.25) and +(-0.5,0) .. (3.5,-4);

            \draw[mid arrow=>,thick] (2.25,-3.25) -- (2.25,-2.25);

            \node[draw,thick,fill=white,minimum width = 0.5cm, minimum height = 0.5cm] at (0.5,1.25) {$y_i$};
            \node[draw,thick,fill=white,minimum width = 0.5cm, minimum height = 0.5cm] at (1.5,0.5) {$y^i$};

            \node[draw,thick,fill = white,minimum width = 1cm, minimum height = 0.5cm] at (2.5,-2) {$\Lambda$};

            \node[draw,thick,fill = white,minimum width = 1cm, minimum height = 0.5cm] at (2.5,-3.5) {$\Lambda$};
        \end{tikzpicture}   
    \end{center}

    The right-hand side is given by 

    
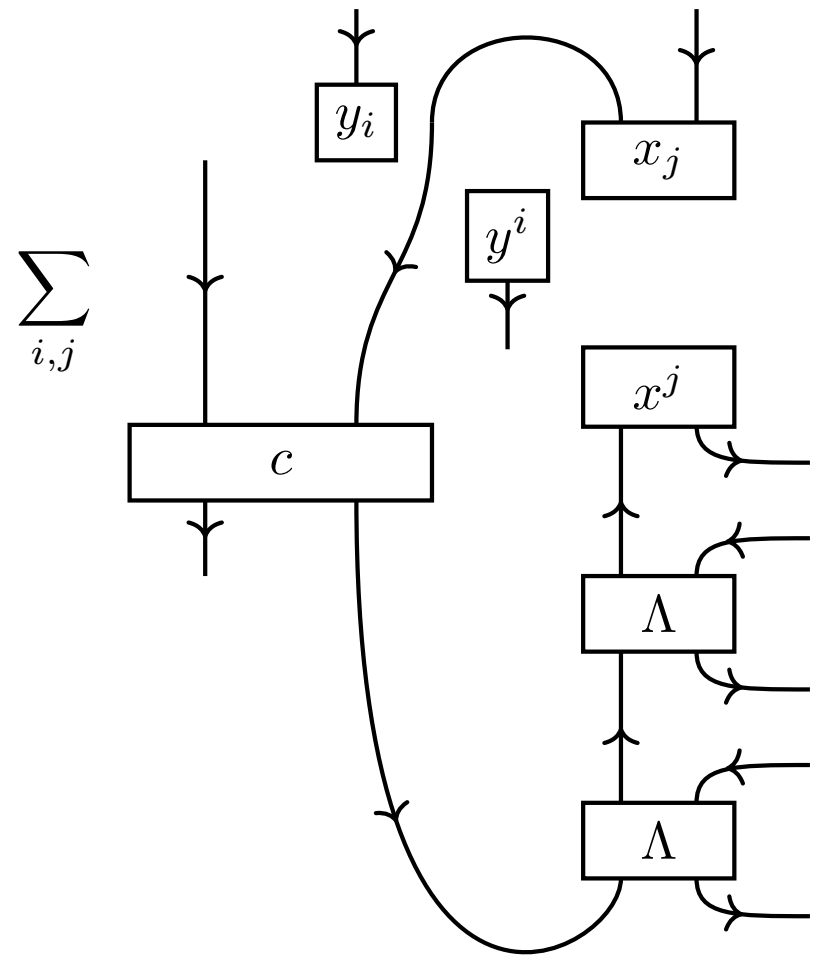
\begin{center}
        \begin{tikzpicture}
            \node at (-1.5,0) {$\displaystyle\sum_{i,j}$};
            \draw[mid arrow=>,thick] (-0.5,1) -- (-0.5,-0.75);
            \draw[mid arrow=>,thick] (-0.5,-1.25) -- (-0.5,-1.75);
            \draw[mid arrow=>,thick] (2.75,2) -- (2.75,1.25);
            
            \draw[mid arrow=>,thick] (1.5,0.25) -- (1.5,-0.25);
            \draw[mid arrow=>,thick] (0.5,-1.25) .. controls +(0,-4) and +(0,-0.5) .. (2.25,-3.75);
            \draw[mid arrow=>,thick] (0.5,2) -- (0.5,1.5);

            \draw[thick] (2.25,1.25) .. controls +(0,0.75) and +(0,0.75) .. (1,1.25);
            \draw[mid arrow=>,thick] (1,1.25) .. controls +(0,-1) and +(0,1) .. (0.5,-0.75);

            \draw[mid arrow=>,thick] (2.75,-0.75) .. controls +(0,-0.25) and +(-0.5,0) .. (3.5,-1);
            \draw[mid arrow=<,thick] (2.75,-1.75) .. controls +(0,0.25) and +(-0.5,0) .. (3.5,-1.5);

            \draw[mid arrow=>,thick] (2.75,-2.25) .. controls +(0,-0.25) and +(-0.5,0) .. (3.5,-2.5);
            \draw[mid arrow=<,thick] (2.75,-3.25) .. controls +(0,0.25) and +(-0.5,0) .. (3.5,-3);

            \draw[mid arrow=>,thick] (2.75,-3.75) .. controls +(0,-0.25) and +(-0.5,0) .. (3.5,-4);

            \draw[mid arrow=>,thick] (2.25,-3.25) -- (2.25,-2.25);
            \draw[mid arrow=>,thick] (2.25,-1.75) -- (2.25,-0.75);

            \node[draw,thick,fill = white,minimum width = 2cm, minimum height = 0.5cm] at (0,-1) {$c$};
            \node[draw,thick,fill = white,minimum width = 1cm, minimum height = 0.5cm] at (2.5,1) {$x_j$};
            \node[draw,thick,fill = white,minimum width = 1cm, minimum height = 0.5cm] at (2.5,-0.5) {$x^j$};

            \node[draw,thick,fill=white,minimum width = 0.5cm, minimum height = 0.5cm] at (0.5,1.25) {$y_i$};
            \node[draw,thick,fill=white,minimum width = 0.5cm, minimum height = 0.5cm] at (1.5,0.5) {$y^i$};

            \node[draw,thick,fill = white,minimum width = 1cm, minimum height = 0.5cm] at (2.5,-2) {$\Lambda$};

            \node[draw,thick,fill = white,minimum width = 1cm, minimum height = 0.5cm] at (2.5,-3.5) {$\Lambda$};
        \end{tikzpicture}   
    \end{center}

    Then performing a similar procedure by moving $x^j$ past each $\Lambda$ and $x_j$ to the left yields the chromatic identity, after which we obtain the diagram

    
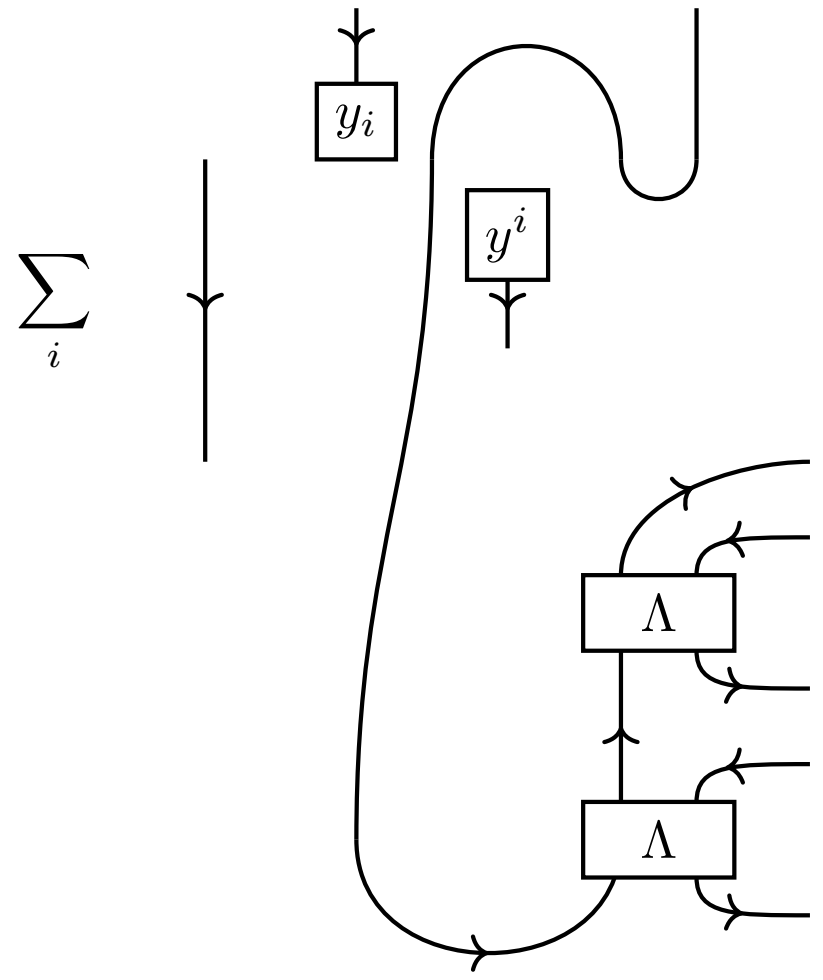
\begin{center}
        \begin{tikzpicture}
            \node at (-1.5,0) {$\displaystyle\sum_i$};
            \draw[mid arrow=>,thick] (-0.5,1) -- (-0.5,-1);
            \draw[thick] (2.75,2) -- (2.75,1);
            \draw[thick] (2.75,1) .. controls +(0,-0.35) and +(0,-0.35) .. (2.25,1);
            
            \draw[mid arrow=>,thick] (1.5,0.25) -- (1.5,-0.25);
            \draw[thick] (2.25,1) .. controls +(0,1) and +(0,1) .. (1,1);
            \draw[thick] (1,1) .. controls +(0,-2) and +(0,2) .. (0.5,-3.5);
            \draw[mid arrow=>,thick] (0.5,-3.5) .. controls +(0,-1) and +(0,-1) .. (2.25,-3.5);
            \draw[mid arrow=>,thick] (0.5,2) -- (0.5,1.5);

            \draw[mid arrow=>,thick] (2.25,-1.75) .. controls +(0,0.5) and +(-0.5,0) .. (3.5,-1);
            \draw[mid arrow=<,thick] (2.75,-1.75) .. controls +(0,0.25) and +(-0.5,0) .. (3.5,-1.5);

            \draw[mid arrow=>,thick] (2.75,-2.25) .. controls +(0,-0.25) and +(-0.5,0) .. (3.5,-2.5);
            \draw[mid arrow=<,thick] (2.75,-3.25) .. controls +(0,0.25) and +(-0.5,0) .. (3.5,-3);

            \draw[mid arrow=>,thick] (2.75,-3.75) .. controls +(0,-0.25) and +(-0.5,0) .. (3.5,-4);

            \draw[mid arrow=>,thick] (2.25,-3.25) -- (2.25,-2.25);

            \node[draw,thick,fill=white,minimum width = 0.5cm, minimum height = 0.5cm] at (0.5,1.25) {$y_i$};
            \node[draw,thick,fill=white,minimum width = 0.5cm, minimum height = 0.5cm] at (1.5,0.5) {$y^i$};

            \node[draw,thick,fill = white,minimum width = 1cm, minimum height = 0.5cm] at (2.5,-2) {$\Lambda$};

            \node[draw,thick,fill = white,minimum width = 1cm, minimum height = 0.5cm] at (2.5,-3.5) {$\Lambda$};
        \end{tikzpicture}   
    \end{center}

    as desired.
\end{proof}

\begin{cor}\label{cor:CGPT_equals_RTcenter}
    For a chromatic category $\mcal C$ and a closed connected oriented 3-manifold $M$,
    \[Z_{\mathrm{CGPT}}(M;\mcal C) = B_{\mcal Z(\mcal C)}(M)\]
\end{cor}

\section{Obtaining Invariants for Ribbon Chromatic Categories}
\label{RibbonInvariants_Sec}

Given a ribbon chromatic category $\mcal D$, we can use the natural embedding into the Drinfeld center $\mcal Z(\mcal D)$ to obtain $B_{\mcal D}(M)$ using Heegaard-Link diagrams. 

Let $B:\mcal D\to \mcal Z(\mcal D)$ be the functor defined by sending an object $X$ to $(X,c_{X,\bullet})$, where $c_{X,Y}:X\otimes Y\to Y\otimes X$ is the braiding in $\mcal D$. Then $B$ defines a section of the forgetful functor $U:\mcal Z(\mcal D)\to \mcal D$.

Given a Heegaard-Link diagram and a chromatic morphism $c$ in $\mcal D$, use $B(c)$ to color the link with $B(G)$, where $G$ is a projective generator of $\mcal D$, now treated as an object in $\mcal Z(\mcal D)$ using the braiding. Then since $B$ is a section, $UB(G) = G$ and $UB(c) = c$. Finally, we use the ribbon graph $O_G = \mathrm{tr}_G(h)$ for some $h:G\to G$ with $\mathrm{t}^{\mcal D}_G(h) = 1$ to define an invariant of Heegaard-Link diagrams $Z_{\mathrm{ExRib}}(M,L;\mcal D)$ analogously to $Z_{\mathrm{ExCGPT}}(M,L;\mcal C)$.

\begin{thm}\label{thm:ExRib_is_RT}
    Let $M$ be presented via surgery on $S^3$ along a framed link $L$. Then \[\frac{Z_{\mathrm{ExRib}}(S^3,L;\mcal D)}{\Delta_+^r\Delta_-^s} = B_{\mcal D}(M)\] where $(r,s)$ is the signature of the link $L$.
\end{thm}

\begin{proof}
    Using a similar argument as in Theorem \ref{thm:ExCGPT_is_RT}, we see that at a ``bridge'' we have the braiding (or inverse braiding) of the color of the link with itself (or its dual). The result is a link diagram colored by $G$ using the chromatic morphism $c$ in $\mcal D$. This is exactly $\mcal F'(\widetilde{L})$.
\end{proof}

\section{Discussion and Examples}
\label{Discussion_Xmps_Sec}

\subsection{On the relationship between the topological quantum field theories}

The authors expect that this relationship generalizes to a correspondence between the $(3+1)$ Skein topological quantum field theory of \cite{Geer2026skein3plus1} and the non-compact $(2+1)$ Skein topological quantum field theory of \cite{Geer2023noncompact} where the latter is a boundary theory for the former based on the Drinfeld center. We leave the exploration of this relationship to future work.

\subsection{An example from Hopf algebras}

We conclude this paper with an example coming from Hopf algebras.

\begin{xmp}
    Let $H$ be a finite-dimensional unimodular spherical Hopf algebra with pivot $g$. The non-degenerate m-trace on $\mathrm{Rep}(H)$ is given uniquely by its evaluation on endomorphisms of $H$. Given a symmetrized integral $\mu:H\to \mbb k$, the m-trace is given by \[\mu(f(1)), \text{ for any }f:H\to H\]

    For the Drinfeld double $D(H)$ of $H$, the symmetrized integral is given by $\Lambda \otimes \mu$, where $\Lambda$ is the two-sided integral in $H$. Thus the non-degenerate m-trace on $\mcal Z(\mathrm{Rep}(H)) \cong \mathrm{Rep}(D(H))$ is given by \[(\Lambda\otimes \mu)(f(\epsilon\otimes 1)),\text{ for any }f:D(H)\to D(H)\] (see \cite{beliakova2021modified}).

    We verify the m-trace on $\mathrm{Rep}(D(H))$ using the definition in Section \ref{Center_Sec}.

    Recall that for a Hopf algebra $H$, the category of finite-dimensional left-left Yetter-Drinfeld modules ${}^H_H\mcal{YD}$ is equivalent to $\mathrm{Rep}(D(H))$. These modules consist of $H$-modules $M$ together with a left $H$-coaction $\Delta_M:M \to H\otimes M$ which we denote by $\Delta_M(m) = m_{(-1)}\otimes m_{(0)}$ satisfying
    \[(h.m)_{(-1)}\otimes (h.m)_{(0)} = h_{(1)}m_{(-1)}S(h_{(3)})\otimes h_{(2)}.m_{(0)}\]

    For a Hopf algebra $H$, the functors $R$ and $L$ are defined using the adjoint representation $H_{\mathrm{ad}}$ and the coadjoint representation $H^*_{\mathrm{coad}}$, respectively. Specifically, for an $H$-representation $V$,
    \[R(V) = H_{\mathrm{ad}}\otimes_{\Bbbk} V ,\qquad\qquad L(V) = H^*_{\mathrm{coad}}\otimes_{\Bbbk} V\]

    The counit for $R$ is $\pi_{V,\mbb 1}(h\otimes v) = \epsilon(h) v$ and the unit for $L$ is $i_{V,\mbb 1}(v) = \epsilon\otimes v$

    The unimodular isomorphism $\Omega_V:L(V)\to R(V)$ is given by\[\Omega_V(f\otimes v) = \Lambda_{(1)}f(\Lambda_{(2)}) \otimes v\] where $\Lambda$ is the two-sided integral in $H$.

    Then we can compute. For an endomorphism $f:L(H)\to L(H)$, we have:
    \begin{align*}
        \pi_{H,\mbb 1} U(\Omega_H)U(f)i_{H,\mbb 1} &= \pi_{H,\mbb 1}\Omega_H(f(\epsilon\otimes ?))\\
        &= \pi_{H,\mbb 1}(\Lambda_{(1)} \otimes (\Lambda_{(2)}\otimes \Id_H)(f(\epsilon\otimes ?)))\\
        &= (\Lambda \otimes \Id_H)(f(\epsilon\otimes ?))
    \end{align*}

    Then the non-degenerate m-trace on $\mathrm{Rep}(H)$ yields
    \[\t^{\mathrm{Rep}(H)}_H(\pi_{H,\mbb 1} U(\Omega_H)U(f)i_{H,\mbb 1}) = (\Lambda\otimes \mu)(f(\epsilon\otimes 1))\]
    as desired.
\end{xmp}

Using the above example, we expect that when $H$ is a semisimple involutory Hopf algebra, this construction recovers the invariant $Z_{\mathrm{ExKup}}$ defined in \cite{bridges20263man}. Additionally, when $H$ is unimodular and unibalanced, we expect to recover the renormalized Hennings invariant. We do not prove this here.

\vspace{0.5cm}
\noindent\textbf{Acknowledgements.}
N.B. and N.G. would like to thank Bertrand Patureau-Mirand for helpful discussions. 

\vspace{0.5cm}
\noindent\textbf{AI Disclosure.}
The statements of all theorems were proposed entirely by the authors with no outside help. Portions of the proof of Theorems \ref{thm:mtraceforcenter} and \ref{thm:chromaticforcenter} were initially drafted by the authors but then completed using OpenAI's ChatGPT 6 Astra. Additionally, the statement for defining the modified trace in the center using the end and coend was proposed by the authors with no outside help and OpenAI's ChatGPT 6 Astra helped complete the proof by proposing using the natural isomorphism $\nu:L\to R$ and using a module category argument to show $\epsilon'$ can be written using $\Omega$. Finally, the idea in the proof of Theorem \ref{thm:chromatic_halfbraiding_identity} to rewrite the cutting morphism using the cutting morphism in the center was generated by OpenAI's ChatGPT 6 Astra. The authors independently checked all arguments and rewrote the arguments to suit the flow of the paper. Additionally, OpenAI's ChatGPT 6 Astra and Anthropic's Claude Opus 5.5 were used to proofread the article for typos and small errors. The authors take full responsibility for the claims made and proofs of the arguments.

\bibliographystyle{alpha}
\newcommand{\etalchar}[1]{$^{#1}$}

\end{document}